\documentclass[11pt]{article}
\usepackage[colorlinks=true,linkcolor=blue,citecolor = blue]{hyperref}
\usepackage {verbatim}
\usepackage[commandnameprefix=always]{changes}
\usepackage{mathrsfs}
\usepackage[english]{babel}
\usepackage[nottoc]{tocbibind}
\usepackage{amssymb}
\usepackage{amsthm}
\usepackage{amsmath}
\usepackage{graphicx}
\usepackage{etoolbox,xstring,mfirstuc,textcase}
\usepackage{empheq}
\usepackage{indentfirst}
\usepackage{cite} 
\usepackage{cases}
\usepackage{graphics}
\usepackage{xcolor}  
\usepackage{amsmath,bm}
\usepackage{changes}
\usepackage{lineno}
\usepackage{cases}
\usepackage{hyperref}
\newtheorem{theorem}{\textbf{Theorem}}[section]
\newtheorem{lemma}{\textbf{Lemma}}[section]
\newtheorem{proposition}{\textbf{Proposition}}[section]
\newtheorem{corollary}{\textbf{Corollary}}[section]
\newtheorem{remark}{\textbf{Remark}}[section]
\newtheorem{definition}{\textbf{Definition}}[section]
\allowdisplaybreaks[4] 
\def\be{\begin{equation}}
	\def\ee{\end{equation}}
\def\bt{\begin{theorem}}
	\def\et{\end{theorem}}

\def\bl{\begin{lemma}}
	\def\el{\end{lemma}}
\def\br{\begin{remark}}
	\def\er{\end{remark}}
\def\bp{\begin{proposition}}
	\def\ep{\end{proposition}}
\def\bc{\begin{corollary}}
	\def\ec{\end{corollary}}
\def\bd{\begin{definition}}
	\def\ed{\end{definition}}

\begin{document}
	
\title{Global Well-posedness and Long-time Behavior of the Two-dimensional General Ericksen--Leslie System in the Isotropic Case under a Magnetic Field }
	
\author{
	{Qingtong Wu}
	\footnote{School of Mathematical Sciences, Fudan University, Handan Road 220, Shanghai 200433, China. Email:  
       \textit{qtwu570@gmail.com}}
        }

\date{}
	
\maketitle
	

\begin{abstract}
\noindent 
  This paper establishes the global well-posedness and long-time dynamics of the general Ericksen--Leslie system for isotropic nematic liquid crystals under a constant magnetic field.
 On the two-dimensional torus $\mathbb{T}^2$, a liquid crystal molecule coincides with itself under rotations by integer multiples of $\pi$, which results in special boundary conditions.
  We prove the existence of global-in-time strong solutions by developing novel high-order energy estimates and employing compactness techniques.
  A key challenge lies in controlling the orientation of the liquid crystal molecules. 
  After achieving a uniform bound for the molecular orientation angle in $\mathbb{S}^1$, we further characterize the long-time behavior of the solutions. 
  This is accomplished by applying the Lojasiewicz--Simon inequality, which reveals the convergence of the solutions as time approaches infinity.
 \medskip \\
\medskip\noindent	
\end{abstract}

\tableofcontents
	
\section{ Introduction}
\setcounter{equation}{0}
\numberwithin{equation}{section}
Liquid crystal is a soft matter between an isotropic fluid and an anisotropic crystalline solid\cite{Wang2021Acta}.
Liquid crystal phase's types include the nematic phase, 
smectic phase  and cholesteric phase. 
The dynamic continuum theory for liquid crystal flows was established by Ericksen and Leslie, 
thus leading to the Ericksen--Leslie system.
Ericksen (1961)\cite{Ericksen1961} proposed a continuum theory applicable to various liquid crystal states,
and Leslie (1968) \cite{Leslie1968} presented a theory specifically for nematic liquid crystals.
The Ericksen--Leslie system describes the flow characteristics of nematic liquid crystal materials
and emphasizes the effect of the orientation of liquid crystal molecules on the fluid velocity.

\subsection{The general Ericksen--Leslie system}
Suppose that $Q=[0,1]\times[0,1]$ is a unit square in $\mathbb{R}^2$,
and the torus $\mathbb{T}^2$  can be obtained by using the gluing map in topology.
In this article, 
we consider the general  Ericksen--Leslie system\cite{Lin2000ARMA,Lin2014CMP}
in $ Q \times \left[0,\infty\right)$,
which
reads as follows,
\begin{equation}\label{general_Ericksen--Leslie}
\left\{
\begin{aligned}
&\partial_t \mathbf{v}+\mathbf{v} \cdot \nabla \mathbf{v}
=-\nabla p+\frac{\gamma}{\operatorname{Re}} \Delta \mathbf{v}
+\frac{1-\gamma}{\operatorname{Re}} \nabla \cdot \sigma,\\
&\nabla \cdot \mathbf{v}=0, \\
&\mathbf{d} \times\left(\mathbf{h}-\gamma_1 \mathbf{N}+\gamma_2 \mathbf{D} \mathbf{d}\right)
=0, \\
&|\mathbf{d}|=1.
\end{aligned}
\right.
\end{equation}
Here the vector field  
$\mathbf{v}=\left(v_1,v_2\right)^T :
Q \times  \left[0,    \infty \right) \rightarrow \mathbb{R}^2$ 
represents  the velocity of the fluid,
and
the unit vector field  
$\mathbf{d}=\left(d_1, d_2\right)^T  :Q\times \left[0,    \infty \right) \rightarrow \mathbb{S}^1$ 
is the averaged orientation of liquid crystal
molecules.
The scalar function
$p$ is  the hydrodynamic  pressure,
and the constants ${\operatorname{Re}}$  (Reynolds number) and $\gamma \in(0,1)$ are parameters of the system.

The second-order tensor $\mathbf{N}$  is  
the rigid rotation part of the director changing rate by fluid vorticity
\begin{equation}\label{N_def}
\mathbf{N}
=\partial_t \mathbf{d}
+\mathbf{v} \cdot \nabla \mathbf{d} 
-\mathbf{\Omega} \mathbf{d}.
\end{equation} 
We use matrix multiplication to compute
$\mathbf{\Omega} \mathbf{d}$ and $\mathbf{D} \mathbf{d}$.
The second-order tensor $ \mathbf{D} $ denotes the symmetric part of the strain tensor $ \nabla \mathbf{v}$, which represents both stretching and compression.
$$
\mathbf{D}=\frac{1}{2}\left[\nabla \mathbf{v}+(\nabla \mathbf{v})^T\right]=\frac{1}{2}\left(\frac{\partial {v}_i}{ \partial x_j}+\frac{\partial {v}_j}{\partial x_i}\right). 
$$
The second-order tensor $\boldsymbol{\Omega}$ denotes  the skew-symmetric part of the strain tensor $ \nabla \mathbf{v}$,
which 
represents rotation.
$$\boldsymbol{\Omega}=\frac{1}{2}\left[\nabla \mathbf{v}-(\nabla \mathbf{v})^T\right]=\frac{1}{2}\left(\frac{\partial {v}_i}{\partial x_j}-\frac{\partial {v}_j}{\partial x_i}\right).
$$

The stress tensor $\sigma$ can be divided into two parts,
$$
\sigma=\sigma^L+\sigma^E,
$$
where $\sigma^L$ is the viscous stress tensor
\begin{equation}\label{sigma^L_def}
\begin{aligned}
\sigma^L
=& \;
\alpha_1\left(\mathbf{d}\otimes \mathbf{d}: \mathbf{D}\right) \mathbf{d}\otimes\mathbf{d} +\alpha_2\mathbf{N}\otimes \mathbf{d}
+\alpha_3 \mathbf{d}\otimes \mathbf{N} \\
&+\alpha_4 \mathbf{D}
+\alpha_5 \left(\mathbf{D} \mathbf{d}\right)\otimes \mathbf{d} 
+\alpha_6 \mathbf{d}\otimes\left(\mathbf{D}  \mathbf{d}\right)\\
=&\;
\alpha_1\left(d_k d_p D_{kp}\right) d_i d_j +\alpha_2 N_i d_j
+\alpha_3 d_i N_j \\
&+\alpha_4 D_{ij}
+\alpha_5\left(D_{ik} d_k\right)d_j 
+\alpha_6 d_i\left(D_{jk} d_k\right).
\end{aligned}
\end{equation}
We use the Einstein summation convention to represent $\sigma^L$. 
The viscous constants $\alpha_i$ $(i=1,    \ldots,  6)$  are called Leslie coefficients. 

The term $\sigma^E$ is the elastic stress tensor
\begin{equation}\label{sigma^E_def}
\sigma^E
=-\frac{\partial E_F}{\partial\left(\nabla \mathbf{d}\right)} \;\left(\nabla \mathbf{d}\right)^T,
\end{equation}
where $E_F=E_F\left(\mathbf{d},     \nabla \mathbf{d}\right)$ is the Oseen--Frank bulk energy density function\cite{alouges1997minimizing}
with the form
\begin{equation}\label{Oseen}
\begin{aligned}
E_F
= & \;\frac{k_1}{2}\left(\nabla \cdot \mathbf{d}\right)^2
+\frac{k_2}{2}\left|\mathbf{d} \cdot\left(\nabla \times \mathbf{d}\right)\right|^2 
+\frac{k_3}{2}\left|\mathbf{d} \times\left(\nabla \times \mathbf{d}\right)\right|^2\\
& +\frac{\left(k_2+k_4\right)}{2}
\left[\operatorname{tr}\left(\nabla \mathbf{d}\right)^2
-\left(\nabla \cdot \mathbf{d}\right)^2\right] .
\end{aligned}
\end{equation}
The Oseen--Frank theory\cite{oseen1933theory,frank1958liquid} uses a free energy functional to describe the elastic deformations of the molecular orientation,
providing a continuum model for the static equilibrium of nematic liquid crystals.
Here $k_i$ $(i=1,2,3,4)$ are elastic constants depending on the density and the temperature of the fluid.
The first three terms in (\ref{Oseen}) represent three different kinds of deformation, which are splay, twist and bend\cite{Wang2021Acta}, respectively,
and the last term of (\ref{Oseen}) is a
null Lagrangian  \cite{WXL2013ARMA},
$$\int_Q \frac{\left(k_2+k_4\right)}{2}\left[\operatorname{tr}(\nabla \mathbf{d})^2
-(\nabla \cdot \mathbf{d})^2\right] \mathrm{~d}x
=   
\frac{\left(k_2+k_4\right)}{2}\int_Q 
\nabla\cdot
\left[
(\nabla\mathbf{d})\mathbf{d}
-(\nabla \cdot \mathbf{d})\mathbf{d}
\right]
 \mathrm{~d}x.
$$ 
Without loss of generality, 
we suppose that $k_1=k_2=k_3=1$ in the isotropic case, 
and from (\ref{Oseen}),
we can get
$$
E_F
= \frac{1}{2}\left|\nabla\mathbf{d}\right|^2.
$$

We study the  Fréedericksz transition \cite{CKY2018SIAM,Sagues1985,kim2021SIAM,kim2020Non} for the
general Ericksen--Leslie system under an external magnetic field $\mathbf{H}=\mathrm{H}(1,    0)^T \in \mathbb{R}^2$, 
where $\mathrm{H}$ is a non-negative constant. 
The free energy functional 
for nematic liquid crystals under a magnetic field \cite{CKY2018SIAM,Aramaki2012,Atkin1997,Cohen1991,Haraux1993} is given by
\begin{equation}\label{functional_magnetic_field}
E_M
=\int_{Q}
E_F+
\frac{1}{2} \left[ \left|\mathbf{H}\right|^2-\left(\mathbf{H}\cdot \mathbf{d}\right)^2\right]
\mathrm{~d}x
=
\int_{Q}
 \frac{1}{2}\left|\nabla\mathbf{d}\right|^2+
\frac{1}{2} \left[ \left|\mathbf{H}\right|^2-\left(\mathbf{H}\cdot \mathbf{d}\right)^2\right]
\mathrm{~d}x
.    
\end{equation}
The vector field $\mathbf{h}$ is defined as follows,
\begin{equation}\label{h_def}
\mathbf{h}
=-\frac{\delta E_M}{\delta \mathbf{d}}
=\nabla \cdot \frac{\partial E_M}{\partial \nabla \mathbf{d}}-\frac{\partial E_M}{\partial \mathbf{d}}
=\Delta \mathbf{d}+  
\left|\nabla\mathbf{d}\right|^2 \mathbf{d}
+
\left(\mathbf{H}\cdot \mathbf{d}\right)\mathbf{H}
-\left(\mathbf{H}\cdot \mathbf{d}\right)^2\mathbf{d}.
\end{equation}

The Leslie coefficients in
system (\ref{general_Ericksen--Leslie}) satisfy 
Parodi's relation \cite{Parodi1970},
$$
 \alpha_2+\alpha_3=\alpha_6-\alpha_5.
$$
Here the Parodi's relation, 
providing a constraint on the Leslie coefficients in the stress tensor of nematic liquid crystals,
is derived from the Onsager reciprocal relations, which express the equality of certain cross coefficients between flows and forces in thermodynamic systems that are not in equilibrium.

The Leslie coefficients and $\gamma_1, \gamma_2$ in \eqref{general_Ericksen--Leslie} obey the following relations\cite{Leslie1968}, 
\begin{equation}\label{gamma_def}
\begin{aligned}
	\gamma_1=\alpha_3-\alpha_2, \quad \gamma_2=\alpha_5-\alpha_6,
\end{aligned}
\end{equation}
which are necessary conditions to ensure that the equation of motion holds identically.

To simplify the third equation \eqref{general_Ericksen--Leslie},
denoted by $\mu_1=\dfrac{1}{\gamma_1}$ 
and $\mu_2=\dfrac{\gamma_2}{\gamma_1}$.
From the Ericksen--Leslie system (\ref{general_Ericksen--Leslie}),
we have the following equivalent system 
in $ Q \times \left[0,\infty\right)$,
\begin{equation}\label{Ericksen_Leslie_new_formulation}
\left\{
\begin{aligned}
\partial_t \mathbf{v}+\mathbf{v} \cdot \nabla \mathbf{v}
=
&
-\nabla p
+\frac{\gamma}{\operatorname{Re}} \Delta \mathbf{v}
+\frac{1-\gamma}{\operatorname{Re}}\nabla \cdot\left(
\sigma^L
+\sigma^E\right),    
\\
\nabla \cdot \mathbf{v}
= 
 & 
\; 0,    \\
\partial_t \mathbf{d}+\mathbf{v} \cdot \nabla \mathbf{d}
=
&
\;\mu_1\left[\mathbf{h}-\left(\mathbf{h}\cdot\mathbf{d}\right)\mathbf{d} \right]
+\left(\mathbf{I}-\mathbf{d}\otimes\mathbf{d}\right)\left(\boldsymbol{\Omega}  \mathbf{d}
+\mu_2 \mathbf{D} \mathbf{d}\right).
\end{aligned}
\right.
\end{equation}
The initial conditions are as follows
\begin{equation}\label{initial_condition_1}
\mathbf{v}|_{t=0} = \mathbf{v}_0,    \quad
\nabla\cdot \mathbf{v}_0 = 0,    \quad
\mathbf{d}|_{t=0} = \mathbf{d}_0,    \quad
| \mathbf{d}_0| = 1,    \quad
 x\in Q.
\end{equation}
For incompressible nematic liquid crystal fluids  in torus $\mathbb{T}^2$, 
we have the following boundary conditions before utilizing the gluing map in topology,
\begin{equation}\label{boundary_condition_1}
\mathbf{v}(x+\mathbf{e}_i,    t) = \mathbf{v}\left(x,    t\right),   \quad  
\mathbf{d}(x+\mathbf{e}_i,    t) = (-1)^{a_i}\mathbf{d}\left(x,    t\right),      \quad  
\left(x,    t\right)\in \partial Q \times  \left[0,    \infty \right).
\end{equation}
The unit vectors $\mathbf{e}_i \;(i=1,2)$ form the canonical basis of $\mathbb{R}^2$,  
and the integers $a_i\; (i=1,2)$ are constants depending on
the initial condition
$$
\mathbf{d}_0(x+\mathbf{e}_i) = \left(-1\right)^{a_i}\mathbf{d}_0(x),    \quad
x\in \partial Q.
$$ 
Because of the symmetry about the center of liquid crystal molecules, when a liquid crystal molecule rotates by an angle of \( a_i\pi \) about its center, it coincides with its original state.

\subsection{Brief review of previous results}

The Ginzburg--Landau energy functional was first introduced in 1950 \cite{GinzburgLandau1950}. 
Let \( \epsilon > 0 \) and \( \Omega \subset \mathbb{R}^n \) (\( n = 2, 3 \)) be a bounded domain, the Ginzburg--Landau energy functional for \( u: \Omega \rightarrow \mathbb{R}^3 \) is defined as follows,
\[
E_{\epsilon}(u; \Omega) = \int_{\Omega} \frac{1}{2} |\nabla u|^2 + \frac{1}{4 \epsilon^2} (1 - |u|^2)^2 \, \mathrm{d}x.
\]
Under the framework of Ginzburg--Landau approximation to the simplified Ericksen--Leslie system, Lin and Liu (1995) in \cite{Lin1995CPAM}
studied the existence of the global weak solutions and classical solutions in 2D and 3D bounded regions.
Wu (2010)  in \cite{Wu2010DCDS} demonstrated the long-term behavior of solution.

For the general Ericksen--Leslie system
under the framework of Ginzburg-Landau approximation, 
Lin and Liu (2000) in \cite{Lin2000ARMA} proved the existence of the global weak solutions  and asymptotic stability of classical solutions  and weak solutions in 2D and 3D bounded regions.
Sun and Liu (2009) in \cite{Sun2009DCDS} established the well-posedness of classical solutions in 2D and 3D bounded regions under the condition, where  Leslie stress tensor \(\sigma^L = \left[ \Delta \mathbf{d} - \frac{1}{\epsilon^2} \left( |\mathbf{d}|^2 - 1 \right) \mathbf{d} \right] \otimes \mathbf{d}\).
Feng et al. (2020) in \cite{Hong2020SIMA}
proved that the
solutions of the Ginzburg--Landau approximation system converge smoothly to the solution of the
Ericksen--Leslie system.
Wu and De Anna (2023) in \cite{Wu2023cvpde} demonstrated the uniqueness of global weak solutions under periodic boundary conditions in two dimensions.
Wu et al. (2012) in \cite{Wu2012CVPDE} studied the asymptotic behavior of global classical solutions 
to the nematic liquid crystal ﬂows  under kinematic transports for molecules of different shapes
with periodic boundary conditions.
Wu et al. (2013) in \cite{WXL2013ARMA} established the long-time behavior of global solutions under periodic boundary conditions in three dimensions with large viscosity and in two dimensions under periodic boundary conditions.

In the simplified Ericksen--Leslie system of fixed liquid crystal  modulus $|\mathbf{d}|=1$,
Lin et al. (2010) in \cite{Lin2010ARMA} proved the existence of global weak solutions in a bounded domain of \( \mathbb{R}^2 \) in their 2010 paper, and showed that  weak solutions are smooth except at finitely many singular points. 
Lin and Wang (2010) in \cite{Lin2010CAMS} demonstrated the uniqueness of weak solutions in a bounded domain of \( \mathbb{R}^2 \). Hong (2011) in \cite{HMC2011CVPDE} established the existence of global weak solutions, and these weak solutions are smooth except at finitely many singular points. 
Xu and Zhang  (2012) in  \cite{ZZZ2012JDE} proved the regularity and uniqueness of weak solutions on \( \mathbb{R}^2 \) . 
Fan and  Li (2015) in \cite{LiJK2015JMAA} presented a blow-up criterion for global strong solutions in three dimensions. 
Lin and Wang proved (2016) in \cite{Lin2016CPAM}  the existence of global weak solutions in a bounded domain of \( \mathbb{R}^3 \). 
Huang et al. (2016) in \cite{huangtao2016ARMA} provided examples of singularities occurring in finite time. 
Gong et al. (2017) in \cite{LiJK2017JDE} demonstrated the non-uniqueness of weak solutions in a bounded domain of \( \mathbb{R}^3 \).
Chen et al. (2018) in \cite{CKY2018SIAM} proved the long-time behavior of global strong solutions for the system under the influence of a magnetic field in a bounded domain of \( \mathbb{R}^2 \).

For the general Ericksen--Leslie system of fixed liquid crystal molecular modulus, 
Wang et al. (2013) in \cite{ZZZ2013ARMA} established the local well-posedness and global well-posedness for small initial values of the system on \( \mathbb{R}^3 \). 
Huang et al. (2014) in   \cite{Lin2014CMP} proved the existence and regularity of global weak solutions for the system on \( \mathbb{R}^2 \). 
Gong et al. (2015) in \cite{Gong2015Non} demonstrated the existence of strong solutions for the general form of the Ericksen--Leslie system. 
Wang and Xu (2014)  in \cite{WangchangyouXvxiangJFA2014}  proved the existence of global weak solutions for the system on a two-dimensional sphere. 
Huang et al. (2016) in \cite{huangtao2016ARMA} among others, provided examples of singularities occurring in finite time.
Huang and Qu (2013) in \cite{Huang2023MMAS} proved the long-time behavior of global weak solutions on a two-dimensional sphere.

\subsection{Notations}

We recall the well established functional settings
for periodic problems
\cite{Temam1995navier}:
$$
\begin{aligned}
H^m(Q) & = \left\{ u \mid u \in H^m_{\text{loc}}\left(\mathbb{R}^2, \mathbb{R}\right) \right\}, \;\text{ where }\; L^2(Q) = H^0(Q), \\
H^m_p(Q) & = \left\{ u \in H^m_{\text{loc}}\left(\mathbb{R}^2, \mathbb{R}\right) \mid u(x + \mathbf{e}_i) = u(x) \right\}, \\
\dot{H}^m_p(Q) & = H^m_p(Q) \cap \left\{ u : \int_Q u(x) \, \mathrm{d}x = 0 \right\}, \\
H & = \left\{ \mathbf{v} \in L^2_p(Q) \mid \nabla \cdot \mathbf{v} = 0 \right\}, \;\text{ where }\; L^2_p(Q) = H^0_p(Q), \\
V & = \left\{ \mathbf{v} \in \dot{H}^1_p(Q) \mid \nabla \cdot \mathbf{v} = 0 \right\}, \\
V' & =  \text{the  dual space of}\; V.
\end{aligned}
$$

Suppose that $\mathbf{d}=(\sin\theta,\cos\theta)^{T} \in \mathbb{S}^1$.
To describe nematic liquid crystals  in $\mathbb{T}^2$, we need the following definitions for $\mathbf{d}$ and $\theta$,  respectively, which are generally not periodic or in $\mathbb{T}^2$.
$$
\begin{aligned}
H^m_{sp}(Q)  :=\left\{u \in H^m_{loc}\left(\mathbb{R}^2, \mathbb{R}\right) \mid u\left(x+\mathbf{e}_i\right)=(-1)^{a_i} u(x)\;(1\leqslant i\leqslant n),\; 
\text{for} \;
x\in \partial Q,\;
a_i \in \mathbb{Z} \right\},\\
H^m_{sq}(Q)  :=\left\{u \in H^m_{loc}\left(\mathbb{R}^2, \mathbb{R}\right) \mid u\left(x+\mathbf{e}_i\right)= u(x)+a_i\pi\;(1\leqslant i\leqslant n),\;
\text{for} \;
x\in \partial Q,\;
a_i \in \mathbb{Z} \right\}.
\end{aligned}
$$

\paragraph{}
We need the following notation:

(i) If $X(Q)$ is a function space defined on $Q$, we simply denote $X(Q)$ by $X$.

(ii) 
 For any Banach space $X$, use $X$ to denote the vector space $(X)^r$, where $r$ is a positive integer.

(iii)
The notation $A \lesssim B$ or $A \leqslant C B$ means there exists a constant $C > 0$ such that $A \leqslant C B$.

(iv)
For positive constants $c$ and $m$, the inequality $c \ll m$ indicates that $c$ is much smaller than $m$.

(v)
For tensors $A=\left(a_{ij}\right)_{n\times n}$ and  $B=\left(b_{ij}\right)_{n\times n}$,
we have 
$A: B
= 
a_{ij}b_{ij}.$

(vi) 
  The dual product refers to a bilinear map between a Banach space \( X \) and its dual space \( X^* \)
   \[
   \langle \cdot, \cdot \rangle_{X^*, X }: X^* \times X \to \mathbb{R}.
   \]
   For every  $ f \in X^* $ and any $ x \in X$,
   we have \( \langle f, x \rangle_{X^*, X } = f(x).\)

\paragraph{}
Let \( 0 < T < \infty \). The local weak solution \( (\mathbf{v}, \mathbf{d}) \) to system (\ref{general_Ericksen--Leslie}) satisfies the initial conditions (\ref{initial_condition_1}) and the boundary conditions (\ref{boundary_condition_1}), and
\[
\mathbf{v} \in L^2\left(0, T; V\right), \quad \mathbf{d} \in L^2\left(0, T; H^2_{sp}\right).
\]
For any \(\eta \in C^\infty([0, T])\) satisfying \(\eta(T) = 0\), and for any \(\xi \in V\) and \(\tilde{\xi} \in H^1_{sp}\), it holds that

$$
\begin{aligned}
-\eta(0) \int_{Q}  \mathbf{v}_0 \cdot \xi 
\mathrm{~d}x
=& -\int_0^T \eta^{\prime} \int_{Q}  \mathbf{v}\cdot \xi
\mathrm{~d}x \mathrm{~d} t
+\int_0^T \eta\int_{Q} \mathbf{v} \otimes \mathbf{v}:  \nabla \xi
\mathrm{~d}x \mathrm{~d} t
\\
&+\frac{\gamma}{\operatorname{Re}}\int_0^T \eta\int_{Q}\nabla \mathbf{v}: \nabla \xi
\mathrm{~d}x \mathrm{~d} t\\
&+\frac{1-\gamma}{\operatorname{Re}}\int_0^T \eta\int_{Q} \sigma:  \nabla \xi 
\mathrm{~d}x \mathrm{~d} t,\\
-\eta(0) \int_{Q}  \mathbf{d}_0\cdot \tilde{\xi} 
\mathrm{~d}x 
=&
- \int_0^T \eta^{\prime} \int_{Q}  \mathbf{d}\cdot  \tilde{\xi} 
\mathrm{~d}x \mathrm{~d} t
+\int_0^T  \eta\int_{Q}  \mathbf{v} \cdot \nabla \mathbf{d}\cdot \tilde{\xi} 
\mathrm{~d}x \mathrm{~d} t\\
&+\mu_1\int_0^T \eta \int_{Q} \mathbf{h}+|\nabla \mathbf{d}|^2 \mathbf{d}-(\mathbf{H}\cdot\mathbf{d})^2 \mathbf{d}\cdot \tilde{\xi}
\mathrm{~d}x \mathrm{~d} t
\\
&+\int_0^T \eta  \int_{Q}  (\mathbf{I}-\mathbf{d}\otimes\mathbf{d})\left(\boldsymbol{\Omega}  \mathbf{d}+\mu_2 \mathbf{D}  \mathbf{d}\right)\cdot \tilde{\xi}
\mathrm{~d}x \mathrm{~d} t.
\end{aligned}
$$
For both {isotropic} and {anisotropic} states, the meanings of both \(\sigma\) and \(\mathbf{h}\) are different, but the forms of the equations for \((\mathbf{v}, \mathbf{d})\) remain the same. Therefore, the definition of weak solutions is identical in both cases.

\paragraph{}
Let \( 0 < T < \infty \). 
A {local strong solution} \( (\mathbf{v}, \mathbf{d}) \) to system (\ref{general_Ericksen--Leslie}) 
with the initial conditions (\ref{initial_condition_1}) and the boundary conditions (\ref{boundary_condition_1}) is defined as follows.
 For  \( (x, t) \) in \( Q \times [0, T] \) $a.e.$, the solution \( (\mathbf{v}, \mathbf{d}) \) satisfies (\ref{general_Ericksen--Leslie}) with the initial conditions (\ref{initial_condition_1}) and the boundary conditions (\ref{boundary_condition_1}) and 
\[
\begin{aligned}
& \mathbf{v} \in L^{\infty}\left( 0, T ; V\right) \cap L^2\left(0, T; H^2_p\right), \\
& \mathbf{d} \in L^{\infty}\left(0, T; H^2_{sp}\right) \cap L^2\left(0, T ; H^3_{sp}\right).
\end{aligned}
\]

A {global strong solution} \( (\mathbf{v}, \mathbf{d}) \) to system (\ref{general_Ericksen--Leslie}) 
with the initial conditions (\ref{initial_condition_1}) and the boundary conditions (\ref{boundary_condition_1}) is defined as follows.
 For  \( (x, t) \) in \( Q \times \left[0, \infty\right) \) $a.e.$, the solution \( (\mathbf{v}, \mathbf{d}) \) satisfies (\ref{general_Ericksen--Leslie}) with the initial conditions (\ref{initial_condition_1}) and the boundary conditions (\ref{boundary_condition_1}) and 
$$
\begin{aligned}
& \mathbf{v} \in L^{\infty}\left( 0, \infty ; V\right) 
   \cap L^2\left(0, \infty ; H^2_p\right), \\
& \mathbf{d} \in L^{\infty}\left(0, \infty; H^2_{sp}\right) 
    \cap L^2_{loc}\left(0, \infty ; H^3_{sp}\right).
\end{aligned}
$$

\paragraph{}
Let $\left\{\lambda_j \mid 0=\lambda_1<\lambda_2\leqslant\lambda_3\leqslant\cdots\leqslant\lambda_n\leqslant\cdots, j=1,2,3, \ldots\right\}$ be the set of eigenvalues  of the periodic boundary value problem 
\begin{equation}\label{periodic_boundary_eigenvalue}
 \left\{\begin{array}{rlr}
-\Delta \psi_j (x)& =\lambda_j  \psi_j(x) & \text { in } Q,\\
\psi_j(x+\mathbf{e}_i) & =\psi_j(x) & \text { on } \partial Q.
\end{array}\right.   
\end{equation}
Let $\lambda_2>0$ be the principal eigenvalue,
and the following conclusion can be drawn in\cite{Evans2010partial},
\begin{equation}\label{lambda_2_def}
\lambda_2
=\inf\limits_{\substack{\psi \in H^1_p,\\
\psi\not\equiv \int_{Q}\psi
} }
\frac{\int_{Q}\left|\nabla\psi\right|^2\mathrm{~d}x}{\int_{Q}\left|\psi-\int_{Q}\psi\mathrm{~d}x\right|^2\mathrm{~d}x}.
\end{equation}

\subsection{Main results}

The novelty of this paper lies in the incorporation of a magnetic field term into the general Ericksen–Leslie system.
The primary mathematical challenge is that the well-posedness, energy estimates, and long-time behavior become significantly more complex due to the system's nonlinear structure.
revious studies on the Fréedericksz transition were limited to simplified Ericksen–Leslie systems \cite{CKY2018SIAM,kim2020Non,kim2021SIAM}.
Therefore, a key mathematical contribution of this work is to provide a mathematical analysis of the Fréedericksz transition within the framework of the general Ericksen–Leslie system.

The general Ericksen–Leslie system with the constraint of constant molecular length has intricate coupling relationships.
Previous results on high-order estimates and well-posedness of global strong solutions were only available under the Ginzburg–Landau approximation\cite{WXL2013ARMA}.
Therefore, it is important to handle the case where $|\mathbf{d}| = 1$.

Obtaining an estimate for $\|\theta\|$  with respect to time t directly from the basic energy equality and high-order estimates is not feasible.
First, a maximum principle for $\theta$ is unavailable due to the boundary conditions and nonlinearities in the system.
Second, energy and high-order estimates only provide bounds for $\mathbf{d}$.
Therefore, achieving a uniform-in-time estimate for $\|\theta\|$ is a major difficulty in this work.
Besides, the Lojasiewicz–Simon inequality can only be applied to the equation for $\theta$, not to the equation for $\mathbf{d}$, which is restricted to the submanifold \( \mathbf{S}^1 \).
  
The director field $\mathbf{d}$ is not necessarily periodic on $Q$, and the corresponding angle $\theta$ (via $\mathbf{d} = (\sin\theta,\cos\theta)^T$) is also not necessarily periodic on $Q$.
This poses a significant difficulty for the well-posedness of both global strong solutions and steady-state solutions.
To consider  isotropic liquid crystals,
suppose that $\mathbf{d}=(\sin\theta,\cos\theta)^T$ and $k_1=k_3=1$.
If $\mathbf{v}=0$ and $\theta=\theta(x)$,
then the system (\ref{general_Ericksen--Leslie})
with initial conditions (\ref{initial_condition_1}) and 
boundary conditions (\ref{boundary_condition_1})
can be reduced to the following steady-state problem including the elliptic sine-Gorden  equation,
\begin{equation}\label{iso_elliptic_thm_first}
\left\{
\begin{aligned}
-\Delta \theta
=& \;
\frac{\mathrm{H}^2}{2}
\sin 2\theta
,
\quad
&& x\in Q,
\\
\theta(x+\mathbf{e}_i) 
=& \;
\theta(x)+a_i\pi,
\quad
&& x\in \partial Q.
\end{aligned}\right.
\end{equation}
If $a_1=a_2=0$, then $\theta$ staisfies periodic boundary conditions and shares similar conclusions  of \cite{CKY2018SIAM}, where the principal eigenvalue of the sine-Gordon equation serves as the critical magnetic field strength to characterize the solution branches.
However, function $\theta$ not exhibit a counterpart to the $(0, \pi)$  interval observed in \cite{CKY2018SIAM}.
Besides, the condition $|a_1|+|a_2|>0$ leads to the novel and challenging scenario characterized by modulo-$\pi$ periodic boundary conditions for $\theta$.

\begin{remark}
If $\theta(x)$ is a solution of (\ref{iso_elliptic_thm_first}), 
for each constant vector $w \in \mathbb{R}^2$ and any $k\in \mathbb{Z}$,
the function $\theta(x+w)+k\pi$ can be another solution of (\ref{iso_elliptic_thm_first}).
Therefore, the uniqueness of $\theta(x)$ in (\ref{iso_elliptic_thm_first}) means that 
if $\theta_1(x)$ and $\theta_2(x)$ are two solutions satisfying (\ref{iso_elliptic_thm_first}), there must have a constant vector $w \in \mathbb{R}^2$ and 
an integer constant $k$
such that 
$\theta_1(x)=\theta_2(x+w)+k\pi$.
Moreover, by the standard bootstrap method, it can be concluded that every solution to (\ref{iso_elliptic_thm_first}) is a smooth function.    
\end{remark}

\begin{theorem}\label{Theorem_1.1}
For steady-state problem (\ref{iso_elliptic_thm_first}), we have the following conclusions.

(i)
Assume $a_1=a_2=0$.
If $\mathrm{H}^2>\lambda_2$, there exists at least one nonconstant solution to (\ref{iso_elliptic_thm_first}).
If $\mathrm{H}^2\leqslant\lambda_2$, 
all solutions to  (\ref{iso_elliptic_thm_first}) are
constants.

(ii)
Assume $|a_1|+|a_2|>0$,
then for any constant $\mathrm{H}$,
there exists at least one  solution to (\ref{iso_elliptic_thm_first}), and all solutions to  (\ref{iso_elliptic_thm_first}) are
nonconstant solutions.
If $\mathrm{H}^2\leqslant\lambda_2$, 
the solution to  (\ref{iso_elliptic_thm_first}) is unique.
\end{theorem}

Subtract the smooth solution to (\ref{iso_elliptic_thm_first}) from the molecular angle $\theta$, transforming it into a function with periodic boundary conditions.
This facilitates the proof of existence for the global strong solutions.
Reformulate the equation for the director $\mathbf{d}$ into an equation for the angle $\theta$.
Apply the Galerkin method to the $\theta$-equation to construct solutions.

\begin{theorem}\label{Theorem_1.2}
Assume that
$\mathbf{d}\in \mathbb{S}^1$, $ k_1=k_3=1$ and
initial data $\left(\mathbf{v}_0, \mathbf{d}_0\right) \in V\times H_{sp}^2$, 
there exists a unique global strong solution 
$\left(\mathbf{v}, \mathbf{d}\right) $ to 
the general system(\ref{general_Ericksen--Leslie}) with the initial conditions (\ref{initial_condition_1}) and the boundary conditions (\ref{boundary_condition_1}).
Furthermore, 
$\left(\mathbf{v}, \mathbf{d}\right) \in C^\infty(Q\times (0,\infty))$.
\end{theorem}

To prove Theorem \ref{Theorem_1.2}, we get the following estimates,
$$
\begin{aligned}
    \|\nabla^2\mathbf{d}\|
\leqslant&\;
    C\left(\| \Delta\mathbf{d}+|\nabla\mathbf{d}|^2\mathbf{d} \|
+1
\right),
\\
    \|\nabla^3\mathbf{d}\|
\leqslant&\;
    C\left(
    \| \Delta \mathbf{d} + |\nabla\mathbf{d}|^2 \mathbf{d} \|^2
    +
    \| \nabla(\Delta\mathbf{d}+|\nabla\mathbf{d}|^2\mathbf{d}) \|
    +1 \right),
\end{aligned}
$$
where $C$ is a constant depending on $\|\nabla\mathbf{d}\|^2$.
After establishing theorems for high-order energy estimates.
Leverage the coupling within the system to complete the proof of these high-order estimates, which are crucial for establishing the well-posedness of the global strong solution $ (\mathbf{v}, \mathbf{d})$.

\paragraph{}
However, $\|\theta(t)\|$ does not have a uniform boundary regarding time.
After obtaining the boundness of $(\mathbf{v}, \theta)$ by the relationship between $\mathbf{d}$ and $\theta$ and \L ojasiewicz--Simon inequality,
we can get the long-time behavior of $(\mathbf{v}, \theta)$, which illustrates 
the long-time behavior of $(\mathbf{v}, \mathbf{d})$.

\begin{theorem}\label{Theorem1.3}
Assume that
$\mathbf{d} \in \mathbb{S}^1$,
$k_1=k_3=1$
and initial data $\left(\mathbf{v}_0, \mathbf{d}_0\right) \in V\times H_{sp}^2$.
Then there exist a smooth solution $\theta_{\infty}(x)$ 
and a constant $\nu \in \left(0,\frac{1}{2}\right]$ such that
the global solution
$(\mathbf{v}, \mathbf{d})$ staisfies the following inequalities.

(i) If \(\nu \in \left(0, \frac{1}{2}\right)\), then for all \(t > 0\),
\begin{equation}
\|\mathbf{v}(t)\|_{H^1} + \left\|\theta(t) - \theta_{\infty}\right\|_{H^2} \lesssim (1 + t)^{-\frac{\nu}{1 - 2\nu}}.
\end{equation}

(ii) If \(\nu = \frac{1}{2}\), then for all \(t > 0\), there exists a constant \(c > 0\) such that
\begin{equation}
\|\mathbf{v}(t)\|_{H^1} + \left\|\theta(t) - \theta_{\infty}\right\|_{H^2} \lesssim e^{-ct}.
\end{equation}
The constant \(c\) depends on $Q$,  $ |\mathbf{H}|$, $\left\|\mathbf{d}_0\right\|_{H^2}$,  
$\left\|\mathbf{v}_0\right\|_{H^1}$
and
coefficients of the system (\ref{general_Ericksen--Leslie}).

(iii)
If $ \mathrm{H}^2<\lambda_2 $, then 
$\nu=\frac{1}{2}$.
\end{theorem}

\paragraph{}
The paper is organized as follows.
In Section 2, we provide necessary lemmas, 
derive  the general Ericksen--Leslie system for  \((\mathbf{v}, \mathbf{d})\) and 
\((\mathbf{v}, \theta)\), respectively,
and
obtain the basic energy law. 
In Section 3, 
we discuss the existence and uniqueness of both constant and nonconstant solutions to the steady-state problem for the cases of periodic and modulo-$\pi$ periodic boundary conditions, respectively. 
In Section 4,
we prove higher-order estimates to obtain uniform a priori estimate.
In Section 5,
we utilize the semi-Galerkin approximation method with employing uniform a priori estimates to demonstrate the well-posedness of the global strong solutions to the general Ericksen--Leslie system, 
In Section 6,
we obtain a uniform bound on  \(|\theta|\)  with respect to time  $t>0$ and  the long-time behavior of  \((\mathbf{v}, \theta)\).

\section{Basic energy law}

\subsection{\texorpdfstring{System for $\left(\mathbf{v}, \mathbf{d}\right)$}{System for} 
}

\setcounter{equation}{0}
\numberwithin{equation}{section}

\begin{lemma}\label{lem_inequalities}
Suppose that  $\Omega\subset \mathbb{R}^2$ .
We report the Ladyzhenskaya\cite{Ladyzhenskaya1963mathematical,Ladyzhenskaya1968linear},
Agmon\cite{Stein1970Agmon},
Gagliardo-Nirenberg\cite{Zheng2004nonlinear} 
and trace interpolation inequalities\cite{Zheng2004nonlinear}.
$$
\begin{aligned}
\|f\|_{L^4(\Omega)} 
\leqslant & \;
C\|f\|_{L^2(\Omega)}^{\frac{1}{2}}\|f\|_{H^1(\Omega)}^{\frac{1}{2}}, \quad &&  f \in H^1(\Omega),  \\
\|f\|_{L^{\infty}(\Omega)} 
\leqslant & \;
C\|f\|_{L^2(\Omega)}^{\frac{1}{2}}\|f\|_{H^2(\Omega)}^{\frac{1}{2}}, \quad &&  f \in H^2(\Omega), \\
\|f\|_{L^p(\Omega)} 
\leqslant & \;
C p^{\frac{1}{2}}\|f\|_{L^2(\Omega)}^{\frac{2}{p}}\|f\|_{H^1(\Omega)}^{1-\frac{2}{p}},  \quad &&  f \in H^1(\Omega),\; 2 \leqslant p<\infty,\\
\|f\|_{L^2(\partial \Omega)} 
\leqslant & \;
C\|f\|_{L^2(\Omega)}^{\frac{1}{2}}\|f\|_{H^1(\Omega)}^{\frac{1}{2}}, \quad &&  f \in H^1(\Omega).
\end{aligned}
$$
The positive constant $C$ depends on $\Omega$.
\end{lemma}

\begin{lemma}\label{Lemma_N}
Suppose that $\mathbf{d}_0 \in \mathbb{S}^1$
satisfies the compatibility condition
$$\mathbf{d}_0(x+\mathbf{e}_i, t)
=(-1)^{a_i}\mathbf{d}_0(x, t), 
\quad
 (x,t) \in \partial Q \times \left[0,\infty\right),
$$ 
where constants $a_i\; (i=1,2)$ are integers.
From  initial  conditions (\ref{initial_condition_1})
and  boundary conditions (\ref{boundary_condition_1}),
\begin{equation}\label{Ericksen_Leslie_d_initial_boundary_condition}
\left\{
\begin{aligned}
\mathbf{d}(x,  0)
=&\;\mathbf{d}_0(x),  \quad        
&&x\in Q,
\\
\mathbf{d}(x+\mathbf{e}_i, t)
=& (-1)^{a_i}\mathbf{d}(x, t),   \quad 
&&(x,t) \in   \partial Q \times \left[0,\infty\right). 
\end{aligned}\right.
\end{equation}
Thanks to the system (\ref{general_Ericksen--Leslie}), we have the following equation about $\mathbf{d}\in \mathbb{S}^1$
in $ Q \times \left[0,\infty\right)$,
\begin{equation}\label{Ericksen_Leslie_old_d}
\begin{aligned}
\mathbf{d} \times\left(\mathbf{h}-\gamma_1 \mathbf{N}
+\gamma_2 \mathbf{D} \mathbf{d}\right)
=0 .
\end{aligned}
\end{equation}
The above problem is equivalent to the following equation
about $\mathbf{d}\in \mathbb{S}^1$
in $ Q \times \left[0,\infty\right)$,
\begin{equation}\label{Ericksen_Leslie_new_d_new_N}
\begin{aligned}
\mathbf{N}
=
\mu_1
\left[\mathbf{h}-\left(\mathbf{h}\cdot\mathbf{d}\right)\mathbf{d} \right]
+
\mu_2
\left[
\mathbf{D}  \mathbf{d}-(\mathbf{D}: \mathbf{d}\otimes\mathbf{d} ) \mathbf{d}\right].
\end{aligned}
\end{equation}
\end{lemma}

\paragraph{Proof.}
The idea of Lemma \ref{Lemma_N}  originates from \cite{ZZZ2013ARMA}. For the completeness of this paper, we provide the proof here, which mainly follows the approach in \cite{Wang2024DCDS}.

We introduce the Lagrange multiplier $\lambda$, 
then
the  equation o(\ref{Ericksen_Leslie_old_d}) is equivalent to
\begin{equation}\label{lambda_d_N_h}
\mathbf{h} - \gamma_1 \mathbf{N} + \gamma_2 \mathbf{D} \mathbf{d}
=\lambda \mathbf{d}.
\end{equation}
Since $\mathbf{d} \in \mathbb{S}^1$ and $\mathbf{N} \cdot \mathbf{d}=0$,
after
multiplying  $\mathbf{d}$ on both sides of (\ref{lambda_d_N_h}),
we have
$$
\lambda=
\mathbf{h}\cdot \mathbf{d} + 
\gamma_2(\mathbf{D}: \mathbf{d}\otimes\mathbf{d}).
$$
Substituting $\lambda$ back into (\ref{lambda_d_N_h}),
we can obtain the equation (\ref{Ericksen_Leslie_new_d_new_N}). 

On the other hand, 
by substituting $\mathbf{N}$ from (\ref{Ericksen_Leslie_new_d_new_N})  into (\ref{Ericksen_Leslie_old_d}),
we can obtain the proof of Lemma \ref{Lemma_N}. 
$\Box$

\paragraph{}
Substituting the definition of $\mathbf{N}$ from equation (\ref{N_def}) into the expression for $\mathbf{N}$ in  (\ref{Ericksen_Leslie_new_d_new_N}),
we have
\begin{equation}\label{First_d_equation}
\begin{aligned}
 \partial_t \mathbf{d}+\mathbf{v} \cdot \nabla \mathbf{d}
=
\mu_1\left[\mathbf{h}-\left(\mathbf{h}\cdot\mathbf{d}\right)\mathbf{d} \right]
+(\mathbf{I}-\mathbf{d}\otimes\mathbf{d})\left(\boldsymbol{\Omega} \mathbf{d}
+\mu_2 \mathbf{D}  \mathbf{d}\right).  
\end{aligned}
\end{equation}
The requirement $\mu_1>0$ is essential to ensure the necessary condition for the dissipation law of the director field $\mathbf{d}$.

\paragraph{}
Substituting the expression  of $\mathbf{N}$ from equation (\ref{Ericksen_Leslie_new_d_new_N})   into 
the expression  of $\sigma$ from equation (\ref{sigma^L_def}), 
we have
\begin{equation}\label{sigma^L_first}
\begin{aligned}
\sigma^L
= 
& \;
\beta_1 \left(\mathbf{D}: \mathbf{d}\otimes\mathbf{d} \right) \mathbf{d}\otimes\mathbf{d}
+\frac{1}{2}\left(-1-\mu_2\right) \mathbf{h}\otimes\mathbf{d}
+\frac{1}{2}\left(1-\mu_2\right) \mathbf{d}\otimes\mathbf{h}\\
& \;
+\mu_2 \left(\mathbf{h}\cdot\mathbf{d} \right) \mathbf{d}\otimes\mathbf{d} 
+\beta_2 \mathbf{D}
+\frac{\beta_3}{2}\left[\mathbf{d}\otimes(\mathbf{D} \mathbf{d}) +(\mathbf{D}\mathbf{d})\otimes\mathbf{d}\right]\\
= 
& \;
\beta_1(\mathbf{D}: \mathbf{d}\otimes\mathbf{d}) \mathbf{d}\otimes\mathbf{d}
+\beta_2 \mathbf{D}
+\frac{\beta_3}{2}(\mathbf{d}\otimes(\mathbf{D} \mathbf{d}) +(\mathbf{D}\mathbf{d})\otimes\mathbf{d})\\
& \; +\frac{1}{2}\left(-1-\mu_2\right) \left[\mathbf{h}-\left(\mathbf{h}\cdot\mathbf{d}\right)\mathbf{d} \right]\otimes\mathbf{d}
+\frac{1}{2}\left(1-\mu_2\right) \mathbf{d}\otimes\left[\mathbf{h}-\left(\mathbf{h}\cdot\mathbf{d}\right)\mathbf{d} \right],
\end{aligned}
\end{equation}
where
\begin{equation}\label{beta_def}
    \beta_1=\alpha_1+\frac{\gamma_2^2}{\gamma_1},  \quad \beta_2=\alpha_4,  \quad \beta_3=\alpha_5+\alpha_6-\frac{\gamma_2^2}{\gamma_1} .
\end{equation}

\begin{lemma}\label{beta_lemma}
For any symmetric trace free matrix $\mathbf{D}$ and  vector field $\mathbf{d}:Q\times \left[0,    \infty \right) \rightarrow \mathbb{S}^1$,
the following dissipation relation holds
\begin{equation}\label{beta_inequality_first}
\beta_1(\mathbf{d}\otimes\mathbf{d}: \mathbf{D})^2
+\beta_2 \mathbf{D}: \mathbf{D}
+\beta_3|\mathbf{D}  \mathbf{d}|^2 \geqslant 0,
\end{equation}
 if and only if
\begin{equation}\label{beta_inequality_second}
    \beta_1 + 2\beta_2+\beta_3 \geqslant 0, \quad 2\beta_2+\beta_3 \geqslant 0.
\end{equation}    
\end{lemma}

\paragraph{Proof.} 
We assume that $\mathbf{d}=(d_1,d_2)^{T}$ and $\mathbf{D}=(D_{ij})_{2\times 2}$ with $D_{11}+D_{22}=0$ and $D_{12}=D_{21}$.
By rotation invariance, we  assume that $d_1=1$ and $d_2=0$.
Then we have
$$
\begin{aligned}
\beta_1(\mathbf{d}\otimes\mathbf{d}: \mathbf{D})^2+\beta_2 \mathbf{D}: \mathbf{D}+\beta_3|\mathbf{D} \mathbf{d}|^2 
=(\beta_1+2\beta_2+\beta_3)D_{11}^2
+( 2\beta_2+\beta_3)D_{12}^2.
\end{aligned}
$$
Therefore, we get (\ref{beta_inequality_second}) from (\ref{beta_inequality_first}).

On the other hand,
we have
$$
\begin{aligned}
\beta_1 \left(\mathbf{d}\otimes\mathbf{d}: \mathbf{D}\right)^2
=& \;
\beta_1
\left[
\left(d_1^2-d_2^2\right)D_{11}
+
2d_1 d_2 D_{12}
\right]^2,
\\
\beta_2 \left(\mathbf{D}: \mathbf{D} \right)
=&\;
2\beta_2  \left(D_{11}^2+D_{12}^2\right),
\\
\beta_3|\mathbf{D} \mathbf{d}|^2
=& \;
\beta_3
\left(D_{11}^2+D_{12}^2\right).
\end{aligned}
$$
Denoted by
$\mathbf{d}=(\sin\theta, \cos\theta)^T$, 
we have
$$
\beta_1 \left(\mathbf{d}\otimes\mathbf{d}: \mathbf{D}\right)^2
=
\beta_1
\left( \cos 2\theta D_{11}
-
\sin 2\theta D_{12}
\right)^2.
$$
By Cauchy's inequality,
$$
0\leqslant
\left( \cos 2\theta D_{11}
-
\sin 2\theta D_{12}
\right)^2
\leqslant
\left(D_{11}^2+D_{12}^2\right)
\left( \cos^2 2\theta +
\sin^2 2\theta 
\right)
=
\left(D_{11}^2+D_{12}^2\right).
$$
Suppose that  $m=m(x, t)\in [0, 1]$ satisfies 
$$
\left( \cos 2\theta D_{11}
-
\sin 2\theta D_{12}
\right)^2
=
m 
\left(D_{11}^2+D_{12}^2\right).
$$
Then inequality (\ref{beta_inequality_first})  is 
equivalent to satisfying  the following inequality for any 
$ m\in [0, 1]$,
$$
0\leqslant
\left( m(x, t) \beta_1+2\beta_2+\beta_3   \right) \left(D_{11}^2+D_{12}^2\right).
$$
By setting $m=0$ and $m=1$, we can obtain the inequality (\ref{beta_inequality_second}).
$\Box$

\begin{remark}\label{Remark_iso_beta_first}
Suppose that $\mathbf{D}=\mathbf{D}_{n\times n}$ is a symmetric trace free matrix.\\
(1) 
For $n=2$ and $\mathbf{d}:\mathbb{R}^2 \rightarrow \mathbb{S}^1$,
 the  dissipation relation (\ref{beta_inequality_first}) 
is equivalent to (\ref{beta_inequality_second}).
\\
(2)
For $ n=2$  and $\mathbf{d}:\mathbb{R}^2 \rightarrow \mathbb{S}^2$ in \cite{Wangmeng2014cvpde},
 the  dissipation relation (\ref{beta_inequality_first}) 
is equivalent to 
$$
\begin{aligned}
&\beta_2 \geqslant 0,\quad \beta_1 + 2\beta_2 + \beta_3 \geqslant 0,\quad \beta_1 < 0; \\
\text{or}\quad & \beta_2 \geqslant 0, \quad 2\beta_2 + \beta_3 \geqslant 0, \quad \beta_1 \geqslant 0.
\end{aligned}
$$
\\
(3)
For $n=3$ and $\mathbf{d}:\mathbb{R}^3 \rightarrow \mathbb{S}^2$ in \cite{ZZZ2013ARMA},
 the  dissipation relation (\ref{beta_inequality_first}) 
is equivalent to 
$$
\beta_2 \geqslant 0, \quad 2\beta_2 + \beta_3 \geqslant 0, \quad \frac{3}{2}\beta_2 + \beta_3 + \beta_1 \geqslant 0.
$$
\end{remark}

\paragraph{}
The system (\ref{general_Ericksen--Leslie}) takes the form
\begin{equation}\label{Ericksen--Leslie_new_2}
\left\{
\begin{aligned}
&\partial_t \mathbf{v}+\mathbf{v} \cdot \nabla \mathbf{v}
=-\nabla p
+\frac{\gamma}{\operatorname{Re}} \Delta \mathbf{v}
+\frac{1-\gamma}{\operatorname{Re}}\nabla \cdot\left(\sigma_1
+\sigma_2
+\sigma^E\right), \\
&\nabla \cdot \mathbf{v}=0 \\
&\partial_t \mathbf{d}+\mathbf{v} \cdot \nabla \mathbf{d}
=\mu_1\left(\Delta\mathbf{d}+f(\mathbf{d},\nabla\mathbf{d}) \right)
+(\mathbf{I}-\mathbf{d}\otimes\mathbf{d})\left(\boldsymbol{\Omega}  \mathbf{d}
+\mu_2 \mathbf{D} \mathbf{d}\right),
\end{aligned} 
\right.
\end{equation}
where $|\mathbf{d}|=1$ and
$$
\begin{aligned}
& \sigma_1 =\beta_1(\mathbf{D}: \mathbf{d}\otimes\mathbf{d}) \mathbf{d}\otimes\mathbf{d}
+\beta_2 \mathbf{D}
+\frac{\beta_3}{2}\left[\mathbf{d}\otimes(\mathbf{D} \mathbf{d}) 
+(\mathbf{D}\mathbf{d})\otimes\mathbf{d}\right], \\
& \sigma_2
=\frac{1}{2}\left(-1-\mu_2\right) \left[\mathbf{h}-\left(\mathbf{h}\cdot\mathbf{d}\right)\mathbf{d} \right]\otimes\mathbf{d}
+\frac{1}{2}\left(1-\mu_2\right) \mathbf{d}\otimes\left[\mathbf{h}-\left(\mathbf{h}\cdot\mathbf{d}\right)\mathbf{d} \right].
\end{aligned}
$$
To deal with $\sigma_2$,
we need the following coupling relationship between $\sigma_2$ and $\mathbf{h}$.

\begin{lemma}\label{sigma_2_lem}
Suppose that  $(\mathbf{v},\mathbf{d})$ are smooth functions with respect to $x\in Q$,
then it holds that
\begin{equation}\label{sigma_2_first_equation}
\sigma_2:\nabla \mathbf{v}
+\left[( \mathbf{I}-\mathbf{d}\otimes\mathbf{d})(\boldsymbol{\Omega}  \mathbf{d}
+\mu_2 \mathbf{D} \mathbf{d}) \right]\cdot \mathbf{h}
=0.
\end{equation}
\end{lemma}

The proof of  Lemma \ref{sigma_2_lem} can be obtained by Lemma 3.1 in \cite{ZZZ2013ARMA} and Lemma 3.5 in \cite{Wang2024DCDS}.

For isotropic liquid crystal fluid,
we have $k_1=k_2=k_3=1$,
and from (\ref{functional_magnetic_field}),
the free energy functional $E_M$ for a nematic liquid crystal under a magnetic field  is given by
\begin{equation}\label{isotropic_functional_magnetic_field}
E_M
=
\frac{1}{2}
\int_{Q}
\left[|\nabla \mathbf{d}|^2 + |\mathbf{H}|^2-(\mathbf{H}\cdot \mathbf{d})^2\right]
\mathrm{~d}x,
\end{equation}
where $\mathbf{H}=\mathrm{H}(1,0)^T$ is the external magnetic field
and
$\mathrm{H}$ is a constant. 

For the stress tensor $\sigma^E$,
we have
$$
\left(\sigma^E\right)_{i j}
=-(\nabla \mathbf{d} \odot \nabla \mathbf{d})_{i j}
=-\partial_i \mathbf{d}\cdot \partial_j \mathbf{d}.
$$
By (\ref{h_def}) and (\ref{isotropic_functional_magnetic_field}),
the molecular field $\mathbf{h}$ is defined as follows
\begin{equation}\label{h_isotropic_def}
\mathbf{h}
 =-\frac{\delta E_M}{\delta \mathbf{d}}
 =\nabla \cdot \frac{\partial E_M}{\partial \nabla \mathbf{d}}-\frac{\partial E_M}{\partial \mathbf{d}}
 =\Delta \mathbf{d} + (\mathbf{H}\cdot \mathbf{d}) \mathbf{H} .
\end{equation}
Then we have 
$$
\mathbf{h}-(\mathbf{h}\cdot\mathbf{d})\mathbf{d}
=
\Delta \mathbf{d} 
+
|\nabla\mathbf{d}|^2\mathbf{d}
+
(\mathbf{H}\cdot \mathbf{d}) \mathbf{H}
-
(\mathbf{H}\cdot \mathbf{d})^2 \mathbf{d}.
$$
To simplify the form of  $\mathbf{h}-(\mathbf{h}\cdot\mathbf{d})\mathbf{d}$,   
define
$$
f(\mathbf{d}, \nabla\mathbf{d}) 
=|\nabla\mathbf{d}|^2\mathbf{d}
+(\mathbf{H}\cdot\mathbf{d})\mathbf{H} 
-(\mathbf{H}\cdot\mathbf{d})^2 \mathbf{d}.    
$$
We have 
$$
\mathbf{h}-(\mathbf{h}\cdot\mathbf{d})\mathbf{d}
=
\Delta \mathbf{d} 
+
f(\mathbf{d}, \nabla\mathbf{d}) .
$$

From (\ref{Ericksen_Leslie_new_formulation}),
we have
the following system 
about $(\mathbf{v},   \mathbf{d})$ in $Q\times\left[0,\infty\right)$,
where $\mathbf{d}\in\mathbb{S}^1$.
\begin{equation}\label{iso_Ericksen_Leslie_new_formulation}
\left\{
\begin{aligned}
\partial_t \mathbf{v}+\mathbf{v} \cdot \nabla \mathbf{v}
=&-\nabla p
+\frac{\gamma}{\operatorname{Re}} \Delta \mathbf{v}
+\frac{1-\gamma}{\operatorname{Re}}\nabla \cdot\left(
\sigma_1+\sigma_2
+\sigma^E\right),     \\
\nabla \cdot \mathbf{v}=&\;0,    \\
\partial_t \mathbf{d}+\mathbf{v} \cdot \nabla \mathbf{d}
=&\;\mu_1\left(\Delta \mathbf{d} 
+
f(\mathbf{d}, \nabla\mathbf{d})  \right)
+(\mathbf{I}-\mathbf{d}\otimes\mathbf{d})\left(\boldsymbol{\Omega}  \mathbf{d}
+\mu_2 \mathbf{D} \mathbf{d}\right),
\end{aligned}
\right.
\end{equation}
where
$$
\begin{aligned}
& \sigma_1=\;\beta_1(\mathbf{D}: \mathbf{d}\otimes\mathbf{d}) \mathbf{d}\otimes\mathbf{d}
+\beta_2 \mathbf{D}
+\frac{\beta_3}{2}(\mathbf{d}\otimes(\mathbf{D} \mathbf{d}) 
+(\mathbf{D}\mathbf{d})\otimes\mathbf{d}),  \\
& \sigma_2
=\;\frac{1}{2}\left(-1-\mu_2\right) (\Delta\mathbf{d}+f(\mathbf{d}, \nabla\mathbf{d}) )\otimes\mathbf{d}
+\frac{1}{2}\left(1-\mu_2\right) \mathbf{d}\otimes(\Delta\mathbf{d}+f(\mathbf{d}, \nabla\mathbf{d}) ).
\end{aligned}
$$
Then we have the following basic energy law.

\begin{proposition}
Suppose that  $(\mathbf{v},\mathbf{d})$ are smooth solutions to (\ref{iso_Ericksen_Leslie_new_formulation}) with boundary conditions (\ref{boundary_condition_1}),
then it holds that
\begin{equation}\label{iso_basic_energy_law}
\begin{aligned}
&\frac{1}{2} \frac{\mathrm{d}}{\mathrm{d} t} 
\int_Q
|\mathbf{v}|^2
+\frac{1-\gamma}{\operatorname{Re}}
(|\nabla \mathbf{d}|^2 +|\mathbf{H} |^2 -(\mathbf{H}\cdot\mathbf{d})^2 ) 
\mathrm{~d} x\\
= \;&  
-\int_Q
\frac{\gamma}{\operatorname{Re}}
|\nabla \mathbf{v}|^2
+
\frac{(1-\gamma)\mu_1}{\operatorname{Re}}
\left[
|\mathbf{h}|^2
-(\mathbf{h} \cdot \mathbf{d})^2 \right]
\mathrm{d} x \\
&  
-\int_Q
\frac{1-\gamma}{\operatorname{Re}}
\left[\beta_1(\mathbf{d}\otimes\mathbf{d}: \mathbf{D})^2
+\beta_2 \mathbf{D}: \mathbf{D} 
+\beta_3|\mathbf{D}\mathbf{d}|^2
\right]
\mathrm{d} x.
\end{aligned}
\end{equation}
\end{proposition}

\paragraph{Proof.}
The proof of basic energy law is inspired by \cite{ZZZ2013ARMA,Wang2024DCDS}.

Multiply $\mathbf{v}$ on both sides of the first equation of \eqref{Ericksen_Leslie_new_formulation} and integral it in $Q$
\begin{equation}
    -\frac{1}{2} \frac{\mathrm{d}}{\mathrm{d} t}\int_Q|\mathbf{v}|^2 \mathrm{d} x 
=\int_Q \frac{ \gamma}{\operatorname{Re}}|\nabla \mathbf{v}|^2+\frac{1-\gamma}{\operatorname{Re}}[\left(\sigma_1+\sigma_2\right): \nabla \mathbf{v}+\sigma^E: \nabla \mathbf{v}] \mathrm{~d} x.
\end{equation}
For the Ericksen stress term, we have
$$
\int_Q  
\left[\sigma^E: \nabla \mathbf{v}-(\mathbf{v} \cdot \nabla \mathbf{d}) \cdot \Delta \mathbf{d} 
\right] \mathrm{~d}x
=0.
$$
From Lemma \ref{beta_lemma} and 
$ \sigma_1: \nabla \mathbf{v}=\sigma_1:\mathbf{D}$,
we have
$$
\int_Q
\sigma_1: \nabla \mathbf{v}
\mathrm{~d} x
=
\int_Q
\frac{1-\gamma}{\operatorname{Re}}
\left[\beta_1(\mathbf{d}\otimes\mathbf{d}: \mathbf{D})^2
+\beta_2 \mathbf{D}: \mathbf{D} 
+\beta_3|\mathbf{D}\mathbf{d}|^2
\right]
\mathrm{~d} x.
$$

Multiply $\frac{1-\gamma}{\operatorname{Re}}\mathbf{h}$ on both sides of the third equation of \eqref{Ericksen_Leslie_new_formulation} and integral it in $Q$
\begin{equation}
\begin{aligned}
   & \frac{1-\gamma}{\operatorname{Re}}\int_Q 
( \partial_t \mathbf{d}+\mathbf{v} \cdot \nabla \mathbf{d})\cdot  \mathbf{h} \mathrm{~d} x \\
=&\; 
\frac{1-\gamma}{\operatorname{Re}}\int_Q
\mu_1|\mathbf{h}|^2-\mu_1(\mathbf{h}\cdot \mathbf{d})^2 
+[(\mathbf{I}-\mathbf{d}\otimes\mathbf{d})(\boldsymbol{\Omega}  \mathbf{d}+\mu_2 \mathbf{D}\mathbf{d})]\cdot \mathbf{h} \mathrm{~d} x.
\end{aligned}
\end{equation}
where
$$
\begin{aligned}
\frac{1-\gamma}{\operatorname{Re}}\int_Q 
( \partial_t \mathbf{d}+\mathbf{v} \cdot \nabla \mathbf{d})\cdot  \mathbf{h} \mathrm{~d} x 
=& -\frac{\mathrm{d}}{\mathrm{dt}}\frac{1-\gamma}{\operatorname{Re}}\int_Q \frac{1}{2}|\nabla \mathbf{d}|^2+\frac{1}{2}(|\mathbf{H} |^2 -(\mathbf{H}\cdot\mathbf{d})^2 )
	\mathrm{~d} x\\
   & +\frac{1-\gamma}{\operatorname{Re}}\int_Q 
	(\mathbf{v} \cdot \nabla \mathbf{d})\cdot  \Delta\mathbf{d} 
	+\mathbf{v} \cdot \nabla [\frac{(\mathbf{H}\cdot \mathbf{d})^2 }{2}]
	\mathrm{~d} x.
\end{aligned}
$$
It is easy to verify that the boundary terms vanish by the boundary conditions (\ref{boundary_condition_1}). 
From Lemma \ref{sigma_2_lem} about $\sigma_2$ and $\mathbf{h}$,
we can obtain the basic energy law.    $\Box$

From basic energy law   (Theorem \ref{iso_basic_energy_law}),
we can deduce that
$$
\begin{aligned}
& \mathbf{v} \in L^{\infty}(0, T ; H) \cap L^2(0, T ; V), \\
& \mathbf{d} \in L^{\infty}\left(0, T ; H_{sp}^1 \right) \cap L^2_{loc}\left(0, T ;H_{sp}^2\right).  
\end{aligned}
$$

\subsection{\texorpdfstring{System for $\left(\mathbf{v}, \theta\right)$}{System for} 
}

\setcounter{equation}{0}
\numberwithin{equation}{section}

\paragraph{\texorpdfstring{System for $\left(\mathbf{v}, \theta\right)$ }{System for_v_d} }
Suppose that \(\mathbf{d} = (\sin\theta, \cos\theta)^T\),
then we can  derive the equation for \(\theta\),
which is important for the well-posedness ane long-time behavior of  the general Ericksen--Leslie system.

At first, 
we need to prove that  \(\theta_0\) corresponding to \(\mathbf{d}_0 = (\sin\theta_0, \cos\theta_0)^T\) is single-valued,
which implies that if there exists  \(x_0 \in Q\) such that \(\theta_0(x_0)\) is determined, then   \(\theta_0(x)\) is determined for all \(x \in Q\).

\begin{theorem}\label{thm_theta_single_value}
Suppose that $\mathbf{d}_0 \in H^2_{sp}$ and 
\(\mathbf{d}_0 = (\sin\theta_0, \cos\theta_0)^T\),
then $\theta_0$ is single-valued.    
\end{theorem}

\paragraph{Proof.}
Without loss of generality, we only consider  the interior of \( Q = [0, 1] \times [0, 1] \).

For any \( x_1 \in \partial Q \), replace the function \( \theta_0 \) and the region \( Q \) with
\[
\widetilde{\theta}_0 (x) = \theta\left(x - x_1 + \frac{1}{2}(\mathbf{e}_1 + \mathbf{e}_2)\right), \quad \widetilde{Q} = \left\{ y \;\bigg|\; y = x - x_1 + \frac{1}{2}(\mathbf{e}_1 + \mathbf{e}_2), \; x \in Q \right\}.
\]
After the translation, the function \( \widetilde{\theta}_0 \) still satisfies the boundary conditions in the region \( \widetilde{Q} \), and \( x_1 \) is in the interior of \( \widetilde{Q} \).

Since $Q\in \mathbb{R}^2$ is a simply connected domain, every closed curve can be continuously deformed to any point within the region enclosed by the closed curve.
Suppose that there exist two curves $r_1, r_2$ starting from $x_0$ to $y_0 \in Q$ such that $\theta_0(y_0)$ differs on $r_1, r_2$, then a smooth closed curve $C_{r_1, r_2}$ can be formed from $r_1, r_2$. 
Without loss of generality, assume that $C_{r_1, r_2}$ is a smooth simple closed curve.

The notation \( \mathrm{deg}\big(\mathbf{d}_0, C_{r_1, r_2}\big) \) represents the number of counterclockwise rotations that \( \mathbf{d}_0 \) makes around the closed curve \( C_{r_1, r_2} \). The topological degree\cite{Zhang2005NonlinearAnalysis}  \( \mathrm{deg}\big(\mathbf{d}, C_{r_1, r_2}\big) \in \mathbb{Z} \) is a homotopy invariant.
The region enclosed by the closed curve \(C_{r_1, r_2}\)
is denoted by $R\big(C_{r_1, r_2}\big)$, 
the diameter of which is denoted by $l<1$.

Denoted by \(\mathbf{d}_0 (x)= \left( d_1 (x), d_2 (x) \right) \in \mathbb{S}^1\), and let \( x_2 \in R\big(C_{r_1, r_2}\big) \). 
When \(d_2=0\), define additionally
\[
\arctan(\infty) = \frac{\pi}{2}, \quad \arctan(-\infty) = -\frac{\pi}{2}.
\]
Thus, we have
\[
\theta_0 = \arg\left(\mathbf{d}_0(x) - x_2 \right) = \arctan \left( \frac{d_1}{d_2} \right) + k\pi.
\]
Here, the integer \(k\) can vary to ensure that the function \(\theta_0\) is continuous.

By trace interpolation inequalities from Lemma (\ref{lem_inequalities}), if \(\mathbf{d}_0(x) \in H^2\), then \(\mathbf{d}_0(x) \in H^{\frac{3}{2}}\big( C_{r_1, r_2} \big)\). Therefore, there exist  constants \(K\), which is independent of   
  the closed curve   \(C_{r_1, r_2} \), and $s\in (0,1)$ such that 
    \[
    \|\nabla\mathbf{d}_0\|^2_{H^1\big( C_{r_1, r_2} \big)} 
    \leqslant 
    K l^s
    \|\mathbf{d}_0\|^2_{H^{\frac{3}{2}}\big( R\big(C_{r_1, r_2}\big) \big)} 
    \leqslant 
    K l^s
    \|\mathbf{d}_0\|^2_{H^{2}(Q)}.
    \]

For any point \(x \in C_{r_1, r_2}\), the vector \(\mathbf{r}(x)\) is the unit tangent vector to the curve \(C_{r_1, r_2}\) at point \(x\), and \(\mathbf{r}\) rotates counterclockwise along \(C_{r_1, r_2}\). Let \(\left|C_{r_1, r_2}\right|\) denote the length of the curve \(C_{r_1, r_2}\), we obtain
\[
\begin{aligned}
\left|
2\pi \mathrm{deg}\big(\mathbf{d}_0, C_{r_1, r_2}\big) 
\right|
&= \left| \int_{C_{r_1, r_2}} \nabla \left( \arctan \frac{d_1}{d_2} \right) \cdot \mathbf{r} \mathrm{~d}s \right| \\
&= \left| \int_{C_{r_1, r_2}} \frac{d_2 \left(\mathbf{r} \cdot \nabla d_1\right) - d_1 \left(\mathbf{r} \cdot \nabla d_2\right)}{d_1^2 + d_2^2} \mathrm{~d}s \right|\\
&\leqslant
\left| {C_{r_1, r_2}} \right|^{\frac{1}{2}}
\left| \int_{C_{r_1, r_2}} \left|\nabla\mathbf{d}_0\right|^2 \mathrm{~d}s \right|^{\frac{1}{2}}.
\end{aligned}
\]
where \(d_1^2 + d_2^2 = 1\).

 If there is a singularity of \(\theta_0\) within the region enclosed by the simple closed curve \(C_{r_1, r_2}\), then there exists a smooth curve \(C_{r_1, r_2}'\) surrounding one of the singularities such that 
\[
    \mathrm{deg}\big(\mathbf{d}_0, C_{r_1, r_2}'\big) \neq 0, \quad \left|C_{r_1, r_2}'\right| \leqslant \frac{1}{ K l^s
    \|\mathbf{d}_0\|^2_{H^{2}(Q)}}. 
 \]
We obtain
\[
    (2\pi)^2
    \leqslant 
     \left| 2\pi\mathrm{deg}\big(\mathbf{d}_0, C_{r_1, r_2}'\big) \right|^2 
     \leqslant 
     \left|C_{r_1, r_2}'\right|
     \left| \int_{C_{r_1, r_2}'} \left|\nabla\mathbf{d}_0\right|^2 \mathrm{~d}s \right|
     \leqslant 
      K l^s
    \|\mathbf{d}_0\|^2_{H^{2}(Q)}
    \left|C_{r_1, r_2}'\right|
    \leqslant
    1.
\]
Therefore, no such \(y_0 \in Q\) exists, and \(\theta_0\) is single-valued.  $\Box$

\paragraph{}
Let $\mathbf{d} = (\sin\theta, \cos\theta)^T$. 
By Theorem \ref{thm_theta_single_value},
for a fixed time $t>0$, if $\mathbf{d}(x, t) \in H^{2}_{sp}$, then  $\theta(x, t)$ is single-valued.
Therefore,
we can derive the general system for $(\mathbf{v},\theta)$.

Suppose that
$\mathbf{d^{\perp}}=(\cos\theta,-\sin\theta)^{T}$,
then
\begin{equation}
\begin{aligned}
\Delta\mathbf{d}+f(\mathbf{d},\nabla\mathbf{d}) 
=
&\;
\mathbf{h}+
|\nabla \mathbf{d}|^2 \mathbf{d}
-(\mathbf{H}\cdot\mathbf{d})^2 \mathbf{d}\\
=  
&\;
\left(\Delta \theta + \frac{\mathrm{H}^2}{2} \sin 2\theta  \right) \mathbf{d^{\perp}}.
\end{aligned}
\end{equation}
Furthermore,
we can obtain the following identities,
$$
\begin{aligned}
|\mathbf{h}|^2-|\mathbf{h}\cdot\mathbf{d} |^2 
=
&\;
 |\mathbf{h}-(\mathbf{h}\cdot\mathbf{d})\cdot\mathbf{d} |^2
\\
=&\;
|\mathbf{h}+|\nabla \mathbf{d}|^2 \mathbf{d}-(\mathbf{H}\cdot\mathbf{d})^2 \mathbf{d}|^2 \\
=&\;
\left|\Delta\theta+\frac{\mathrm{H}^2}{2}\sin2\theta\right|^2 \\
=&\;|\Delta\mathbf{d}+f(\mathbf{d},\nabla\mathbf{d}) |^2 .  
\end{aligned}
$$
Thanks to
$
\mathbf{I}-\mathbf{d}\otimes\mathbf{d
}=\mathbf{d}^{\perp}\otimes\mathbf{d}^{\perp},
$
then we have
$$
\left(\mathbf{d}^{\perp}\otimes\mathbf{d}^{\perp}\right)
\left(\boldsymbol{\Omega}  \mathbf{d}
+\mu_2 \mathbf{D} \mathbf{d}\right)
=
\mathbf{d}^{\perp}
\left( 
\boldsymbol{\Omega}:\mathbf{d}^{\perp}\otimes\mathbf{ d} 
+\mu_2 \mathbf{D}:\mathbf{d}^{\perp}\otimes\mathbf{d}   \right).
$$

The system \eqref{iso_Ericksen_Leslie_new_formulation} can be written as follows,
\begin{equation}\label{iso_theta_Ericksen_Leslie_new_formulation}
\left\{\begin{aligned}
\partial_t \mathbf{v}+\mathbf{v} \cdot \nabla \mathbf{v}
=&
-\nabla p
+\frac{\gamma}{\operatorname{Re}} \Delta \mathbf{v}
+\frac{1-\gamma}{\operatorname{Re}} \nabla \cdot \sigma_1
+
\frac{1-\gamma}{\operatorname{Re}}  
\left[
 -\Delta \theta \cdot \nabla \theta
-\nabla\left(\frac{|\nabla \theta|^2}{2}\right) \right] \\ 
&+
\frac{1-\gamma}{\operatorname{Re}}  \nabla \cdot\left[
\frac{1}{2}\left(-1-\mu_2\right) 
\mathbf{ d^{\perp}}\otimes\mathbf{d}
\left(\Delta\theta+\frac{\mathrm{H}^2}{2}\sin 2\theta\right) 
\right]\\
&+\frac{1-\gamma}{\operatorname{Re}}  
\nabla \cdot\left[
\frac{1}{2}\left(1-\mu_2\right) \mathbf{d}\otimes\mathbf{d^{\perp}} 
\left(\Delta\theta+\frac{\mathrm{H}^2}{2}\sin 2\theta\right) 
\right], \\
\nabla \cdot \mathbf{v}=&0,\\
\partial_t \theta
+\mathbf{v} \cdot \nabla \theta
=&
\mu_1\left(\Delta \theta + \frac{\mathrm{H}^2}{2} \sin 2\theta  \right) 
+\left(\boldsymbol{\Omega}:\mathbf{d}^{\perp}\otimes\mathbf{ d} 
+\mu_2 \mathbf{D}:\mathbf{d}^{\perp}\otimes\mathbf{d} \right).
\end{aligned}\right.
\end{equation}	

From $\mathbf{d}_0 = (\sin \theta_0,     \cos \theta_0)^T $ and initial conditions (\ref{initial_condition_1}) about $(\mathbf{v},\mathbf{d})$,
we have the following initial conditions about 
$(\mathbf{v},\theta)$
\begin{equation}\label{iso_theta_initial_condition}
\mathbf{v}|_{t=0} = \mathbf{v}_0,    \quad 
\nabla\cdot \mathbf{v}_0 = 0,    \quad
\theta|_{t=0} = \theta_0,    \quad
 x\in Q,
\end{equation}
where constants $a_i\; (i=1,    2)$ are integers,
depending on
$$
\theta_0(x+\mathbf{e}_i) = \theta_0(x) + a_i \pi,    \quad
x\in \partial Q.
$$
By boundary conditions (\ref{boundary_condition_1}),
we have the following initial conditions about 
$(\mathbf{v},\theta)$
\begin{equation}\label{iso_theta_boundary_condition}
\mathbf{v}(x+\mathbf{e}_i,    t) = \mathbf{v}\left(x,    t\right),    \quad 
\theta(x+\mathbf{e}_i,    t) = \theta\left(x,    t\right) + a_i \pi,      \quad 
\left(x,    t\right)\in \partial Q \times  \left[0,    \infty \right).    
\end{equation}

Moreover, we have
$$
\begin{aligned}
&\mathbf{d^{\perp} }(\boldsymbol{\Omega}:\mathbf{d}^{\perp}\otimes\mathbf{ d} )
=(\mathbf{I}-\mathbf{d}\otimes\mathbf{d}) (\boldsymbol{\Omega}\cdot\mathbf{d} )
=\boldsymbol{\Omega}\cdot\mathbf{d} -\mathbf{d}\otimes\mathbf{d}(\boldsymbol{\Omega}\cdot\mathbf{d})
=\boldsymbol{\Omega}\cdot\mathbf{d},\\
&\mathbf{d^{\perp} }(\mathbf{D}:\mathbf{d}^{\perp}\otimes\mathbf{d})
   =(\mathbf{I}-\mathbf{d}\otimes\mathbf{d}) (\mathbf{D}\cdot\mathbf{d} )
   =\mathbf{D}\cdot\mathbf{d} -\mathbf{d}\otimes\mathbf{d}(\mathbf{D}\cdot\mathbf{d})
   =\mathbf{D}\cdot\mathbf{d} -(\mathbf{D}:\mathbf{ d}\otimes\mathbf{d})\mathbf{d}.
\end{aligned}
$$
Because
$
\mathbf{D}=\dfrac{1}{2}\left[\nabla \mathbf{v}+(\nabla \mathbf{v})^T\right], \;\boldsymbol{\Omega}=\dfrac{1}{2}\left[\nabla \mathbf{v}-(\nabla \mathbf{v})^T\right]
$,
then 
\begin{equation}
\boldsymbol{\Omega}\mathbf{d}
       +\mu_2 \mathbf{D}\mathbf{d}
= \frac{1}{2}(1+\mu_2)(\nabla\mathbf{v}\cdot\mathbf{d})
+\frac{1}{2}(-1+\mu_2)((\nabla\mathbf{v})^{T}\cdot\mathbf{d}).
\end{equation}
By comparing the coefficients of  $\sigma_2$ and $\boldsymbol{\Omega}\mathbf{d}
+\mu_2 \mathbf{D}\mathbf{d}$, 
which are $(1+\mu_2)$ and $(-1+\mu_2)$, 
we can obtaion the
 coupling relationship.

\section{Steady-state problem}

\setcounter{equation}{0}
\numberwithin{equation}{section}

If $\mathbf{v}=0$ and $\theta=\theta(x)$,
then the system (\ref{iso_Ericksen_Leslie_new_formulation})
with initial conditions (\ref{iso_theta_initial_condition}) and 
boundary conditions (\ref{iso_theta_boundary_condition})
can be reduced to the following steady-state problem including the elliptic sine-Gorden  equation,
\begin{equation}\label{elliptic_first}
\left\{
\begin{aligned}
-\Delta \theta
=& \;
\frac{\mathrm{H}^2}{2}
\sin 2\theta
,
\quad
&& x\in Q,
\\
\theta(x+\mathbf{e}_i) 
=& \;
\theta(x)+a_i\pi,
\quad
&& x\in \partial Q.
\end{aligned}\right.
\end{equation}
Using the values of $a_i$ 
as classification conditions, 
we distinguish between periodic and modulo-$\pi$ periodic boundary conditions to identify different steady-state solutions. 
We need to prove the existence of solutions to
(\ref{elliptic_first})
and discuss the existence of non-trivial solutions based on the critical value $\lambda_2$ in (\ref{lambda_2_def}). 
The regularity of steady-state solutions can be obtained through the bootstrap method.

\subsection{Periodic boundary conditions}
We consider the  following elliptic equation under periodic boundary conditions
\begin{equation}\label{elliptic_periodic_boundary}
\left\{
\begin{aligned}
-\Delta \theta(x)
=& \;
\frac{\mathrm{H}^2}{2}
\sin 2\theta(x)
,
\quad
&& x\in Q,
\\
\theta(x+\mathbf{e}_i) 
=& \;
\theta(x),
\quad
&& x\in \partial Q.
\end{aligned}\right.
\end{equation}

The existence of solutions for the elliptic Sine--Gordon equation under Neumann boundary conditions was introduced in paper \cite{Ni2004}.

Suppose that 
$\psi:Q\rightarrow\mathbb{R}$ and
$\psi(x)=2\theta(x)$,
then we have
\begin{equation}\label{iso_psi_periodic_problem}
 \left\{\begin{aligned}
-\Delta \psi(x) & = \;\mathrm{H}^2 \sin \psi(x),   && x\in  Q, \\
\psi(x+\mathbf{e}_i) & = \;\psi(x),   && x\in \partial Q.
\end{aligned}\right.   
\end{equation}
The equation in (\ref{iso_psi_periodic_problem}) corresponds to the Euler--Lagrange equation of the following energy functional
$$
   I(\psi)=\int_Q \frac{1}{2} |\nabla\psi|^2 +\mathrm{H}^2(-1+\cos \psi) \mathrm{~d}x. 
$$

For $\mathrm{H}^2=0$, it is obvious that $\psi=C$ is the solution to  (\ref{iso_psi_periodic_problem}) for any constant $C$.
The following Lemma \ref{iso_elliptic_periodic_existence_first} gives the existence of (\ref{iso_psi_periodic_problem}),
when $\mathrm{H}^2>0$.

\begin{lemma}\label{iso_elliptic_periodic_existence_first}
Suppose that $\mathrm{H}^2>0$ and  $a_1=a_2=0$,
then the energy functional
$I(\psi)$ has a local minimum at  $\psi=\pi$.
Moreover,
for every integer $k$,
$\psi=(2 k-1) \pi$
are also local minimum points.
\end{lemma}

The proof can be seen in Lemma 2.2 of  \cite{Ni2004}, 
which is under Neumann boundary conditions and has a similar progress to the  periodic boundary conditions.

\begin{lemma}\label{lem_iso_lambda_2_mountain_pass_type_solution}
If $\mathrm{H}^2>\lambda_2$ and  $a_1=a_2=0$, 
then energy functional
$I(\psi)$ has a nonconstant critical point of mountain-pass type in $H^1_p$. 
Moreover, 
$\psi(x)+2 k \pi$ 
are also solutions of (\ref{iso_psi_periodic_problem}) for any integer $k$.  
\end{lemma}

\paragraph{Proof.}
The proof is inspired by Theorem 2.3 in  \cite{Ni2004}.
It can be verified that  $I(\psi)$  satisfies the Palais--Smale condition. By Lemma \ref{iso_elliptic_periodic_existence_first} and Mountain Pass Lemma, 
energy functional $I(\psi)$ has at least one  critical point $\omega_0 \in H^1_p$ such that
$$
 I\left(\omega_0 \right)=\inf _{\gamma \in \Gamma} \sup _{\psi \in \gamma([0, 1])} I(\psi)\geqslant\delta_0>0, 
$$
where $\Gamma=\left\{\gamma \in C\left(\mathbb{R},  H^1_p \right) \mid \gamma(0)=\pi,  \gamma(1)=-\pi\right\}$.

We use  $H^{-1}$ to denote the dual space of $H^1_p$. 
With a similar progress, we have
\begin{equation}\label{iso_morse_1}
<I^{\prime \prime}(0) \psi_2,  \psi_2>_{H^{-1}, H^1_p}=\left(\lambda_2-
\mathrm{H}^2 \right) \int_Q \psi_2^2 \mathrm{~d} x <0.
\end{equation}
Thanks to  $\phi_1=1$ and $\lambda_1=0$, 
we can obtain
\begin{equation}\label{iso_morse_2}
<I^{\prime \prime}(0) \psi_1,  \psi_1>_{H^{-1}, H^1_p}=-\mathrm{H}^2 |Q|  = -\mathrm{H}^2  <0 .  
\end{equation}

Since  $\psi_1$ and $\psi_2$
are orthogonal in $H^1_p$, 
by (\ref{iso_morse_1}) and (\ref{iso_morse_2}),
the Morse index of $I^{\prime \prime}(0) $ is at least 2.
Because  $\omega_0$
is a mountain-pass type critical point of $I$, the Morse index of $I$ is at most 1. Therefore,  $\omega_0 \neq 0$, which implies that $I$ has a critical point $\omega_0$.  
This completes the proof of Lemma \ref{lem_iso_lambda_2_mountain_pass_type_solution}.
$\Box$

\begin{remark}\label{remark_3_1}
If $\psi(x)$ is a solution of (\ref{iso_psi_periodic_problem}), 
for each constant vector $w \in \mathbb{R}^2$ and any $k\in \mathbb{Z}$,
the function $\psi(x+w)+2k\pi$ can be another solution of (\ref{iso_psi_periodic_problem}).
Therefore, the uniqueness of $\psi (x)$ in (\ref{iso_psi_periodic_problem}) means that 
if $\psi_1(x)$ and $\psi_2(x)$ are two solutions satisfying (\ref{iso_psi_periodic_problem}), there must have a constant vector $w \in \mathbb{R}^2$ and 
an integer constant $k$
such that 
$\psi_1(x)=\psi_2(x+w)+2k\pi$.

In fact, the critical point $-\omega_0$ can be another solution to the first equation in (\ref{iso_psi_periodic_problem}) and  under Neumann boundary condition, which can be seen in Theorem 2.3 of \cite{Ni2004}.
However,
for the nonconstant solution $\omega_0 (x)$ to
(\ref{iso_psi_periodic_problem}) under 
periodic boundary conditions,
we wonder whether there exists a constant vector 
$w \in \mathbb{R}^2$ and a constant $k\in \mathbb{Z}$ 
such that
\begin{equation}\label{Omega_0_steady}
\omega_0(x)= -\omega_0(x+w) +2k\pi,
\end{equation}
which means solutions $\omega_0$ and $-\omega_0$ to
(\ref{iso_psi_periodic_problem})
can be regarded as a single solution.
Let $\omega_1(x)=\omega_0(x-\frac{1}{2}w)+2k_0\pi$,
 we have
$$
\omega_1(x-\frac{1}{2}w)
+\omega_1(x+\frac{1}{2}w)
=
2(k-2 k_0)\pi.
$$
Therefore, 
the constant $k$ cannot be omitted, and when $k$ is an even integer, it is only necessary to prove that $\omega_0$ is a centrosymmetric function.
We intend to tackle this issue elsewhere in the future.
\end{remark}

From Lemma \ref{lem_iso_lambda_2_mountain_pass_type_solution},
it can be seen that if  $\mathrm{H}^2>\lambda_2$, 
then there exists at least one nonconstant
solution to (\ref{iso_psi_periodic_problem}).
We can obtain that if $\mathrm{H}^2 \leqslant \lambda_2$,
the solutions to (\ref{iso_psi_periodic_problem}) are constants.

\begin{theorem}\label{Theorem3.2}
Suppose that $\mathrm{H}^2\leqslant \lambda_2$  and 
 $a_1=a_2=0$, 
 then the solutions to (\ref{iso_psi_periodic_problem}) are constants.
\end{theorem}

\paragraph{Proof.}
If  $\mathrm{H} =0 $ and 
$\psi=\psi(x)$ is a nonconstant solution, then 
$$
    0
    <
    \int_Q |\nabla\psi|^2 \mathrm{~d}x
    =0.
$$
Therefore, when $\mathrm{H} =0$,
the solutions to (\ref{iso_psi_periodic_problem}) are constants.

If $0<\mathrm{H}^2 < \lambda_2$ and 
$\psi=\psi(x)$ is a nonconstant solution,
\begin{equation}
    0
    <
    \int_Q |\nabla\psi|^2 \mathrm{~d}x
    =\mathrm{H}^2 \int_Q \psi\sin\psi \mathrm{~d}x.
\end{equation}
Since $\mathrm{H}^2\int_Q \sin\psi \mathrm{~d}x= -\int_Q \Delta\psi \mathrm{~d}x=0$,
then
$$
\begin{aligned}
    \mathrm{H}^2 \int_Q \psi\sin\psi \mathrm{~d}x
    =
    \mathrm{H}^2 \int_Q \left(\psi-\int_Q \psi \mathrm{~d}x  \right)\left(\sin\psi- \sin\left(\int_Q \psi \mathrm{~d}x\right)\right) \mathrm{~d}x.
\end{aligned}
$$
Thanks to the following estimate with a trigonometric identity,
we have
\begin{equation}\label{lemma3.2_inequ_1}
\begin{aligned}
\int_Q  \left(\sin\psi_1- \sin\psi_2\right) 
\left(\psi_1- \psi_2 \right) \mathrm{~d}x
\leqslant  \int_Q \left(\psi_1- \psi_2 \right)^2 \mathrm{~d}x    .
\end{aligned}
\end{equation}
Let $\psi_1=\psi$ and $\psi_2=\int_Q \psi \mathrm{~d}x$,
we have
\begin{equation}\label{lemma4.2_inequ_2}
\begin{aligned}
    \mathrm{H}^2 \int_Q \psi\sin\psi \mathrm{~d}x
    \leqslant
    \mathrm{H}^2 \int_Q \left(\psi-\int_Q \psi \mathrm{~d}x \right)^2 \mathrm{~d}x
    <
    \lambda_2 \int_Q  \left(\psi-\int_Q \psi \mathrm{~d}x \right)^2 \mathrm{~d}x.
\end{aligned}    
\end{equation}
Since $\lambda_2$ is the principal eigenvalue  (\ref{lambda_2_def}),
then $\psi \equiv \int_Q \psi \mathrm{~d}x $ is a constant.

If $\mathrm{H}^2=\lambda_2$ and 
$\psi=\psi(x)$ is a nonconstant solution, 
consider the conditions for equality in the inequality (\ref{lemma3.2_inequ_1})
$$
\begin{aligned}
\int_Q  \left(\sin\psi_1- \sin\psi_2\right) 
\left(\psi_1- \psi_2 \right) \mathrm{~d}x
=
\int_Q \left(\psi_1- \psi_2 \right)^2 \mathrm{~d}x,
\end{aligned}
$$
where
$\psi_1=\psi$ and $\psi_2=\int_Q \psi \mathrm{~d}x$.
Then we have
$$
    2\sin\left(\frac{\psi_1- \psi_2 }{ 2}\right)
    =
    \psi_1- \psi_2,
    \quad
    \cos\left(\frac{\psi_1+ \psi_2 }{ 2}\right)=1.
$$
we have $(\psi_1+ \psi_2)$ is a constant,
which
illustrates $\psi=\psi(x)$ is a constant.

To summarize,
the solutions to (\ref{iso_psi_periodic_problem}) are constants.
If $0<H^2\leqslant$,
we can obtain that
for every integer $k$,
the solutions to (\ref{iso_psi_periodic_problem}) are 
$\psi = k\pi$. 
$\Box$

\begin{remark}
When  the elliptic sine-Gorden  equation is under Dirichlet boundary condition in \cite{CKY2018SIAM,Ni2004}
\begin{equation}\label{Dirichlet_elliptic_iso}
\left\{
\begin{aligned}
-\Delta \theta
=& \;
\frac{\mathrm{H}^2}{2}
\sin 2\theta
,
\quad
&& x\in \Omega,
\\
\theta
=& \;
0,
\quad
&& x\in \partial \Omega.
\end{aligned}\right.
\end{equation}
where $\Omega$ is a smooth bounded domain,
it is really important that there exist a refined maximum principle for $\theta$, which is also helpful for discussing the  uniqueness of the solution to  (\ref{Dirichlet_elliptic_iso}) in a general sense.

However,
inspired by Theorem 2.4 in  \cite{Ni2004},
for  any integer $k$, 
problem (\ref{iso_psi_periodic_problem}) does not have any nonconstant solution $\psi$ such that
\begin{equation}\label{Maximum Principle}
 k \pi<\psi(x)<( k+1) \pi.    
\end{equation}
Suppose that $0<\psi_0<\pi$ is a nonconstant solution to (\ref{iso_psi_periodic_problem}) 
and  $\psi=\pi-\psi_0$,
then we have
$$
-\int_Q |\nabla\psi |^2 \mathrm{~d}x
=
\mathrm{H}^2 \int_Q \psi \sin\psi \mathrm{~d}x \geqslant 0.
$$
Since  $\int_Q |\nabla\psi |^2 \mathrm{~d}x=\int_Q \psi \sin\psi \mathrm{~d}x=0$,
we have  $\psi=0$ or $\psi=\pi$,  which is a contradiction.
Since $\psi_1=-\psi_0$  is also the solution to (\ref{iso_psi_periodic_problem}),
then if  $-\pi<\psi_1<0$ is a nonconstant solution to (\ref{iso_psi_periodic_problem}), which means $0<\psi_0<\pi$,
we can obtain that $\psi_1=0$ or $\psi_1=-\pi$.
Therefore,
we can get (\ref{Maximum Principle}).
\end{remark}

\subsection{\texorpdfstring{Modulo-$\pi$ periodic boundary conditions}{Modulo pi-periodic boundary conditions}}

We consider the  following elliptic equation under 
modulo-$\pi$ periodic boundary conditions,
when \( |a_1|+|a_2|>0  \),
\begin{equation}\label{elliptic_modulo_periodic_boundary}
\left\{
\begin{aligned}
-\Delta \theta(x)
=& \;
\frac{\mathrm{H}^2}{2}
\sin 2\theta(x)
,
\quad
&& x\in Q,
\\
\theta(x+\mathbf{e}_i) 
=& \;
\theta(x) +a_i\pi,
\quad
&& x\in \partial Q.
\end{aligned}\right.
\end{equation}
Suppose that 
$\psi:Q\rightarrow\mathbb{R}$ and
$\psi(x)=2\theta(x)$,
then we have
\begin{equation}\label{iso_modulo_periodic_problem}
 \left\{\begin{aligned}
-\Delta \psi (x) = & \; \mathrm{H}^2 \sin \psi(x), \quad    && x \in Q, \\
\psi\left(x+\mathbf{e}_i\right) =  & \; \psi(x)+2 a_i\pi,  \quad   && x \in \partial Q.
\end{aligned}\right.   
\end{equation}

We utilize the direct methods  \cite{Struwe2008variational}
 to prove the existence of the  problem (\ref{iso_modulo_periodic_problem}).

\begin{theorem}\label{Theorem_3.3}
Suppose that $|a_1|+|a_2|>0$, then there exists at least one solution to (\ref{iso_modulo_periodic_problem}).
\end{theorem}

\paragraph{Proof.}
Assume that $\phi_0:Q\rightarrow\mathbb{R}$ and  $\phi_0(x_1,x_2)= 2 a_1\pi x_1+2 a_2\pi x_2$,
then we  define
\begin{equation}
I_1(\psi)=\int_Q |\nabla\psi|^2+\mathrm{H}^2 \left(1+\cos\left(\psi+\phi_0\right)\right)
\mathrm{~d}x.    
\end{equation}

Let $E=M=\dot{H}^1_p$ which is a reflexive Banach space,
and 
\( M \) is a nonempty weakly closed subset of \( E \).
Thus,
we have $0 \in M$, which means that $M$ is nonempty
and \( I_1: M \to \mathbb{R} \) is  not identically \( +\infty \).
It is obvious that $I_1(\psi) \geqslant 0$ for any  $\psi\in M$, so $I_1$ is lower bounded. 
For $ \psi\in \dot{H}^1_p$, 
since there exist two constants $c_1,c_2$  such that 
$c_1 \|\psi\|_{H^1_p}
\leqslant 
\|\nabla\psi\|_{L^2_p}
\leqslant
c_2 \|\psi\|_{H^1_p},
$
then $I_1(\psi)$ is coercive.

Now we prove that 
\( I_1: M \to \mathbb{R} \) is weakly lower semicontinuous.
Since
 $\int_{Q}|\nabla\psi|^2 \mathrm{~d}x$ is convex about $\psi$, then $\int_{Q}|\nabla\psi|^2 \mathrm{~d}x$ is weakly lower semicontinuous. 
Next, we consider $\int_{Q} \left( 1+\cos\left(\psi+\phi_0\right)\right) \mathrm{~d}x$.
Suppose that the sequence \( \{\psi_n\} \subset M\) such that \( \psi_n \rightharpoonup \psi_0\)
in $\dot{H}^1_p$,
then we have $\psi_n \rightarrow \psi_0$ in $L^2$.
By Egorov's theorem, for any $ \epsilon>0$, 
there exists $Q_\epsilon \subset Q$ such that $\left|Q-Q_\epsilon\right| <\epsilon$ and $\psi_n \rightrightarrows \psi_0$ in $Q_\epsilon$.
As 
$n\rightarrow \infty$,
we have
$$
\begin{aligned}
&\int_{Q_\epsilon} \left( 1+\cos\left(\psi_n+\phi_0\right)\right)-\left( 1+\cos\left(\psi_0+\phi_0\right)\right) \mathrm{~d}x \\ 
=&
\int_{Q_\epsilon} 2\cos\left(\frac{\psi_n+\psi_0 }{2} +\phi_0\right)\sin\left(\frac{ \psi_n-\psi_0 }{2} \right) \mathrm{~d}x \\
\geqslant &
\int_{Q_\epsilon} -2\left|\cos\left(\frac{\psi_n+\psi_0 }{2} +\phi_0\right)\right|
\left|\sin\left(\frac{ \psi_n-\psi_0 }{2} \right)\right| \mathrm{~d}x
\rightarrow 0.
\end{aligned}   
$$
Since $\int_{Q}( 1+\cos\psi) \mathrm{~d}x \geqslant 0$,
we have
\begin{equation}
\begin{aligned}
    \liminf\limits_{n\rightarrow \infty} \int_{Q} \left( 1+\cos\left(\psi_n+\phi_0\right)\right) \mathrm{~d}x
 &   \geqslant
    \liminf\limits_{n\rightarrow \infty} \int_{Q_\epsilon}
    \left( 1+\cos\left(\psi_n+\phi_0\right)\right) \mathrm{~d}x\\
  &  \geqslant
    \liminf\limits_{n\rightarrow \infty} \int_{Q_\epsilon}
   \left( 1+\cos\left(\psi_0+\phi_0\right)\right) \mathrm{~d}x.
    \end{aligned}
\end{equation}
Let $\epsilon\rightarrow 0$, then we can obtain that
$\int_{Q}( 1+\cos\psi) \mathrm{~d}x$ is weakly lower semicontinuous.
Therefore, \( I_1: M \to \mathbb{R} \) is weakly lower semicontinuous.

 By utilizing the direct methods, 
we can get $\psi\in M$ satisfying 
$-\Delta\psi=\mathrm{H}^2\sin\left( 
\psi+\phi_0  \right).
$
Let $\Psi=\psi+\phi_0$,
then we have $\Delta\psi=\Delta\Psi$.
Therefore, we can obtain  that
$\Psi$ is a solution  to (\ref{iso_modulo_periodic_problem}).    $\Box$

\begin{remark}
For $\phi_0:Q\rightarrow\mathbb{R}$ and  $\phi_0(x_1,x_2)= 2 a_1\pi x_1+2 a_2\pi x_2$,
we can derive the equivalent form of the problem (\ref{iso_modulo_periodic_problem}),
\begin{equation}\label{iso_psi_modulo_periodic_refine_problem}
 \left\{\begin{aligned}
-\Delta \psi(x) & = \;\mathrm{H}^2 \sin \left(\psi(x)+\phi_0(x)\right),   && x\in  Q, \\
\psi(x+\mathbf{e}_i) & = \;\psi(x),   && x\in \partial Q.
\end{aligned}\right.   
\end{equation}
The equation in (\ref{iso_psi_modulo_periodic_refine_problem}) corresponds to the Euler--Lagrange equation of the following energy functional
$$
   I(\psi)=\int_Q \frac{1}{2} |\nabla\psi|^2 +\mathrm{H}^2\left(1+\cos\left(\psi+\phi_0\right)\right) \mathrm{~d}x. 
$$
Since \( |a_1|+|a_2|>0  \),
then all of the solutions to (\ref{iso_modulo_periodic_problem}) are not constants.
\end{remark}

\paragraph{}
If $\mathrm{H}=0$, for any constant $C$, we have $\phi_0(x_1,x_2)= C+ 2 a_1\pi x_1+2 a_2\pi x_2$ is the solution to (\ref{iso_modulo_periodic_problem}).
After obtaining  the existence of the solutions to (\ref{iso_modulo_periodic_problem}), 
we can consider the uniqueness in case of $0<\mathrm{H}^2\leqslant \lambda_2$.
The uniqueness is defined in Remark \ref{remark_3_1}.

\begin{theorem}\label{Theorem_unique_first_1}
    Suppose  that $|a_1|+|a_2|>0$  and $0<\mathrm{H}^2\leqslant \lambda_2$, then there exists a unique nonconstant  solution to (\ref{iso_modulo_periodic_problem}).
\end{theorem}

\paragraph{Proof.}
For $0<\mathrm{H}^2 < \lambda_2$,
assume that  $\psi_1=\psi_1(x)$ and $\psi_2=\psi_2(x)$ are  different solutions to (\ref{iso_modulo_periodic_problem}),
which means that we cannot find both a  constant $\tau_0$ and an integer constant $k$ such 
that  $\psi_1(x)=\psi_2(x+\tau_0)+2k\pi$.
Since $\psi_1-\psi_2$ is periodic, then 
\begin{equation}\label{3.15}
    0
    <
    \int_Q |\nabla(\psi_1-\psi_2)|^2 \mathrm{~d}x
    =\mathrm{H}^2 \int_Q (\psi_1-\psi_2)(\sin\psi_1-\sin\psi_2) \mathrm{~d}x.
\end{equation}

Without loss of generality, assume that $a_1>0$.
Denoted by $m(\tau)=\int_Q \psi_2(x+\tau\mathbf{e}_1) \mathrm{~d}x$, where $\tau\in \mathbb{R}^2$,
then $m(0)=\int_Q \psi_2(x) \mathrm{~d}x$ and
\begin{equation}
    m(1)=\int_Q \psi_2(x+\mathbf{e}_1) \mathrm{~d}x
        =\int_Q \psi_2(x)+2a_1\pi \mathrm{~d}x
        =m(0)+2a_1\pi.
\end{equation}
Since $a_1>0$ is an integer constant and $\psi_2$ is a smooth function, 
then there exists a constant $\tau\in [0,1]$ and  an integer constant  $ k_0$ such that             
 $m(\tau)-m(0)=\int_Q \psi_1(x)-\psi_2(x) \mathrm{~d}x+2k_0\pi$.

Let $\psi_3(x)=\psi_2(x+\tau\mathbf{e}_1)-2 k_0\pi$, then
$\int_Q \psi_1(x)-\psi_3(x) \mathrm{~d}x=0$ and
\begin{equation}
 \int_Q (\psi_1-\psi_3)^2 \mathrm{~d}x
 =\int_Q \left(\psi_1-\psi_3 - \int_Q \psi_1(x)-\psi_3(x) \mathrm{~d}x \right)^2 \mathrm{~d}x.
\end{equation}
Since we have
\begin{equation}
    \begin{aligned}
    \int_Q |\nabla(\psi_1-\psi_3)|^2 \mathrm{~d}x
    =\mathrm{H}^2 \int_Q (\psi_1-\psi_3)(\sin\psi_1-\sin\psi_3) \mathrm{~d}x
    < \lambda_2  \int_Q (\psi_1-\psi_3)^2 \mathrm{~d}x.
    \end{aligned}
\end{equation}
Then $\psi_1   \equiv   \psi_3$,
since $\lambda_2$ is the principal eigenvalue satisfying (\ref{lambda_2_def}).

On the other hand, if $\mathrm{H}^2=\lambda_2$, 
consider the equality conditions,
$$
    2\sin\left(\frac{\psi_1-\psi_3 }{ 2}\right)
    =
    \psi_1-\psi_3,\quad
    \cos\left(\frac{\psi_1+\psi_3 }{ 2}\right)=1.
$$
When $|a_1|+|a_2|>0$, 
problem (\ref{iso_modulo_periodic_problem}) doesn't have a trivial solution.
Thus,
we can obtain
$\psi_1   =  \psi_3 $.
Therefore,
we get the uniqueness of the nonconstant solution to (\ref{iso_modulo_periodic_problem}).   $\Box$

Theorem \ref{Theorem_1.1} is proved.

\section{Uniform a priori estimates}

\setcounter{equation}{0}
\numberwithin{equation}{section}

To prove the global existence of solutions to the general Ericksen--Leslie system  (\ref{iso_Ericksen_Leslie_new_formulation}) with initial conditions (\ref{initial_condition_1}) and boundary conditions (\ref{boundary_condition_1}), we require uniform a priori estimates, which are higher-order estimates. Then we can 
take the limit as \( m \to \infty \) to obtain a global strong solution to the general Ericksen--Leslie system system (\ref{iso_Ericksen_Leslie_new_formulation})  
and extend the local solution to a global solution in \([0, \infty)\).
For simplicity, we use $\left(\mathbf{v}, \mathbf{d}\right)$  instead of the smooth solution $\left(\mathbf{v}_m, \mathbf{d}_m\right)$ in Section \ref{section5} to the approximate problem (\ref{iso_approximate_problem_d_first})--(\ref{iso_approximate_boundary_d}).

Define the higher-order energy
\begin{equation}\label{iso_higher_order_energy}
E_H(t) := 
\int_Q |\nabla \mathbf{v}|^2 
+\frac{1-\gamma}{\operatorname{Re}} 
|\Delta\mathbf{d}+f(\mathbf{d},\nabla\mathbf{d}) |^2 
\mathrm{~d} x.
\end{equation}

\subsection{\texorpdfstring{Estimate for $\sigma_1$}{estimate for} 
}

\begin{theorem}\label{iso_higher-order_sigma_1}
Suppose that
$\mathbf{v} $  and  $\mathbf{d}\in\mathbb{S}^1$
are smooth functions, then
\begin{equation}
\int_Q \nabla \mathbf{v}:\nabla(\nabla\cdot \sigma_1) \mathrm{~d} x 
\leqslant 
\epsilon_0\left(\|\Delta\mathbf{v}\|^2 
+\|\nabla\Delta\mathbf{d}\|^2
\right) +C_0 E_H(t),
\end{equation}
where  $\epsilon_0$ and  $C_0=C_0(\epsilon_0)$
are positive  constants,
depending on  $Q$, $|\mathbf{H}|$, $\left\|\mathbf{d}_0\right\|_{H^2}$,
$\left\|\mathbf{v}_0\right\|_{H^1}$ and coefficients of the system (\ref{iso_Ericksen_Leslie_new_formulation}).
\end{theorem}

\paragraph{Proof.} 
Since $\nabla \mathbf{v}=\mathbf{D}+\boldsymbol{\Omega} $, 
we have
\begin{equation}\label{iso_higher-order_sigma_1_directly_1}
\begin{aligned}
\int_Q \nabla \mathbf{v}:\nabla(\nabla\cdot \sigma_1) \mathrm{~d} x
=&-\beta_1 \int_Q \Delta \mathbf{v}\cdot[\nabla\cdot((\mathbf{D}: \mathbf{d}\otimes\mathbf{d}) \mathbf{d}\otimes\mathbf{d} )] \mathrm{~d} x\\
&
-\beta_2 \int_Q \Delta \mathbf{v}\cdot(\nabla\cdot\mathbf{D})\mathrm{~d} x\\
&-\frac{ \beta_3 }{2} \int_Q \Delta \mathbf{v}\cdot[\nabla\cdot(\mathbf{d}\otimes( \mathbf{D}  \mathbf{d})+(\mathbf{D}  \mathbf{d}\otimes\mathbf{d}) )] \mathrm{~d} x\\
:=&\;M_{11}+M_{12}+M_{13}.   
\end{aligned}
\end{equation}

For $M_{11}$,
$$
\begin{aligned}
M_{11}
=&-\beta_1 \int_Q (\nabla \cdot\mathbf{D})\cdot[\nabla\cdot((\mathbf{D}: \mathbf{d}\otimes\mathbf{d}) \mathbf{d}\otimes\mathbf{d} )] \mathrm{~d} x\\
\leqslant& -\beta_1
\int_Q  \sum_{k=1}^2 |\mathbf{D}_{x_k}:\mathbf{d}\otimes\mathbf{d}|^2 \mathrm{~d} x+\epsilon_1\left(\|\Delta\mathbf{v}\|^2 
+\|\nabla\Delta\mathbf{d}\|^2
\right)
+C_2 E_H(t).
\end{aligned}
$$

For $M_{12}$,
$$	
\begin{aligned}
M_{12}
\leqslant
-\beta_2 \int_Q |\nabla\mathbf{D}|^2
\mathrm{~d} x.
\end{aligned}
$$

For $M_{13}$,
$$	
\begin{aligned}
	M_{13}
	=&\;\frac{ \beta_3 }{2} \int_Q \Delta(\mathbf{D}+\boldsymbol{\Omega}):[\mathbf{d}\otimes( \mathbf{D}  \mathbf{d})+(\mathbf{D}  \mathbf{d})\otimes\mathbf{d})] \mathrm{~d} x\\
	=&\;\beta_3  \int_Q (\Delta\mathbf{D}\cdot \mathbf{d})  \cdot ( \mathbf{D} \mathbf{d} ) \mathrm{~d} x\\
	\leqslant &-\beta_3  
 \int_Q   \sum_{k=1}^2  |\mathbf{D}_{x_k}\cdot \mathbf{d}|^2   \mathrm{~d} x +\epsilon_1\left(\|\Delta\mathbf{v}\|^2 
+\|\nabla\Delta\mathbf{d}\|^2
\right)
+C_2 E_H(t).
\end{aligned}
$$

Let $\mathbf{D}_{x_k} :=(\partial_k D_{ij})_{2\times 2}$,
which is a symmetric trace free matrix.
From Lemma \ref{beta_lemma}, we have
\begin{equation}\label{iso_higher_sigma_1_2}
\begin{aligned}
	&\beta_1 \int_Q  \sum_{k=1}^2 |\mathbf{D}_{x_k}:\mathbf{d}\otimes\mathbf{d}|^2 \mathrm{~d} x
      +\beta_2 \int_Q |\nabla\mathbf{D}|^2
     \mathrm{~d} x
	+ \beta_3  \int_Q   \sum_{k=1}^2  |\mathbf{D}_{x_k}\cdot \mathbf{d}|^2   \mathrm{~d} x  \\
	=&\sum_{k=1}^2 \int_Q \beta_1|\mathbf{D}_{x_k}:\mathbf{d}\otimes\mathbf{d}|^2 
	+ \beta_2  |\mathbf{D}_{x_k}|^2
	+ \beta_3  |\mathbf{D}_{x_k}\cdot \mathbf{d}|^2   \mathrm{~d} x \geqslant 0.
\end{aligned}
\end{equation}
Let $C_0=2C_2$ and $\epsilon_0=2\epsilon_1$.
By (\ref{iso_higher-order_sigma_1_directly_1}) and (\ref{iso_higher_sigma_1_2}),
Theorem \ref{iso_higher-order_sigma_1} is proved.    $\Box$

\begin{remark}
If $\mathbf{d}:Q\rightarrow\mathbb{S}^2$,
where $Q$ is a  unit cube in  $\mathbb{R}^2$ or $\mathbb{R}^3$,
after changing the dissipation relations with respect to $\beta_i$ in Remark \ref{Remark_iso_beta_first},
the conclusion in Theorem \ref{iso_higher-order_sigma_1} remains valid.
\end{remark}

\subsection{Coercive estimates for  the harmonic map }

\begin{theorem}\label{iso_thm__higher_order_harmonic_map_estimate_2}
Suppose that
 $\mathbf{d}\in\mathbb{S}^1$
is a smooth  function, then
we have the following coercive estimates regarding the harmonic map,
\begin{equation}\label{iso__higher_order_harmonic_map_1}
    \|\nabla^2\mathbf{d}\|
\leqslant
    C\left(\| \Delta\mathbf{d}+|\nabla\mathbf{d}|^2\mathbf{d} \|
+1
\right),
\end{equation}
\begin{equation}\label{iso__higher_order_harmonic_map_2}
    \|\nabla^3\mathbf{d}\|
\leqslant
    C\left(
    \| \Delta \mathbf{d} + |\nabla\mathbf{d}|^2 \mathbf{d} \|^2
    +
    \| \nabla(\Delta\mathbf{d}+|\nabla\mathbf{d}|^2\mathbf{d}) \|
    +1 \right),
\end{equation}
where $C$ is a constant, depending on 
$\left\|\nabla\mathbf{d}\right\|$.
\end{theorem}

\paragraph{Proof.}
Let $\mathbf{d}=(d_1,d_2)^T=(\sin\theta,\cos\theta)^T$
and
$\mathbf{d}^{\perp}=(\cos\theta,-\sin\theta)^T$.

(i) First, we have
$$
\begin{aligned}
 \|\nabla^2\mathbf{d}\|^2
=&\;\int_Q (\partial_j \partial_k d_i)(\partial_j \partial_k d_i) \mathrm{~d} x\\
=&\;\int_Q \partial_j\partial_j d_i  \partial_k\partial_k d_i  \mathbf{~d}x\\
\lesssim &\;\int_Q  |\Delta\mathbf{d}|^2 
\mathbf{~d}x.
\end{aligned}
$$
Similarly, it can be deduced that 
$$ 
\|\nabla^2\theta\| \lesssim\|\Delta\theta\|=\| \Delta\mathbf{d}+|\nabla\mathbf{d}|^2\mathbf{d} \|.
$$
Thanks to  
$\Delta\mathbf{d}+|\nabla\mathbf{d}|^2\mathbf{d}
=(\Delta\theta )\mathbf{d}^{\perp}$
and
$
\left(  (\Delta\theta )\mathbf{d}^{\perp},  |\nabla\mathbf{d}|^2\mathbf{d}   \right)
=0$,
then we have
$$
\|\Delta\mathbf{d}\|^2
= \| \Delta\mathbf{d}+|\nabla\mathbf{d}|^2\mathbf{d} \|^2
+\|  |\nabla\mathbf{d}|^2\mathbf{d} \|^2
=\|\Delta\theta\|^2+\|\nabla\mathbf{d}\|^4_{L^4}.
$$
By $\nabla\mathbf{d}= \mathbf{d}^{\perp}\otimes\nabla\theta$ and Ladyzhenskaya's inequality (Lemma \ref{lem_inequalities}),
$$
\|\nabla\mathbf{d}\|^4_{L^4}
=\|\nabla\theta\|^4_{L^4}
\leqslant
C\|\nabla\theta\|^2 
\left(
\|\nabla^2\theta\|^2+\|\nabla\theta\|^2
\right)
\leqslant
C
\left(
\|\Delta\theta\|^2+1
\right).
$$
We can get (\ref{iso__higher_order_harmonic_map_1}) from the above inequalities.

\paragraph{}
(ii) Second, we can obtain the following estimates
$$
\begin{aligned}
\|\nabla^3\mathbf{d}\|^2
=&\;\int_Q (\partial_j \partial_k\partial_l d_i)(\partial_j \partial_k\partial_l d_i) \mathrm{~d} x\\
=&\;\int_Q (\partial_l\partial_j\partial_j d_i )
(\partial_l\partial_k\partial_k d_i ) \mathbf{~d}x\\
\lesssim 
&\;\int_Q |\nabla\Delta \mathbf{d}|^2
\mathbf{~d} x,
\end{aligned}
$$
and
$$
\begin{aligned}
\|\nabla(\Delta\mathbf{d}+|\nabla\mathbf{d}|^2\mathbf{d})\|^2
=&\;\|\nabla(\mathbf{d}^{\perp}\Delta\theta)\|^2\\
=&\;\|-\mathbf{d}\otimes\nabla\theta \Delta\theta
+\mathbf{d}^{\perp}\otimes\nabla\Delta\theta\|^2\\
=&\;
\|\mathbf{d}\otimes\nabla\theta \Delta\theta\|^2
+\|\mathbf{d}^{\perp}\otimes\nabla\Delta\theta\|^2\\
=&\;\|\nabla\theta\Delta\theta\|^2
+\|\nabla\Delta\theta \|^2.
\end{aligned}
$$
Calculate $\|\nabla(|\nabla\mathbf{d}|^2\mathbf{d})\|^2$  directly, 
then we have
$$
\begin{aligned}
\|\nabla(|\nabla\mathbf{d}|^2\mathbf{d})\|^2
=&
\int_Q \left[\left|\partial_i\mathbf{d}_j\right|^2 \mathbf{d}^{\perp}\otimes\nabla\theta
+\mathbf{d}\otimes \nabla(\partial_i\mathbf{d}_j)^2
\right]^2
\mathrm{~d}x\\
=&\int_Q \left|  \left|  \nabla\theta    \right|^2 \nabla\theta    \right|^2
\mathrm{~d}x
+\int_Q \left| \nabla \left|  \nabla\theta    \right|^2    \right|^2
\mathrm{~d}x.
\end{aligned}
$$
Thanks to the following identity,
$$
\begin{aligned}
\int_Q \left| \nabla \left|  \nabla\theta    \right|^2    \right|^2
\mathrm{~d}x
=&\;4\int_Q \left(    \partial_k\partial_i\theta \right)
\left[ \left( \partial_i\theta     \right)^2   \partial_k\partial_i\theta \right]
\mathrm{~d}x\\
=&-4\int_Q  \partial_i\theta\partial_k
\left[  \left( \partial_i\theta     \right)^2   \partial_i\partial_k\theta \right]
\mathrm{~d}x\\
=&-2\int_Q \left| \nabla \left|  \nabla\theta    \right|^2    \right|^2
\mathrm{~d}x
-4\int_Q \left( \partial_i\theta     \right)^3 
 \partial_i\Delta\theta  
 \mathrm{~d}x,
\end{aligned}
$$
and Gagliardo--Nirenberg inequality (Lemma \ref{lem_inequalities}),
$$
\|\nabla\theta\|_{L^6}
\leqslant
C\|\nabla\theta\|^{\frac{1}{3}}
\left( \|\nabla^2\theta\|+   \|\nabla\theta\|  \right)^{\frac{2}{3}}
\leqslant C\|\nabla\theta\|^{\frac{1}{3}}
\left( \|\Delta\theta\|+   \|\nabla\theta\|  \right)^{\frac{2}{3}},
$$
then we have
$$
\begin{aligned}
\|\nabla(|\nabla\mathbf{d}|^2\mathbf{d})\|^2
=&\;\|\nabla\theta\|^6_{L^6}
+\frac{4}{3}\int_Q \left( \partial_i\theta     \right)^3 
 \partial_i\left(-\Delta\theta  \right)
 \mathrm{~d}x\\
\leqslant &\;
C\left( \|\nabla\theta\|^6_{L^6} 
+\|\nabla\Delta\theta \|^2
\right)\\
\leqslant &\;
C\left( \|\Delta\theta\|^4
+\|\nabla\Delta\theta \|^2
+1
\right).
\end{aligned}
$$
We can get (\ref{iso__higher_order_harmonic_map_2}) from the above inequalities. $\Box$

\subsection{Uniform a priori estimates}

\begin{theorem}\label{iso_Theorem_higher_order_final}
Suppose that
$\mathbf{v} $  and  $\mathbf{d}\in\mathbb{S}^1$
are smooth functions, 
then we have the following higher-order energy estimate
\begin{equation}\label{iso_higher-order_energy_estimate}
\frac{\mathrm{d}}{\mathrm{dt}}E_H
\leqslant
C(E_H^2(t)+E_H(t)).
\end{equation}	
The constant $C$  depends on $Q, |\mathbf{H}|,\left\|\mathbf{d}_0\right\|_{H^2}$, 
$\left\|\mathbf{v}_0\right\|_{H^1}$
and
coefficients of the system (\ref{iso_Ericksen_Leslie_new_formulation}).
\end{theorem}

\paragraph{Proof.}
(i) Show the terms we need to estimate.

Calculate the left hand of \eqref{iso_higher-order_energy_estimate}  directly
\begin{equation}\label{4.18}
\begin{aligned}	\frac{1}{2}\frac{\mathrm{d}}{\mathrm{dt}}E_H(t)
&=\frac{1}{2}\frac{\mathrm{d}}{\mathrm{dt}}\int_Q|\nabla \mathbf{v}|^2
+ \frac{1-\gamma}{\operatorname{Re}} 
|\Delta\mathbf{d} +f(\mathbf{d},\nabla\mathbf{d}) |^2 dx \\
         &=-(\Delta \mathbf{v},\mathbf{v}_t)+\frac{1-\gamma}{\operatorname{Re}} (\Delta\mathbf{d}+f(\mathbf{d},\nabla\mathbf{d}),\Delta\mathbf{d}_t+\partial_t f( \mathbf{d},\nabla\mathbf{d})).
    \end{aligned}
\end{equation}
For $ -(\Delta \mathbf{v},\mathbf{v}_t) $,
we have
\begin{equation}\label{4.19}
\begin{aligned}
    &-(\Delta \mathbf{v},\mathbf{v}_t)+\frac{\gamma}{\operatorname{Re}} ( \Delta \mathbf{v},\Delta \mathbf{v})
    +\frac{1-\gamma}{\operatorname{Re}} (\Delta \mathbf{v},\nabla \cdot \sigma_1)\\
=&\;(\Delta \mathbf{v}, \mathbf{v} \cdot \nabla \mathbf{v})
    +\frac{1-\gamma}{\operatorname{Re}}  (\Delta \mathbf{v}, \Delta \mathbf{d} \cdot \nabla \mathbf{d})\\
&+\frac{1-\gamma}{\operatorname{Re}} \left(\Delta \mathbf{v},\nabla 
   \cdot\left[\frac{1}{2}\left(1+\mu_2\right) (\Delta\mathbf{d}+f(\mathbf{d},\nabla\mathbf{d}) )\otimes\mathbf{d}\right]\right)\\
& +\frac{1-\gamma}{\operatorname{Re}} \left(\Delta \mathbf{v},\nabla \cdot\left[\frac{1}{2}\left(-1+\mu_2\right)\mathbf{d}\otimes (\Delta\mathbf{d}+f(\mathbf{d},\nabla\mathbf{d}) ) \right]\right)\\
:=&\; I_1+I_2+I_3+I_4.
\end{aligned}
\end{equation}
Let $\epsilon$ is small enough, then we have
$$
I_1\leqslant \|\Delta\mathbf{v}\| 
\| \mathbf{v} \|_{L^{\infty}}
\| \nabla \mathbf{v} \|
\leqslant C\|\Delta\mathbf{v}\|^{\frac{3}{2}} 
\| \nabla \mathbf{v} \|
\leqslant \epsilon\|\Delta\mathbf{v}\|^2+C\|\nabla\mathbf{v}\|^4.
$$ 

 For $ \frac{1-\gamma}{\operatorname{Re}} (\Delta\mathbf{d}+f(\mathbf{d},\nabla\mathbf{d}),\Delta\mathbf{d}_t+\partial_t f(\mathbf{d},\nabla\mathbf{d}))$,
we have
\begin{equation}\label{4.20}
    \begin{aligned}
        &\;\frac{1-\gamma}{\operatorname{Re}} (\Delta\mathbf{d}+f(\mathbf{d},\nabla\mathbf{d}),\Delta\mathbf{d}_t+\partial_t f(\mathbf{d},\nabla\mathbf{d}))
        +\frac{\mu_1(1-\gamma)}{\operatorname{Re}}\|\nabla(\Delta \mathbf{d}+f(\mathbf{d},\nabla\mathbf{d}) )\|^2\\
        =&\;\frac{1-\gamma}{\operatorname{Re}} 
        \left(\Delta \mathbf{d}+f(\mathbf{d},\nabla\mathbf{d}),\Delta\left(-\mathbf{v}\cdot\nabla\mathbf{d}\right)\right)\\
        &+\frac{1-\gamma}{\operatorname{Re}} 
        \left(\Delta \mathbf{d}+f(\mathbf{d},\nabla\mathbf{d}),\Delta\left((\mathbf{I}-\mathbf{d}\otimes\mathbf{d})\left(\boldsymbol{\Omega}  \mathbf{d}
+\mu_2 \mathbf{D} \mathbf{d}\right)\right)
        \right)
        \\
        &+\frac{1-\gamma}{\operatorname{Re}} 
        \left(\Delta \mathbf{d}+f(\mathbf{d},\nabla\mathbf{d}),
        \partial_t f(\mathbf{d},\nabla\mathbf{d}) 
        \right)\\
        :=&\; J_1+J_2+J_3.
    \end{aligned}
\end{equation}

\paragraph{}
(ii) 
Show coupling relationships of the higher-order terms.

We need to divide $J_1$ in order to kill $I_2$,
$$
\begin{aligned}
    J_1=&-\frac{1-\gamma}{\operatorname{Re}} 
        \left(\Delta \mathbf{d}+f(\mathbf{d},\nabla\mathbf{d}),\Delta\left(
        \mathbf{v}\cdot\nabla\mathbf{d}\right)\right)\\
       =&-\frac{1-\gamma}{\operatorname{Re}} 
        \left(\Delta \mathbf{d}+f(\mathbf{d},\nabla\mathbf{d}),\left(\Delta
        \mathbf{v}\right)\cdot\nabla\mathbf{d}\right)
        -\frac{1-\gamma}{\operatorname{Re}} 
        \left(\Delta \mathbf{d}+f(\mathbf{d},\nabla\mathbf{d}),
        \mathbf{v}\cdot\nabla\Delta\mathbf{d}\right)\\
        &-\frac{2(1-\gamma)}{\operatorname{Re}} 
        \left(\Delta \mathbf{d}+f(\mathbf{d},\nabla\mathbf{d}),
        \nabla\mathbf{v}\cdot\nabla^2\mathbf{d}\right)\\
       :=&\;J_{11}+J_{12}+J_{13}, \\
    I_2+J_{11}=&\;\frac{1-\gamma}{\operatorname{Re}} 
        \left(f(\mathbf{d},\nabla\mathbf{d}),\left(\Delta
        \mathbf{v}\right)\cdot\nabla\mathbf{d}\right)
        =
        \frac{1-\gamma}{\operatorname{Re}} 
\left((\mathbf{H}\cdot\mathbf{d})\mathbf{H},\left({\Delta}\mathbf{v}\right)\cdot\nabla\mathbf{d}\right)
        := K_1.
        \end{aligned}
$$

We can couple $J_{12}$ with $J_3$,
after calculating them.
$$
\begin{aligned}
J_{12}
=&-\frac{1-\gamma}{\operatorname{Re}} 
        \left(\Delta \mathbf{d}+f(\mathbf{d},\nabla\mathbf{d}),
        \mathbf{v}\cdot\nabla (\Delta\mathbf{d}+f(\mathbf{d},\nabla\mathbf{d}) )\right)\\
        &+\frac{1-\gamma}{\operatorname{Re}} 
        \left(\Delta \mathbf{d}+f(\mathbf{d},\nabla\mathbf{d}),
        \mathbf{v}\cdot\nabla f(\mathbf{d},\nabla\mathbf{d}) \right)\\
=&\frac{1-\gamma}{\operatorname{Re}} 
        (\Delta \mathbf{d}+f(\mathbf{d},\nabla\mathbf{d}),
        \mathbf{v}_i|\nabla\mathbf{d}|^2\partial_i\mathbf{d}
        +\mathrm{H}^2\mathbf{v}_i\partial_i\theta\cos2\theta \mathbf{d}^{\perp}),\\
J_{3}=&\;\frac{1-\gamma}{\operatorname{Re}} 
        (\Delta \mathbf{d}+f(\mathbf{d},\nabla\mathbf{d}),
|\nabla\mathbf{d}|^2\mathbf{d}_t+\mathrm{H}^2\theta_t\cos2\theta\mathbf{d}^{\perp}+2(\nabla\mathbf{d}:\nabla\mathbf{d}_t)\mathbf{d}+\frac{H^2}{2} \sin 2\theta\mathbf{d}_t^{\perp}).
\end{aligned}
$$
Then we get
$$
\begin{aligned}
J_3+J_{12}
=&\;
    \frac{1-\gamma}{\operatorname{Re}} 
        (\Delta \mathbf{d}+f(\mathbf{d},\nabla\mathbf{d}),
|\nabla\mathbf{d}|^2(\mathbf{d}_t+\mathbf{v}\cdot\nabla\mathbf{d}))\\
&+\frac{1-\gamma}{\operatorname{Re}} 
(\Delta \mathbf{d}+f(\mathbf{d},\nabla\mathbf{d}),
\mathrm{H}^2\cos2\theta(\theta_t+\mathbf{v}\cdot\nabla\theta)\mathbf{d}^{\perp})\\
:= &\; K_{21}+K_{22}.
\end{aligned}
$$

 To deal with $J_1$, we need to estimate $J_{13} $. 
By Amgon's inequality (Lemma \ref{lem_inequalities}),
$$
\|\nabla\mathbf{d}\|_{L^{\infty}}\leqslant C( 1+\|\nabla\Delta\mathbf{d}\|^{\frac{1}{2}} )
$$
where C depends on $\|\nabla\mathbf{d}\|$.
From Theorem \ref{iso_thm__higher_order_harmonic_map_estimate_2},
$$
\begin{aligned}
\left|  J_{13}  \right|
\leqslant &\;
C\left[   
\left|\left( \Delta\mathbf{d},   \nabla\mathbf{v}\cdot\nabla^2\mathbf{d}          \right)\right|
+\left|\left( |\nabla\mathbf{d}|^2\mathbf{d},   \nabla\mathbf{v}\cdot\nabla^2\mathbf{d}          \right)\right|
+\|\nabla\mathbf{v}\| \|\nabla^2\mathbf{d}\|
\right]\\
\leqslant &\;
C\left[   
\left\| \Delta\mathbf{d}   \right\|_{L^4}  
\left\| \nabla\mathbf{v}  \right\|_{L^4} 
\left\| \nabla^2\mathbf{d}  \right\| 
+
\left\| \nabla\mathbf{d}   \right\|^2_{L^{\infty}}  
\left\| \nabla\mathbf{v}  \right\|
\left\| \nabla^2\mathbf{d}  \right\|
+
E_H^2+E_H
\right]\\
\leqslant &\;
C\left(  
\left\| \Delta\mathbf{d}   \right\|^{\frac{1}{2}}
\left\| \nabla\Delta\mathbf{d}   \right\|^{\frac{1}{2}}
\left\| \nabla\mathbf{v}  \right\|^{\frac{1}{2}}
\left\| \Delta\mathbf{v}  \right\|^{\frac{1}{2}}
\left\| \Delta\mathbf{d}  \right\| 
+
 \left\| \nabla\Delta\mathbf{d}   \right\|
\left\| \nabla\mathbf{v}  \right\|
\left\| \Delta\mathbf{d}  \right\|
+
E_H^2+E_H
\right)\\
\leqslant &\;
C\left[   
E_H^2+E_H
\right]
+\epsilon 
\left(  
\left\| \nabla \left(\Delta\mathbf{d} +f(\mathbf{d},\nabla\mathbf{d}) \right)  \right\|^2 
+
\left\| \Delta\mathbf{v}  \right\|^2
\right).
\end{aligned}
$$

$J_2$ can be divided
$$
    \begin{aligned}
J_2=&\;\frac{1-\gamma}{\operatorname{Re}} \left(\Delta \mathbf{d}+f(\mathbf{d},\nabla\mathbf{d}),\Delta \left[\frac{1}{2} 
     (1+\mu_2)(\nabla\mathbf{v}\cdot\mathbf{d})  \right]\right)\\
     &+\frac{1-\gamma}{\operatorname{Re}} \left(\Delta \mathbf{d}+f(\mathbf{d},\nabla\mathbf{d}),\Delta \left[\frac{1}{2}(-1+\mu_2)((\nabla\mathbf{v})^{T}\cdot\mathbf{d})  \right]\right)\\
    & -\frac{1-\gamma}{\operatorname{Re}} \left(\Delta \mathbf{d}+f(\mathbf{d},\nabla\mathbf{d}),\Delta \left[(\mathbf{D}:\mathbf{d}\otimes\mathbf{d})\mathbf{d}  \right]\right)\\ 
    :=&\; J_{21}+J_{22}+J_{23},
    \end{aligned}
$$
where
 $J_{21}$ can be coupled with $I_3$.
$$
\begin{aligned}
    I_3+J_{21}=&-\frac{(1-\gamma)(1+\mu_2)}{2\operatorname{Re}}
             \left( \mathbf{d}\cdot\nabla \Delta \mathbf{v},\Delta\mathbf{d}+f(\mathbf{d},\nabla\mathbf{d}) \right)\\
             &+\frac{(1-\gamma)(1+\mu_2)}{2\operatorname{Re}}
             \left( \mathbf{d}\cdot\nabla \Delta \mathbf{v}+\nabla\mathbf{v}\cdot\Delta \mathbf{d},\Delta\mathbf{d}+f(\mathbf{d},\nabla\mathbf{d}) \right)
              \\
             &+\frac{(1-\gamma)(1+\mu_2)}{\operatorname{Re}}
             \left( \partial_k\mathbf{d}_j\partial_k\partial_j\mathbf{v}_i,\Delta\mathbf{d}+f(\mathbf{d},\nabla\mathbf{d}) \right)\\
            = &\;\frac{(1-\gamma)(1+\mu_2)}{2\operatorname{Re}}
             \left( \nabla\mathbf{v}\cdot\Delta \mathbf{d},\Delta\mathbf{d}+f(\mathbf{d},\nabla\mathbf{d}) \right)\\
             &+\frac{(1-\gamma)(1+\mu_2)}{\operatorname{Re}}
             \left( \partial_k\mathbf{d}_j\partial_k\partial_j\mathbf{v}_i,\Delta\mathbf{d}+f(\mathbf{d},\nabla\mathbf{d}) \right)\\
             :=&\; K_{31}+K_{32}.
\end{aligned}
$$
$J_{22}$ can be coupled with $I_4$.
$$
\begin{aligned}
    I_4+J_{22}=&\;\frac{(1-\gamma)(-1+\mu_2)}{2\operatorname{Re}}
             \left[
             -\left( \mathbf{d}\cdot\nabla (\Delta \mathbf{v})^T,\Delta\mathbf{d}+f(\mathbf{d},\nabla\mathbf{d}) \right)
             \right] \\ 
             &+\frac{(1-\gamma)(-1+\mu_2)}{2\operatorname{Re}}
             \left[
             \left( \mathbf{d}\cdot\nabla (\Delta \mathbf{v})^T+(\nabla\mathbf{v})^T\cdot\Delta \mathbf{d},\Delta\mathbf{d}+f(\mathbf{d},\nabla\mathbf{d}) \right)
             \right] 
             \\
             &+\frac{(1-\gamma)(-1+\mu_2)}{\operatorname{Re}}
             \left( \partial_k\mathbf{d}_j\partial_k\partial_i\mathbf{v}_j,\Delta\mathbf{d}+f(\mathbf{d},\nabla\mathbf{d}) \right)\\
            = &\;
            \frac{(1-\gamma)(-1+\mu_2)}{2\operatorname{Re}}
            \left[
             \left( (\nabla\mathbf{v})^T\cdot\Delta \mathbf{d}
   ,\Delta\mathbf{d}+f(\mathbf{d},\nabla\mathbf{d}) \right)\right]\\
             &+
             \frac{(1-\gamma)(-1+\mu_2)}{\operatorname{Re}}
            \left[\left( \partial_k\mathbf{d}_j\partial_k\partial_i\mathbf{v}_j,\Delta\mathbf{d}+f(\mathbf{d},\nabla\mathbf{d}) \right)\right]
             \\
             :=&\; K_{41}+K_{42}.
\end{aligned}
$$

 Since $ \left(\Delta \mathbf{d}+f(\mathbf{d},\nabla\mathbf{d}), \mathbf{d} \right)=0 $, then we have
$$
\begin{aligned}
    J_{23}=&-\frac{1-\gamma}{\operatorname{Re}} \left(\Delta \mathbf{d}+f(\mathbf{d},\nabla\mathbf{d}), 
    \left[\Delta(\mathbf{D}:\mathbf{d}\otimes\mathbf{d})\right]\mathbf{d} \right)\\
    &
    -\frac{1-\gamma}{\operatorname{Re}} \left(\Delta \mathbf{d}+f(\mathbf{d},\nabla\mathbf{d}),
    +(\mathbf{D}:\mathbf{d}\otimes\mathbf{d})(\Delta\mathbf{d})+2\nabla(\mathbf{D}:\mathbf{d}\otimes\mathbf{d})\cdot\nabla\mathbf{d}
    \right)\\
    =&-\frac{1-\gamma}{\operatorname{Re}} \left(\Delta \mathbf{d}+f(\mathbf{d},\nabla\mathbf{d}), 
    (\mathbf{D}:\mathbf{d}\otimes\mathbf{d}) \Delta\mathbf{d}\right)\\
    &-\frac{1-\gamma}{\operatorname{Re}} \left(\Delta \mathbf{d}+f(\mathbf{d},\nabla\mathbf{d}), 
    2\nabla(\mathbf{D}:\mathbf{d}\otimes\mathbf{d})\cdot\nabla\mathbf{d}
    \right)\\
    :=&\; K_{51}+K_{52}.
\end{aligned}
$$

(iii)
Finally, we have the following estimates to finish the proof.

 For $K_1$,
since
$(
|\nabla\mathbf{d}|^2\mathbf{d},
\Delta\mathbf{v} \cdot \nabla\mathbf{d}
)
=
(
|\nabla\theta|^2\mathbf{d},
(\Delta\mathbf{v} \cdot \nabla\theta) \mathbf{d}^{\perp}
)
=0
$,
we have
$$
\begin{aligned}
K_1
=
\frac{1-\gamma}{\operatorname{Re}} 
\left((\mathbf{H}\cdot\mathbf{d})\mathbf{H},\left({\Delta}\mathbf{v}\right)\cdot\nabla\mathbf{d}\right)
\leqslant
C E_H.
\end{aligned}
$$

For $K_{21}$,
$$
\begin{aligned}
|K_{21}|
\leqslant&\;
C\int_Q
|\nabla\mathbf{d}|^2
\left(
|\Delta\mathbf{d}+f(\mathbf{d},\nabla\mathbf{d}) |^2
+|\nabla\mathbf{v}|^2
\right)
\mathrm{~d}x\\
\leqslant&\;
C\|\nabla\mathbf{d}\|^2_{L^{\infty}}
\left(
\|\Delta\mathbf{d}+f(\mathbf{d},\nabla\mathbf{d}) \|^2
+\|\nabla\mathbf{v}\|^2
\right)\\
\leqslant&\;
C
\left( \| \nabla\Delta\mathbf{d} \| + 1 \right)
E_H
\\
\leqslant&\;
C
\left( \| \nabla(\Delta\mathbf{d}+f(\mathbf{d},\nabla\mathbf{d}) ) \| + 1+E_H \right)
E_H\\
\leqslant&\;
C
\left( E_H^2  +E_H \right)
+
\epsilon
\| \nabla(\Delta\mathbf{d}+f(\mathbf{d},\nabla\mathbf{d}) ) \|^2 .
\end{aligned}
$$

For $ K_{22} $,
$$
\begin{aligned}
|K_{22}|
\leqslant &\;
C|(
\Delta\theta+\frac{H^2}{2}\sin2\theta,
\Delta\theta+\frac{H^2}{2}\sin2\theta 
+|\nabla\mathbf{v}|
)|\\
\leqslant &\;
C\left(   \|  \Delta\theta+\frac{H^2}{2}\sin2\theta         \|^2
+\|   \nabla\mathbf{v}        \|^2    \right)
=C \; E_H.
\end{aligned}
$$

Let $ K_3=\max\left\{|K_{31}|,|K_{41}|,|K_{51}|\right\}$,
since we can deal with them by the same estimate.
$$
\begin{aligned}
K_3
\leqslant &\;
C
\|\Delta\mathbf{d}\|
\|\nabla\mathbf{v}\|_{L^4 }
\|  \Delta\mathbf{d}+f(\mathbf{d},\nabla\mathbf{d})     \|_{L^4 }\\
\leqslant &\;
C
\left(
E_H+1
\right)
\left[
E_H^{\frac{3}{2}}\|\Delta\mathbf{v}\|
+E_H \|  \nabla(\Delta\mathbf{d}+f(\mathbf{d},\nabla\mathbf{d}) )  \| \|\Delta\mathbf{v}\|  
\right]^{\frac{1}{2}}\\
\leqslant &\;
C
\left(
E_H+1
\right)
\left[
\left(E_H+E_H^{\frac{1}{2}}\|\Delta\mathbf{v}\| \right)
+
E_H +\epsilon\|  \nabla(\Delta\mathbf{d}+f(\mathbf{d},\nabla\mathbf{d}) )  \| \|\Delta\mathbf{v}\|  
\right]\\
\leqslant &\;
C
\left(
E_H+1
\right)
E_H
+\epsilon\|\Delta\mathbf{v}\|^2
+\epsilon\|  \nabla(\Delta\mathbf{d}+f(\mathbf{d},\nabla\mathbf{d}) )  \|^2.
\end{aligned}
$$

Let $ K_4=\max\left\{|K_{32}|,|K_{42}|\right\} $,
and we have
$$
\begin{aligned}
K_4
\leqslant &\;
C\| \nabla\mathbf{d}  \|_{L^{\infty} }
\|\nabla^2\mathbf{v}\|
\| \Delta\mathbf{d}+f(\mathbf{d},\nabla\mathbf{d})  \|
\\
\leqslant &\;
C\left(1+E_H+ \|  \nabla(\Delta\mathbf{d}+f(\mathbf{d},\nabla\mathbf{d}) )  \|  \right)
\|\Delta\mathbf{v}\|
\| \Delta\mathbf{d}+f(\mathbf{d},\nabla\mathbf{d})  \|
\\
\leqslant &\;
C
\left(
E_H+1
\right)
E_H
+\epsilon\|\Delta\mathbf{v}\|^2
+\epsilon\|  \nabla(\Delta\mathbf{d}+f(\mathbf{d},\nabla\mathbf{d}) )  \|^2.
\end{aligned}
$$

For $ K_{52}$,
$$
\begin{aligned}
|K_{52}|
\leqslant &\;
C\int_Q
\left| \Delta\mathbf{d}+f(\mathbf{d},\nabla\mathbf{d})  \right|
\left( \left| \nabla^2\mathbf{v} \right| 
+\left| \nabla\mathbf{v} \right|   \left| \nabla\mathbf{d} \right|   \right)
 \left| \nabla\mathbf{d} \right| 
\mathrm{~d}x\\
\leqslant &\;
C\int_Q
\left| \Delta\mathbf{d}+f(\mathbf{d},\nabla\mathbf{d})  \right|
 \left| \nabla^2\mathbf{v} \right|
 \left| \nabla\mathbf{d} \right| 
\mathrm{~d}x
+C\int_Q
\left| \Delta\mathbf{d}+f(\mathbf{d},\nabla\mathbf{d})  \right|
 \left| \nabla\mathbf{v} \right|   \left| \nabla\mathbf{d} \right|^2   
\mathrm{~d}x\\
\leqslant &\;
C(K_4+|K_{21}|) 
\\
\leqslant &\;
C
\left(
E_H+1
\right)
E_H
+\epsilon\|\Delta\mathbf{v}\|^2
+\epsilon\|  \nabla(\Delta\mathbf{d}+f(\mathbf{d},\nabla\mathbf{d}) )  \|^2.
\end{aligned}
$$
By Theorem \ref{iso_higher-order_sigma_1}, after choosing $\epsilon_0$ and $\epsilon$ small enough,
we can prove the Theorem \ref{iso_Theorem_higher_order_final}. $\Box$

\paragraph{}
From the higher-order estimate in Theorem \ref{iso_Theorem_higher_order_final}, we obtain
\[
\frac{\mathrm{d}\ln (E_H(t) + 1)}{\mathrm{d}t} \leqslant C\;E_H(t).
\]
Thus,
\begin{equation}\label{3.70}
\begin{aligned}
\sup_{0 < t < T} E_H(t)
& \leqslant
(E_H(0) + 1) \exp \left(C \int_0^{T} E_H(s) \,\mathrm{d} s\right) - 1 \\
& \leqslant
(E_H(0) + 1) \exp \left(C \int_0^{\infty} E_H(s) \,\mathrm{d} s\right) - 1.
\end{aligned}
\end{equation}

Since $(E_H(0) + 1) \exp \left(C \int_0^{\infty} E_H(s) \,\mathrm{d} s\right) - 1$ is independent of $T$, 
then \eqref{3.70} holds for any $T \in (0, \infty)$.
Therefore, we obtain the uniform a priori estimate for $(\mathbf{v}_m, \mathbf{d}_m)$
\begin{equation}\label{3.71}
\sup_{t > 0} E_H(t)
\leqslant
(E_H(0) + 1) \exp \left(C \int_0^{\infty} E_H(s) \,\mathrm{d} s\right) - 1.
\end{equation}
Therefore,
\[
\int_0^{\infty} \|\Delta\mathbf{v}_m(s)\|^2 +
\left\|\nabla\left(\Delta\mathbf{d}_m + f(\mathbf{d}_m, \nabla\mathbf{d}_m)\right)(s)\right\|^2
\,\mathrm{d} s
\lesssim
1.
\]
The above uniform a priori estimate is independent of the parameter $m$ and time $t$.

\section{Global well-posedness of the general Ericksen--Leslie system}\label{section5}

\setcounter{equation}{0}
\numberwithin{equation}{section}

\subsection{Galerkin approximation}
Let $\left\{\phi_i\right\} \subset H$ with $\| \phi_i  \|=1$ be the eigenvectors of the Stokes operator in the periodic space with zero mean,
which means
$$
\Delta \phi_i+\nabla P_i=-\kappa_i \phi_i \;\text { in } Q, \quad
\nabla \cdot \phi_i=0 \; \text { in } Q,
\quad
\int_Q \phi_i(x) \mathrm{~d} x=0,
$$
where $P_i \in L_p^2(Q)$ and $0 < \kappa_1 
 \leqslant  \kappa_2\leqslant \cdots$  are eigenvalues\cite{Temam1995navier}. 
 The eigenvectors $\phi_i$ are
smooth and the sequence $\left\{\phi_i\right\}$ is an orthogonal basis of $H$. Denoted by
$$
P_m: H \longrightarrow H_m :=\operatorname{span}\left\{\phi_1, \phi_2, \ldots, \phi_m\right\},
$$
where $P_m$ is a orthonormal projection. 

Suppose that constant $T\in (0,\infty)$.
We consider the following approximate problem in 
 $Q_T=Q\times [0, T]$, 
\begin{equation}\label{iso_approximate_problem_d_first}
\left\{
\begin{aligned}
\partial_t\mathbf{v}_m
= & \; \mathrm{P}_m\left\{-\mathbf{v}_m \cdot \nabla \mathbf{v}_m
+\frac{\gamma}{\operatorname{Re}} \Delta \mathbf{v}_m
+\frac{1-\gamma}{\operatorname{Re}}\nabla\cdot\left[\sigma_{1, m}
+\sigma_{2, m}
+\sigma^E_m\right]\right\},  \\
 \nabla \cdot \mathbf{v}_m=& \; 0, \\
\partial_t \mathbf{d}_m
=& \; -\mathbf{v}_m \cdot \nabla \mathbf{d}_m
+\mu_1\left(\mathbf{h}_m+|\nabla \mathbf{d}_m|^2 \mathbf{d}_m-(\mathbf{H}\cdot\mathbf{d}_m)^2 \mathbf{d}_m\right)\\
& \;+(\mathbf{I}-\mathbf{d}_m \mathbf{d}_m)\left(\boldsymbol{\Omega}_m \mathbf{d}_m+\mu_2 \mathbf{D}_m \mathbf{d}_m\right), \\
|\mathbf{d}_m|=& \; 1,
\end{aligned}
\right.
\end{equation}
where
$$
\begin{aligned}
\boldsymbol{\Omega}_m= & \; \frac{1}{2}\left(\nabla \mathbf{v}_m-\nabla^T \mathbf{v}_m\right),  
\quad \mathbf{D}_m=  \frac{1}{2}\left(\nabla \mathbf{v}_m+\nabla^T \mathbf{v}_m\right),  \\
\mathbf{h}_m= & \; \Delta\mathbf{d}_m+\left(\mathbf{H} \cdot \mathbf{d}_m  \right) \mathbf{H},  
\quad  \sigma^E_m=  \nabla \mathbf{d}_m \odot \nabla \mathbf{d}_m, \\
\sigma_{1, m}  =& \;
\beta_1\left(\mathbf{D}_m: \mathbf{d}_m\otimes\mathbf{d}_m \right) \mathbf{d}_m\otimes\mathbf{d}_m
+\beta_2\mathbf{D}_m
+\frac{\beta_3}{2}\left[\mathbf{d}_m \otimes\left(\mathbf{D}_m \mathbf{d}_m\right)+\left(\mathbf{D}_m \mathbf{d}_m\right)\otimes\mathbf{d}_m\right],  \\
\sigma_{2, m} =& \; \frac{1}{2}\left(-1-\mu_2\right) \mathbf{h}_m\otimes\mathbf{d}_m+\frac{1}{2}\left(1-\mu_2\right) \mathbf{d}_m\otimes\mathbf{h}_m 
+\mu_2\left(\mathbf{h}_m\otimes\mathbf{d}_m\right)
\mathbf{d}_m\otimes\mathbf{d}_m.
\end{aligned}
$$
We have initial conditions
\begin{equation}\label{iso_approximate_initial_d}
    \mathbf{v}_m(x,  0)=\mathrm{P}_m \mathbf{v}_0(x), 
    \quad \mathbf{d}_m(x,  0)=\mathbf{d}_0(x), 
    \quad |\mathbf{d}_0|=1, 
    \quad x\in Q,
\end{equation}
and boundary conditions
\begin{equation}\label{iso_approximate_boundary_d}
    \mathbf{v}_m(x+\mathbf{e}_i,  t) = \mathbf{v}_m(x,  t),  \quad \mathbf{d}_m(x+\mathbf{e}_i,  t) = (-1)^{a_i} \mathbf{d}_m(x,  t),  \quad  \ (x, t) \in \partial Q \times \left[0, T\right].
\end{equation}

Suppose that we have a  solution $\left(\mathbf{v}_m, \mathbf{d}_m\right)$ in ODE system (\ref{iso_approximate_problem_d_first}). Let
\begin{equation}\label{iso_vm_def}
\mathbf{v}_m(x, t)=\sum_{i=1}^m g_m^i(t) \phi_i(x).
\end{equation}
Then the first equation in (\ref{iso_approximate_problem_d_first}) becomes
\begin{equation}\label{iso_g}
    \frac{\mathrm{d}}{\mathrm{d} t} g_m^i(t)
= \frac{\kappa_i\gamma}{\operatorname{Re}} g_m^i(t)
+\sum_{j,  k=1}^{2} A_{j k}^i g_m^k(t) g_m^j(t)
+\sum_{k=1}^{2} B_k^i g_m^k(t)
+C_m^i(t),
\end{equation}
where
$$
\begin{aligned}
A_{j k}^i  =& \; -\int_Q\left(\phi_j(x) \cdot \nabla \phi_k(x)\right) \phi_i(x) \;\mathrm{d} x, \\
B_k^i  =& \; \frac{\gamma-1}{\operatorname{Re}}\int_Q \left\{ \beta_1\left[\left(\frac{\nabla\phi_k+(\nabla\phi_k)^T}{2}\right):\mathbf{d}_m\otimes\mathbf{d}_m\right]^2
\right.\\
& \; +\beta_2\left[\frac{\nabla\phi_k+(\nabla\phi_k)^T}{2}\right]:\left[\frac{\nabla\phi_i+(\nabla\phi_i)^T}{2}\right]\\
& \; \left.+\beta_3\left[\mathbf{d}_m\cdot \left(\frac{\nabla\phi_k+(\nabla\phi_k)^T}{2}\right)\right] \cdot\left[\mathbf{d}_m\cdot \left(\frac{\nabla\phi_i+(\nabla\phi_i)^T}{2}\right)\right]\right\} \;\mathrm{d} x,  \\
C_m^i  =&-\frac{1-\gamma}{\operatorname{Re}}\int_Q \left\{\sum_{k,  l=1}^2
\left[\nabla_k \mathbf{d}_m \cdot \nabla_l \mathbf{d}_m
+\frac{1}{2}(1-\mu_2)d_m^k \left(\Delta d_m^l+f^l\left(\mathbf{d}_m,\nabla\mathbf{d}_m\right)\right)\right.\right.\\
&+\left.\left.\frac{1}{2}(-1-\mu_2) \cdot\left(\Delta d_m^k+f^k\left(\mathbf{d}_m, \nabla\mathbf{d}_m\right)\right)
d_m^l
\right]\nabla_l \phi_i^k\right\}  \;\mathrm{d} x .
\end{aligned}
$$
Here $f\left(\mathbf{d}_m, \nabla\mathbf{d}_m\right)
=\left(f^1\left(\mathbf{d}_m, \nabla\mathbf{d}_m\right), f^2\left(\mathbf{d}_m, \nabla\mathbf{d}_m\right)\right)^T$, and
$$\mathbf{d}_m=\left(d^1_m, d^2_m\right)^T,   \quad
\phi_i(x)=\left(\phi_i^1(x), \phi^2_i(x)\right)^T.
$$

The initial conditions  of the ODE system (\ref{iso_approximate_problem_d_first}) are
\begin{equation}\label{iso_gm_initial}
g_m^i(0)=\left(\mathbf{v}_0, \phi_i\right) \text { for } i=1,2, \ldots, m.    
\end{equation}

\begin{theorem}\label{iso_approxiamte_existence}
Suppose that $\left(\mathbf{v}_0, \mathbf{d}_0 \right)\in V\times H_{sp}^2$. 
For any $m>0$, there exists $T_0 \in (0,T)$ depending on $\left( \mathbf{v}_0, \mathbf{d}_0  \right)$ 
such that the  approximate problem  (\ref{iso_approximate_problem_d_first})
with initial conditions (\ref{iso_approximate_initial_d})
and boundary conditions (\ref{iso_approximate_boundary_d})
admits a unique strong solution $\left(\mathbf{v}_m, \mathbf{d}_m\right)$ in $Q_{T_0}$.

Moreover,  the solution $\left( \mathbf{v}_m, \mathbf{d}_m  \right)\in C^\infty(Q\times(0,T_0))$.
\end{theorem}

\paragraph{Proof.}
(i) 
First, 
for a given $\mathbf{v}_m$, solve the equation for $\mathbf{d}_m$.

Suppose that $\mathbf{v}_m$ is defined in (\ref{iso_vm_def}), satisfying the ODE system (\ref{iso_approximate_problem_d_first})
with initial conditions (\ref{iso_gm_initial}).
Assume that  $M$ is a large constant and $T_0\in (0, T)$ satisfies
$$
\left(\sum_{i=1}^m\left|g_m^i(t)\right|^2\right)^{1 / 2} \leqslant M, \quad
t \in\left[0,  T_0\right].
$$
Since $\mathbf{d}_0\in H^2_{sp}$, 
we have $\theta_0 \in H^2_{sq}$.
Consider the parabolic problem  about $\theta_m$ in $Q_{T_0}$
\begin{equation}\label{iso_approximate_theta_parabolic_problem_first}
\left\{
\begin{aligned}
\partial_t \theta_m
=& \;-\mathbf{v}_m \cdot \nabla \theta_m
+\mu_1\left(
\Delta\theta_m + \frac{\mathrm{H}^2}{2} \sin 2\theta_m \right)&\\
& \;
+\left(\boldsymbol{\Omega}_m:\mathbf{d}_m^{\perp}\otimes\mathbf{d}_m 
+\mu_2 \mathbf{D}_m:\mathbf{d}_m^{\perp}\otimes\mathbf{d}_m \right),  && 
(x, t)\in  Q\times\left[0,  T_0\right],  \\
\theta_m|_{t=0}=& \;
\theta_0, \;
 x\in Q ;\quad
\theta_m(x+\mathbf{e}_i, t)
=\theta_m(x, t)+a_i\pi,   \; &&
(x, t)\in \partial Q\times\left[0,  T_0\right].
\end{aligned}
\right.
\end{equation}

Assume that $\theta_\infty(x)$ is the solution to the elliptic
sine-Gorden equation (\ref{elliptic_first}), and denoted by $\widetilde{\theta}_m(x, t)=\theta_m(x, t)-\theta_\infty(x) $.
Then  $\widetilde{\theta}_m$ is a periodic function, and
we have
$$
\begin{aligned}
\mathbf{d}_m
=
\left(\sin\left( \widetilde{\theta}_m+\theta_\infty  \right), 
\cos\left( \widetilde{\theta}_m+\theta_\infty  \right)
\right)^T, \quad
\mathbf{d}_0 =\left(\sin\left(\theta_0+\theta_\infty\right), \cos\left(\theta_0+\theta_\infty\right)\right)^T.
\end{aligned}
$$
From (\ref{iso_approximate_theta_parabolic_problem_first}),
we can obtain the following parabolic equation
\begin{equation}\label{iso_approximate_theta_parabolic_problem_second2}
  \partial_t\widetilde{\theta}_m
+\mathbf{v}\cdot\nabla\widetilde{\theta}_m
=
\mu_1\Delta\widetilde{\theta}_m
+g\left(\widetilde{\theta}_m\right), 
\end{equation}
where
$$
g\left(\widetilde{\theta}_m\right)=
\frac{\mu_1 \mathrm{H}^2}{2}\sin2\left(\widetilde{\theta}_m+\theta_\infty\right)
+\left(
\boldsymbol{\Omega}_m:\mathbf{d}_m^{\perp}  \mathbf{d}_m
+\mu_2 \mathbf{D}_m:\mathbf{d}_m^{\perp}  \mathbf{d}_m\right)
+\mu_1\Delta\theta_\infty-\mathbf{v}_m\cdot\nabla\theta_\infty.
$$
Since
$\theta_\infty$ is a smooth function, 
then
$ g\left(\widetilde{\theta}_m\right) $ is bounded.

The existence and regularity of  
$\widetilde{\theta}_m$  are referred to  the literature literature \cite{Evans2010partial,Ladyzhenskaya1968linear}, which give detailed proof of well-posedness.  The solution  
$\widetilde{\theta}_m$
of the parabolic equation (\ref{iso_approximate_theta_parabolic_problem_second2}) satisfies
$$
\begin{gathered}
\widetilde{\theta}_m \in L^{\infty}\left(0,  T_0 ; H^2_{sq}\right) \cap L^2\left(0,T_0 ;  H_{sp}^3 \right),  \quad
\partial_t\widetilde{\theta}_m \in L^{\infty}\left(0,  T_0 ; L_{p}^2\right) \cap L^2\left(0,  T_0 ; H_{p}^1\right).
\end{gathered}
$$
Thus, the following parabolic problem for $\mathbf{d}_m \in \mathbb{S}^1$  admits a solution.
\begin{equation}\label{iso_approximate_d_parabolic_problem_second2}
\left\{
\begin{aligned}
\partial_t \mathbf{d}_m
+\mathbf{v}_m \cdot \nabla \mathbf{d}_m
=& \;
\mu_1\left(\mathbf{h}_m+|\nabla \mathbf{d}_m|^2 \mathbf{d}_m-(\mathbf{H}\cdot\mathbf{d}_m)^2 \mathbf{d}_m\right)\\
& \; +(\mathbf{I}-\mathbf{d}_m \mathbf{d}_m)\left(\boldsymbol{\Omega}_m \mathbf{d}_m+\mu_2 \mathbf{D}_m \mathbf{d}_m\right),  \;
&&
(x, t)\in  Q\times\left[0,  T_0\right], 
\\
\mathbf{d}|_{t=0}=\mathbf{d}_0, \quad
 x\in& \; Q ;\quad
\mathbf{d}_m(x+\mathbf{e}_i, t)
=(-1)^{a_i}\mathbf{d}_m(x, t),   \;
&& 
(x, t)\in \partial Q\times\left[0,  T_0\right].
\end{aligned}
\right.
\end{equation}
The solution $\mathbf{d}_m$ to (\ref{iso_approximate_d_parabolic_problem_second2})
satisfies
$$
\begin{gathered}
\mathbf{d}_m \in L^{\infty}\left(0,  T_0 ; H^2_{sp}\right) \cap L^{\infty}\left(0,  T_0 ; H_{sp}^3\right),  \quad
\partial_t\mathbf{d}_m \in L^{\infty}\left(0,  T_0 ; L_{sp}^2\right) \cap L^2\left(0,  T_0 ; H_{sp}^1\right).
\end{gathered}
$$

\paragraph{}
(ii)
Second,
derive upper bounds for the coefficients $C_m^i(t)$
in  (\ref{iso_g}), and use the Leray--schauder fixed-point theorem to establish the existence of solutions to the approximate problem (\ref{iso_approximate_problem_d_first}),
which is similar to the  Lemma 2.2 in \cite{Lin1995CPAM}.

Moreover, the solution $(\mathbf{v}_m,  \mathbf{d}_m)$ to (\ref{iso_approximate_problem_d_first}) 
 satisfies
$$
\begin{aligned}
 \mathbf{v}_m(x,  t) &\in  L^{\infty}\left(0,  T_0 ; V\right) \cap L^2\left(0,  T_0 ; H_{p}^2\right)
 \cap \mathrm{Lip}\big(Q_{T_0}\big),  
\\
\mathbf{d}_m (x, t) &\in L^{\infty}\left(0,  T_0 ; H_{sp}^2\right) \cap L^2\left(0,  T_0 ; H_{sp}^3\right),  
\end{aligned}
$$
where $\mathrm{Lip}\big(Q_{T_0}\big)$ is the 
space of Lipschitz continuous functions
in $Q_{T_0}$.
By the standard bootstrap argument in \cite{Ladyzhenskaya1968linear}
shows that $\mathbf{v}_m $  and  $  \mathbf{d}_m$
are smooth in the interior  of  $Q_{T_0}$,
the detailed proof of which can be seen in \cite{Lin1995CPAM,Lin2000ARMA}.

For any $x_1 \in \partial Q$, 
let 
$$
\begin{aligned}
&\widetilde{\mathbf{v}_m} (x)
=
\mathbf{v}_m\left(x-x_1+\frac{1}{2}\left(\mathbf{e}_1+\mathbf{e}_2\right)\right),  
\\
&\widetilde{\mathbf{d}_m} (x)
=
\mathbf{d}_m\left(x-x_1+\frac{1}{2}\left(\mathbf{e}_1+\mathbf{e}_2\right)\right),
\\
\widetilde{Q}
&=\left\{ y\mid y= x-x_1+\frac{1}{2}\left(\mathbf{e}_1+\mathbf{e}_2\right),  \; x\in Q  \right\}.
\end{aligned}
$$
After the translation transformation,  
$(\widetilde{\mathbf{v}_m}, \widetilde{\mathbf{d}_m})$ still satisfies the initial conditions  (\ref{iso_approximate_initial_d})  and the boundary conditions (\ref{iso_approximate_boundary_d}), and $x_1$ is in the interior of $\widetilde{Q}$.
Therefore,
$\left(\mathbf{v}_m, \mathbf{d}_m\right)$ is smooth in  $Q\times(0,T_0)$,
which means $\left( \mathbf{v}_m, \mathbf{d}_m  \right)\in C^\infty(Q\times(0,T_0))$.

\paragraph{}
(iii)
Third, we consider the uniqueness of the solution $(\mathbf{v}, \mathbf{d})$  to (\ref{iso_approximate_problem_d_first})   in $Q_{T_0}$.

Assume that $( \mathbf{v}_{mj}, \theta_{mj})$, $j=1, 2$
are solutions to the following approximate problem in $Q\times\left[0,  T_0\right]$
with initial values
$( \mathbf{v}_{0j}, \theta_{0j})\in V\times H^2_{sq}$, $j=1, 2$,$Q\times\left[0,  T_0\right]$
\begin{equation}\label{unique1}
\left\{
\begin{aligned}
\partial_t\mathbf{v}_m
= & \; \mathrm{P}_m\left\{-\mathbf{v}_m \cdot \nabla \mathbf{v}_m
+\frac{\gamma}{\operatorname{Re}} \Delta \mathbf{v}_m
+\frac{1-\gamma}{\operatorname{Re}}\nabla\cdot\left[\sigma_{1, m}
+\sigma_{2, m}
+\sigma^E_m\right]\right\},  \\
 \nabla \cdot \mathbf{v}_m=& \; 0, \\
\partial_t \theta_m
=& \;-\mathbf{v}_m \cdot \nabla \theta_m
+\mu_1\left(
\Delta\theta_m + \frac{\mathrm{H}^2}{2} \sin 2\theta_m \right)\\
& \;
+\left(\boldsymbol{\Omega}_m:\mathbf{d}_m^{\perp}\otimes\mathbf{d}_m 
+\mu_2 \mathbf{D}_m:\mathbf{d}_m^{\perp}\otimes\mathbf{d}_m \right).
\end{aligned}
\right.
\end{equation}
The initial conditions are
$$
\mathbf{v}_m\big|_{t=0}=\mathbf{v}_0,\quad
\nabla\cdot \mathbf{v}_0=0, \quad
\theta_{mj}\big|_{t=0}=\theta_{0j},\quad
 x\in Q \times\left[0,  T_0\right].
$$
The boundary conditions are
$$
\mathbf{v}_m(x+\mathbf{e}_i,t)
=\mathbf{v}(x,t),\quad
\theta_{mj}(x+\mathbf{e}_i,t)
=\theta_{mj}(x,t)+a_i\pi,  \quad
(x,t)\in \partial Q\times\left[0,  T_0\right].
$$
By calculations similar to those in Lemma 2.2 of \cite{Wu2012CVPDE}, using \eqref{unique1}, there exist positive constants $\epsilon$ and $C$
such that
$$
\begin{aligned}
&\frac{\mathrm{d}}{\mathrm{d}t} 
\left( 
\left\| \mathbf{v}_1-\mathbf{v}_2 \right\|^2 
+
\left\|  \theta_1-\theta_2   \right\|^2_{H^1}
\right)
 +
\epsilon 
\left( 
\left\| \nabla\left(\mathbf{v}_1-\mathbf{v}_2\right)  \right\|^2 
+
\left\| \Delta\left( \theta_1-\theta_2\right)  \right\|^2
\right)\\
\leqslant
& \;
C \left( 
\left\| \mathbf{v}_1-\mathbf{v}_2  \right\|^2 
+
\left\|   \theta_1-\theta_2 \right\|^2_{H^1}
\right).
\end{aligned}
$$
Applying Gronwall's inequality, for any $t\in \left[0,  T_0\right]$, the following inequality holds,
$$
\begin{aligned}
& 
\left\| \left(\mathbf{v}_1-\mathbf{v}_2\right)(t)  \right\|^2 
+
\left\|  \left(\theta_1-\theta_2\right)(t)   \right\|^2_{H^1}
+
\epsilon
\int_0^t
\left( 
\left\| \nabla\left(\mathbf{v}_1-\mathbf{v}_2\right)(s)  \right\|^2 
+
\left\| \Delta\left( \theta_1-\theta_2\right) (s)  \right\|^2 
\right) \mathrm{d}s\\
\leqslant
& \;
2 e^{Ct} \left( \left\| \left( \mathbf{v}_1-\mathbf{v}_2\right)(0)  \right\|^2 
+
\left\|  \left( \theta_1-\theta_2\right)(0)   \right\|^2_{H^1}
\right).
\end{aligned}
$$
The constant $C$  depends on $Q,  \mathbf{H}, \left\|\mathbf{d}_{0j}\right\|_{H^2}$,  
$\left\|\mathbf{v}_{0j}\right\|_{H^1}$ and
coefficients of the system (\ref{iso_Ericksen_Leslie_new_formulation}),
and  $C$ is independent of  $t$.
Thus, we get the uniqueness of  the solution $(\mathbf{v}, \theta)$ to  (\ref{unique1}) in $Q_{T_0}$. 
Therefore, we obtain the uniqueness of the solution $(\mathbf{v}, \mathbf{d})$  to (\ref{iso_approximate_problem_d_first})   in $Q_{T_0}$.
$\Box$

\begin{proposition}\label{prop3.2}
For  $T_0\in (0, T)$, 
suppose that $(\mathbf{v}_m,  \mathbf{d}_m)$ is a solution to the  approximate problem (\ref{iso_approximate_problem_d_first})--(\ref{iso_approximate_boundary_d})
in $Q_{T_0}$. 
Then $(\mathbf{v}_m,  \mathbf{d}_m)$  is also a solution to the approximate problem  (\ref{iso_approximate_problem_d_first})--(\ref{iso_approximate_boundary_d}) in  $Q_{T}$.
\end{proposition}

The proof of \ref{prop3.2} refers to the proof of Theorem 2.1 in \cite{Lin1995CPAM}.
If $(\mathbf{v}_m,  \mathbf{d}_m)$ is a solution in $Q_{T_0}$,
by the basic energy law (Theorem \ref{iso_basic_energy_law}) and  the extension theorem for solutions of ordinary differential equations, the solution to the approximate problem exists on $Q_T$.

\subsection{Well-posedness of the global strong solution}
For the initial value
$(\mathbf{v}_0, \mathbf{d}_0) \in V \times H_{sp}^2$, 
we need to prove
the general Ericksen--Leslie system  (\ref{iso_Ericksen_Leslie_new_formulation}) with initial conditions (\ref{initial_condition_1}) and boundary conditions (\ref{boundary_condition_1})
 admits a global strong solution $(\mathbf{v}, \mathbf{d})$ in $Q \times [0, \infty)$ satisfying
\[
\begin{aligned}
& \mathbf{v} \in L^{\infty}(0, \infty; V) \cap L^2(0, \infty; H_p^2), \\
& \mathbf{d} \in L^{\infty}(0, \infty; H_{sp}^2) \cap L_{\mathrm{loc}}^2(0, \infty; H_{sp}^3).
\end{aligned}
\]

\paragraph{Proof of Theorem 1.2.}
(i)
The existence of the global strong solution to the general system  follows the idea in Section 4.3 of \cite{WXL2013ARMA}.

From the uniform a priori estimate \eqref{3.71}, by taking $m \rightarrow \infty$, we obtain a local solution $(\mathbf{v}, \mathbf{d})$ of (\ref{iso_Ericksen_Leslie_new_formulation})  in $Q_T$ satisfying
\[
\begin{aligned}
& \mathbf{v} \in L^{\infty}(0, T; V) \cap L^2(0, T; H_p^2), \\
& \mathbf{d} \in L^{\infty}(0, T; H_{sp}^2) \cap L^2(0, T; H_{sp}^3).
\end{aligned}
\]

Taking $T \rightarrow \infty$, the local solution $(\mathbf{v}, \mathbf{d})$ is extended to be a global strong solution of (\ref{iso_Ericksen_Leslie_new_formulation})  satisfying
\[
\begin{aligned}
& \mathbf{v} \in L^{\infty}(0, \infty; V) \cap L^2(0, \infty; H_p^2), \\
& \mathbf{d} \in L^{\infty}(0, \infty; H_{sp}^2) \cap L_{\mathrm{loc}}^2(0, \infty; H_{sp}^3).
\end{aligned}
\]
Thus, we obtain $(\mathbf{v}, \mathbf{d})$ is a global strong solution of (\ref{iso_Ericksen_Leslie_new_formulation}) .

\paragraph{}
(ii) The regularity of the global strong solution to the general
system follows the idea in 
 \cite{Lin2014CMP}, which uses  Serrin's method in \cite{Serrin1962bootstrap}.
 
It is clear that $(\mathbf{v}, \mathbf{d})$ is smooth in the interior of $Q \times (0, \infty)$. Next, we can obtain that $(\mathbf{v}, \mathbf{d})$ is smooth in $Q \times (0, \infty)$ by a translation transformation,
which can be seen in the proof of Theorem \ref{iso_approxiamte_existence}.
Therefore, we can obatin $(\mathbf{v}, \mathbf{d})$ is smooth in $Q \times (0, \infty)$, which means $\left( \mathbf{v}, \mathbf{d}  \right)\in C^\infty(Q\times(0,\infty))$.

\paragraph{}
(iii)
The uniqueness of the solution follows the same process as in Theorem \ref{iso_approxiamte_existence}. 

For any $T \in (0, \infty)$, we obtain the uniqueness of $(\mathbf{v}, \theta)$ in $Q \times [0, T]$, which leads to the uniqueness of $(\mathbf{v}, \mathbf{d})$ in $Q \times [0, T)$.

In summary, the well-posedness of the global strong solution is proved. $\Box$

\begin{remark}
We can also get the well-posedness of the  global weak solution to 
the general Ericksen--Leslie system  (\ref{iso_Ericksen_Leslie_new_formulation}) with initial conditions (\ref{initial_condition_1}) and boundary conditions (\ref{boundary_condition_1}) by the methods in \cite{Lin1995CPAM,Lin2000ARMA}.
\end{remark}

\section{Long-time behavior}

\setcounter{equation}{0}
\numberwithin{equation}{section}

\subsection{\texorpdfstring{$\omega$-limit set }{Omega-limit set}}
For $(\mathbf{v}_0,  \mathbf{d}_0) \in V\times H_{sp}^2$, we can define $\omega$-limit set 
with respect to the initial value and analysis the decay property.

\begin{lemma}\label{iso_lem_Omega_limit}
Suppose that
$(\mathbf{v}_0,  \mathbf{d}_0) \in V\times H_{sp}^2$, 
then the global strong solution to 
(\ref{iso_Ericksen_Leslie_new_formulation}) with initial conditions (\ref{initial_condition_1}) and boundary conditions (\ref{boundary_condition_1})
satisfies
\begin{equation}\label{iso_lem_Omega_limit_equation}
\lim _{t \rightarrow \infty}
\left\{
\|\mathbf{v}(t)\|_{H^1}
+\left\|\left(
\Delta\mathbf{d}+f\left(\mathbf{d}, \nabla\mathbf{d}\right) 
\right)(t)\right\|\right\}
=0 .
\end{equation}
Additionally, we get the uniform boundedness regarding  $t\geqslant 0$,
\begin{equation}\label{iso_lem_Omega_limit_inequation_1}
\sup _{t \geqslant 0}
\left\{\|\mathbf{v}(t)\|_{H^1}+\|\mathbf{d}(t)\|_{H^2}\right\} \leqslant C, 
\end{equation}
and
\begin{equation}\label{iso_lem_Omega_limit_inequation_2}
			\sup _{t \geqslant 0}\left\{\|\mathbf{v}(t)\|_{H^1_p}+\|\nabla\theta(t)\|_{H^1_p}\right\} \leqslant C.
\end{equation}
The constant $C$  depends on $Q, |\mathbf{H}|,\left\|\mathbf{d}_0\right\|_{H^2}$, 
$\left\|\mathbf{v}_0\right\|_{H^1}$
and
coefficients of the system (\ref{iso_Ericksen_Leslie_new_formulation}).
\end{lemma}

\paragraph{Proof.}
Taking \( y(t) = E_H(t) \) in in Lemma 6.2.1 of \cite{Zheng2004nonlinear},  
by higher-order estimates, we can derive  
\[
\lim _{t \rightarrow \infty} E_H(t)
=
\lim _{t \rightarrow \infty}
\left\{\|\nabla \mathbf{v}(t) \|^2
+
\frac{1-\gamma}{\operatorname{Re}}
\left\| \Delta \theta(t)+\frac{\mathrm{H}^2}{2} \sin 2 \theta(t) \right\| ^2\right\}
=0.
\]
Thus, (\ref{iso_lem_Omega_limit_inequation_1}) is proved.

Define  
\[
\mathcal{E}(t)
:=
\int_Q
|\mathbf{v}|^2
+\frac{1-\gamma}{\operatorname{Re}}
(|\nabla \mathbf{d}|^2 +|\mathbf{H} |^2 -(\mathbf{H}\cdot\mathbf{d})^2 ) 
\, \mathrm{d} x.
\]
For any \( t > 0 \), by  basic energy law (\ref{iso_basic_energy_law}), we have  
\[
\begin{aligned}
\mathcal{E}(t)
+
2
\int_0^t \int_Q
\frac{\gamma}{\operatorname{Re}}
|\nabla \mathbf{v}|^2
+\frac{(1-\gamma)\mu_1}{\operatorname{Re}}
\left(|\mathbf{h}|^2
-(\mathbf{h} \cdot \mathbf{d})^2
\right)
\, \mathrm{d} x \, \mathrm{d} s
\leqslant 
\mathcal{E}(0).   
\end{aligned}
\]
Therefore, we can obtain  
\[
\sup _{t \geqslant 0}
\left\{\|\mathbf{v}(t)\|^2
+
\|\mathbf{d}(t)\|_{H^1}^2\right\}
+
\int_0^{\infty} E_H(s) \, \mathrm{d} s \leqslant C.
\]
Since \( E_H \) is bounded, we have  
\[
\|\nabla\mathbf{v}\|
+\frac{1-\gamma}{\operatorname{Re}}
\left\| \Delta \theta+\frac{\mathrm{H}^2}{2} \sin 2 \theta \right\| 
\leqslant C.
\]
This allows us to derive an upper bound estimate for \( \|\Delta \theta\| \):  
\[
\|\Delta \theta\| \leqslant \|\Delta\theta 
+\dfrac{\mathrm{H}^2}{2}\sin(2\theta)  \|
+ \frac{\mathrm{H}^2}{2}
\leqslant 
C.
\]
Thus, the inequality (\ref{iso_lem_Omega_limit_inequation_2}) is proved.

Thanks to 
$$
\|\Delta\mathbf{d}\|
\leqslant
C( \|\Delta \theta\|+1 ), \quad
|\nabla\mathbf{d}|^2=|\nabla\theta|^2, 
$$
then  the inequality (\ref{iso_lem_Omega_limit_inequation_1}) is proved. $\Box$

\paragraph{}
The $\omega$-limit sets of 
$\left(\mathbf{v}_0, \mathbf{d}_0\right) \in V \times H_{sp}^2$ 
can be defined as follows,
\begin{equation}
\begin{aligned}
 \omega\left(\mathbf{v}_0, \mathbf{d}_0\right)
:=
\big\{
\left(\mathbf{v}_{\infty}, \mathbf{d}_{\infty}\right)
\mid &
\text { there exists a sequence }
\left\{t_n\right\} \nearrow \infty
\text { such that } \\
&
\left(\mathbf{v}\left(\cdot,t_n\right), \mathbf{d}\left(\cdot,t_n\right)\right) 
\rightarrow
\left(\mathbf{v}_{\infty}, \mathbf{d}_{\infty}\right) 
\text { in } V \times H_{sp}^2, 
\text { as } t_n \rightarrow\infty\big\},
\end{aligned}  
\end{equation}
where $\mathbf{d}_\infty=(\sin\theta_\infty, \cos\theta_\infty)^T $ and 
 $\theta_\infty$ satisfies (\ref{elliptic_first}).

\begin{theorem}\label{iso_omega_limit_set_property}
Suppose that $(\mathbf{v}_0, \mathbf{d}_0) \in V \times H^2_{sp}$,  
then the set  
$\omega\left(\mathbf{v}_0, \mathbf{d}_0\right)$ 
is a nonempty bounded subset of $V \times H_{sp}^2$.

Furthermore, all asymptotic limit points  
$(\mathbf{v}_{\infty}, \mathbf{d}_{\infty})$  
of the system (\ref{iso_Ericksen_Leslie_new_formulation}) with  (\ref{initial_condition_1}) and  (\ref{boundary_condition_1}) satisfy  
$(\mathbf{v}_{\infty}, \mathbf{d}_{\infty}) \in \Sigma$,  
where  
$$
\begin{aligned}
\Sigma := \left\{ (0, \mathbf{d}) \mid \mathbf{d} = (\sin\theta_{\infty}, \cos\theta_{\infty})^T, \; \theta_{\infty} \text{ is a solution to (\ref{elliptic_first})} \right\},
\end{aligned}
$$  
which means $\omega\left(\mathbf{v}_0, \mathbf{d}_0\right) \subset \Sigma$.
\end{theorem}

\paragraph{Proof.}
The proof is inspired by the Lemma 3.6 in \cite{CKY2018SIAM}.

By Lemma \ref{iso_lem_Omega_limit}, the global solution $(\mathbf{v}, \mathbf{d})$ has a uniform upper bound in $H_p^1(Q) \times H_{sp}^2(Q)$ for $t \geqslant 0$.  
Applying  Rellich--Kondrachov compactness theorem, we conclude that $\omega\left(\mathbf{v}_0, \mathbf{d}_0\right)$ is a nonempty bounded subset of $L^2 \times H_{sp}^1(Q)$.  
From (\ref{iso_lem_Omega_limit_equation}), we know that as $t \rightarrow \infty$, $\mathbf{v}(\cdot, t) \rightarrow 0$ in $V$.

Since $H^2_{sp} \hookrightarrow H^1_{sp}$ is a compact embedding, there exists a monotonically increasing unbounded sequence $\{t_n\}$ and a function $\mathbf{d}_{\infty}$ such that as $t_n \rightarrow \infty$,  
$$
\lim_{n \rightarrow \infty} \left\| \mathbf{d}(t_n) - \mathbf{d}_{\infty} \right\|_{H^1} = 0.
$$  
Let $\mathbf{d}_{\infty}(x) = (\sin\theta_{\infty}(x), \cos\theta_{\infty}(x))^T$ and $\mathbf{d}(x, t_n) = (\sin\theta(x, t_n), \cos\theta(x, t_n))^T$.  
By (\ref{3.71}), $\left\| \Delta\theta(t) + \frac{\mathrm{H}^2}{2} \sin 2\theta(t) \right\|$ has a uniform upper bound.  
Applying  Banach--Alaoglu theorem, the sequence $\left\{ \Delta\theta(t_n) + \frac{\mathrm{H}^2}{2} \sin 2\theta(t_n) \right\}$ has a subsequence (still denoted as the original sequence) that weakly converges to $0$ in $L^2_p$.  
Thus, $\theta_{\infty}$ satisfies  
$$
- \Delta \theta_{\infty} = \frac{\mathrm{H}^2}{2} \sin 2 \theta_{\infty} \quad \text{in } Q.
$$  
Since $\mathbf{d}(x, t)$ and $\mathbf{d}_0$ have the same boundary conditions, the smooth function $\theta_{\infty}$ satisfies the boundary conditions in (\ref{elliptic_first}).

From  
$$
\lim_{n \rightarrow \infty} \left\| \mathbf{d}(t_n) - \mathbf{d}_{\infty} \right\|_{H^1} = 0,
$$  
we deduce that as $n \rightarrow \infty$,  
$$
\left\| \sin 2\theta(t_n) - \sin 2\theta_{\infty} \right\|_{H^1} \rightarrow 0.
$$  

Since  
$$
\begin{aligned}
\left\| \Delta\theta(t_n) - \Delta\theta_{\infty} \right\| \leqslant &  
\left\| \Delta\theta(t_n) + \frac{\mathrm{H}^2}{2} \sin 2\theta(t_n) - \left( \Delta\theta_{\infty} + \frac{\mathrm{H}^2}{2} \sin 2\theta_{\infty} \right) \right\| \\
& + \frac{\mathrm{H}^2}{2} \left\| \sin 2\theta(t_n) - \sin 2\theta_{\infty} \right\|,
\end{aligned}
$$  
we conclude that  
$$
(0, \mathbf{d}_{\infty}) \in \omega\left(\mathbf{v}_0, \mathbf{d}_0\right), \quad \text{and} \quad \lim_{n \rightarrow \infty} \left\| \mathbf{d}(t_n) - \mathbf{d}_{\infty} \right\|_{H^2} = 0.
$$

\subsection{\texorpdfstring{The boundness of $\theta$ }{ The boundness of theta}}
Having obtained the existence of the global  strong solution $\left(\mathbf{v}, \mathbf{d}\right)$  and $\omega$-limit set $\omega\left(\mathbf{v}_0, \mathbf{d}_0\right)$, 
we still cannot derive a uniform estimate for $\theta$ with respect to time $t$,
because we do not have the maximum principle for $\theta$.
However,
the boundness of $\theta$ can be obtained by the relationship between $\mathbf{d}$ and $\theta$ and \L ojasiewicz--Simon inequality which can be applied to  $\theta$,  but not to $\mathbf{d}$.

\begin{lemma}\label{iso_lem_boundness}
Suppose that $(\mathbf{v}_0, \mathbf{d}_0) \in V \times H^2_{sp}$,  
then  
$\mathbf{d} \in C((0, \infty) \times Q)$.
\end{lemma}

\paragraph{Proof.}
For any \( T \in (0, \infty) \), when \( t \in (0, T) \),  
$$
\begin{aligned}
\left\|\partial_t \theta\right\|
\leqslant &
\left\|\mathbf{v} \cdot \nabla \theta\right\|
+
\left|\mu_1\right| \cdot \left\|\Delta \theta + \frac{H^2}{2} \sin 2\theta  \right\|
+
\left\|\left(\boldsymbol{\Omega}+\mu_2\mathbf{D} \right) :\mathbf{d}^{\perp}\otimes\mathbf{ d}\right\| 
\\
\leqslant &
\| \mathbf{v}  \|^2_{L^4} \| \nabla\theta \|^2_{L^4}
+C\left\|\Delta \theta + \frac{H^2}{2} \sin 2\theta  \right\|
+C\| \nabla\mathbf{v} \|.
\end{aligned}
$$
Since \( \|\Delta\theta\| \leqslant C \),  
we have  
$$
\| \nabla\theta \|^2_{L^4} 
\leqslant
C\|\nabla\theta\|
\left(
\left\|\nabla^2\theta\right\|+\|\nabla\theta\|
\right)
\leqslant
C
\left(
\|\Delta\theta\|^2+1
\right)
\leqslant C.
$$
By the Ladyzhenskaya inequality and Poincaré inequality,  
we can derive  
$$
\| \mathbf{v}  \|^2_{L^4}  
\leqslant C \left(\| \nabla \mathbf{v}\|+ \|  \mathbf{v}  \| \right)
\|  \mathbf{v}  \|
\leqslant C \| \nabla \mathbf{v}  \| \|  \mathbf{v}  \|
\leqslant C \| \nabla \mathbf{v}  \|.
$$ 
In summary,  
$$
\left\|\partial_t\theta\right\|
=
\left\|\partial_t \theta\right\|
\leqslant
C\left\|\Delta \theta + \frac{H^2}{2} \sin 2\theta  \right\|
+C\| \nabla\mathbf{v} \|.
$$
Thus, \( \partial_t\mathbf{d} \in L^{\infty}\left(0, \infty;L^2\right) \).

Since \( \mathbf{d} \in L^{\infty}\left(0, \infty;H^2\right) \),  
we can conclude that  
\( \mathbf{d} \in C((0, \infty)\times Q) \hookrightarrow C\left(0, \infty;H^1\right) \).  
Additionally, the corresponding \( \theta \) belongs to \( C_{loc}(\left[0, \infty\right)\times Q) \).  
$\Box$

\paragraph{}
If $\theta \in H^2_{sq}$, then $\sin2\theta\in H^2_{p}$  is a periodic function,
and both $\Delta\theta\in L^2_p$ and $\partial_t \theta\in L^2_p$
are also periodic functions.
Moreover, function $(\theta-\theta_\infty)$ 
is periodic.

We introduce the following \L ojasiewicz--Simon inequality with respect to $\theta$. 

\begin{lemma}\label{iso_lem_Lojasiewicz--Simon}
Suppose that $\theta_{\infty}$ is a solution to the elliptic problem (\ref{elliptic_first}) and is a critical point of the functional  
$$
E(\theta)=\int_Q \frac{1}{2} |\nabla\theta|^2+\frac{\mathrm{H}^2}{4}(1+\cos 2\theta) \, \mathrm{d}x. 
$$  
Then there exist constants $\rho_0 \in (0,1)$, $\eta > 0$, and $\nu \in \left(0, \frac{1}{2}\right]$, such that for any function $\theta \in H^1_{sq}$ satisfying $\left\|\theta - \theta_{\infty}\right\|_{H^1_p} < \rho_0$, the following inequality holds:  
$$
\left\|\Delta\theta + \frac{H^2}{2}\sin 2\theta\right\|_{H^{-1}} \geqslant \eta \left|E(\theta) - E(\theta_{\infty})\right|^{1 - \nu}.
$$  
Denote the dual space of $H^1_p$ by $H^{-1}$. The constants $\rho_0$, $\eta$, and $\nu$ depend on $Q$, $\theta_{\infty}$,  $ |\mathbf{H}|$, $\left\|\mathbf{d}_0\right\|_{H^2}$,  
$\left\|\mathbf{v}_0\right\|_{H^1}$
and
coefficients of the system (\ref{iso_Ericksen_Leslie_new_formulation}).
\end{lemma}

\begin{theorem}\label{iso_thm_boundness_theta}
Suppose that $\left(0, \mathbf{d}_{\infty}\right) \in \omega\left(\mathbf{v}_0, \mathbf{d}_0\right)$, and $\theta$  corresponds to the vector field $\mathbf{d}$. For any constant $\rho \in \left(0, \rho_0\right)$, there exists a constant $T_1 > 0$ and a smooth function $\theta_\infty$ corresponding to the vector field $\mathbf{d}_{\infty}$ such that for all $t \in \left[ T_1, \infty\right)$,
\begin{equation}\label{iso_inequality_boundness_first_1_theta}
\left\|\theta - \theta_{\infty}\right\|_{H^1} < \rho.
\end{equation}
Therefore, function $\theta$ is bounded in $Q \times (0, \infty)$. Additionally, there exists a constant $T_2 > 0$ such that for all $t \in \left[ T_2, \infty\right)$,
\begin{equation}\label{iso_inequality_boundness_first_2_d}
\left\|\mathbf{d} - \mathbf{d}_{\infty}\right\|_{H^1} < \rho.
\end{equation}
\end{theorem}

\paragraph{Proof.}
(i)
For $\mathbf{d}_\infty=(\sin\theta_{\infty}, \cos\theta_{\infty})^T$ in $\omega\left(\mathbf{v}_0, \mathbf{d}_0\right)$,  
we need to prove that
there exist a sequence of positive real numbers $\{t_n\}\nearrow \infty$ and a sequence of integers $\{l_n\}$ such that  
\begin{equation}\label{iso_equ_theta_d_limit_first}
\lim _{n\rightarrow \infty}
\left\| \theta(t_n)-\theta_{\infty}-2 l_n \pi \right\|_{H^1}=0.
\end{equation} 

Since $\left(0, \mathbf{d}_{\infty}\right) \in \omega\left(\mathbf{v}_0, \mathbf{d}_0\right)$, 
there exists a sequence of times $\{t_n\} \nearrow \infty$ such that  
$$
\lim _{n\rightarrow \infty}
\left\| \mathbf{d}(t_n)-\mathbf{d}_{\infty}  \right\|_{H^2}=0.    
$$
By the embedding \( H^2 \hookrightarrow L^\infty \),  
we have  
$$
\lim _{n\rightarrow \infty}
\left\| \mathbf{d}(t_n)-\mathbf{d}_{\infty}  \right\|_{L^\infty}=0.    
$$
Since $\mathbf{d}(x, t) \in C((0, \infty)\times Q)$ and $\mathbf{d}_{\infty}(x)$ is smooth,  
the above convergence also holds in the norm of continuous functions.  
Thus, as \( n \rightarrow \infty \),  
$$
\mathbf{d}(t_n) \rightrightarrows \mathbf{d}_{\infty}, \quad x \in Q.
$$
Therefore, for any \( 0 < \epsilon \ll 1 \),  
there exists a constant \( T(\epsilon) \geqslant 0 \),  
such that for all \( t_n \geqslant T(\epsilon) \),  
$$
\left|\mathbf{d}(x, t_n)-\mathbf{d}_{\infty}(x)\right|<\epsilon, \quad x \in Q.  
$$
Directly computing \( \left|\mathbf{d}(x, t_n)-\mathbf{d}_{\infty}(x)\right| \),  
we obtain  
$$
\begin{aligned}
\left|\mathbf{d}(x, t_n)-\mathbf{d}_{\infty}(x)\right|^2
=&\left(\sin\theta(x, t_n)-\sin\theta_\infty\right)^2+\left(\cos\theta(x, t_n)-\cos\theta_\infty\right)^2\\
=&4\sin^2\left(\frac{\theta(x, t_n)-\theta_\infty}{2}\right).    
\end{aligned}
$$
For \( \epsilon \ll 1 \),  
there exists a sequence of positive integers \(\{ l_n \}\),  
such that for all \( t_n \geqslant T(\epsilon) \),  
$$
\left|\theta(x, t_n)-\theta_\infty -2 l_n \pi \right| <2 \arcsin\left(\frac{\epsilon}{2}\right), \quad x \in Q.
$$
Thus, we can obtain
$$
\lim _{n\rightarrow \infty}
\left|\theta(x, t_n)-\theta_\infty -2 l_n \pi \right|=0.
$$
On the other hand,  
since \( |\nabla\mathbf{d}|^2=|\nabla\theta|^2 \),  
we have  
$$
\left\| \nabla\left(\mathbf{d}(t_n)-\mathbf{d}_{\infty}  \right)\right\| 
=
\left\| \nabla\left(\theta(t_n)-\theta_{\infty} \right)\right\| 
=
\left\| \nabla\left(\theta(t_n)-\theta_{\infty} -2 l_n \pi  \right)\right\|.
$$

\paragraph{}
(ii) For convenience,  denote $(\theta_\infty + 2 l_n \pi)$ by $\theta_\infty(n)$. 
Define
$$
E_1(\theta) := \int_Q \frac{1-\gamma}{\operatorname{Re}} \left( |\nabla \theta|^2 + |\mathbf{H}|^2 - |\mathbf{H}|^2 \sin^2\theta \right) \, \mathrm{d}x.
$$
Thus, $E_1(\theta_\infty(n)) = E_1(\theta_\infty)$, and there exists a sequence of times $\{t_n\}$ such that
$$
\lim_{n \rightarrow \infty} \mathcal{E}(t_n) = E_1(\theta_\infty).
$$
By energy dissipation, for all $t \geqslant 0$,
$$
E_1(\theta_\infty) \leqslant \mathcal{E}(t).
$$
If there exists a constant $T \geqslant 0$ such that $\mathcal{E}(T) = E_1(\theta_\infty)$, then for all $t > T$,
$$
\left\|\nabla \mathbf{v}(t)\right\| + \left\|\Delta \theta(t) + \frac{\mathrm{H}^2}{2} \sin 2 \theta(t)\right\| = 0.
$$
Thus, there exists an integer constant $k$ such that for all $(x, t) \in Q \times (T, \infty)$,
$$
\theta(x, t) = \theta_{\infty}(x) + 2k\pi, \quad \mathbf{v}(x, t) = 0.
$$

If for all $t \geqslant 0$, we have $\mathcal{E}(t) > E_1(\theta_\infty)$, then there exists a positive integer $N$ such that for all $n \geqslant N$,
\begin{equation}\label{iso_ine_thm3.10}
\left\| \theta(t_n) - \theta_{\infty}(n) \right\|_{H^1_p} \leqslant \frac{\rho}{2}.
\end{equation}
For integers $n \geqslant N$, define
$$
\bar{t}_n := \sup \left\{ t \, \bigg| \, t \geqslant t_n, \, \left\| \theta(t_n) - \theta_{\infty}(n) \right\|_{H^1_p} \leqslant \frac{\rho}{2}, \, \text{for all } s \in [t_n, t] \right\}.
$$
It follows that $\bar{t}_n > t_n$.

If there exists $N_1 \geqslant N$ such that $\bar{t}_{N_1} = \infty$, then for $t > t_{N_1}$, inequality (\ref{iso_ine_thm3.10}) holds.
We proceed by a contradiction to show that there must exist $N_1 \geqslant N$ such that $\bar{t}_{N_1} = \infty$. 
Suppose that for all $n \geqslant N$, the time $\bar{t}_n$ is finite, so for all $n \geqslant N$, we have $t_n < \bar{t}_n < \infty$.
By the \L ojasiewicz--Simon inequality, for all $t \in [t_n, \bar{t}_n]$, we have
$$
\begin{aligned}
\eta \left( E_1(\theta(t)) - E_1(\theta_\infty) \right)^{1 - \nu} \leqslant \left\| \Delta \theta(t) + \frac{\mathrm{H}^2}{2} \sin 2 \theta(t) \right\|_{H^{-1}} 
\leqslant \left\| \Delta \theta(t) + \frac{\mathrm{H}^2}{2} \sin 2 \theta(t) \right\|.
\end{aligned}
$$
Here, $\eta$ and $\nu$ are positive constants. Since $\nu \in \left(0, \frac{1}{2}\right)$, for all $t \in [t_n, \bar{t}_n]$, we can derive
$$
\begin{aligned}
\eta \left( \mathcal{E}(t) - E_1(\theta_\infty) \right)^{1 - \nu} &\leqslant \eta \left( \frac{1}{2} \|\mathbf{v}(t)\|^2 + \left| E_1(\theta(t)) - E_1(\theta_\infty) \right| \right)^{1 - \nu} \\
&\lesssim \left( \|\mathbf{v}(t)\|^2 + \left\| \Delta \theta(t) + \frac{\mathrm{H}^2}{2} \sin 2 \theta(t) \right\|^{\frac{1}{1 - \nu}} \right)^{1 - \nu}.
\end{aligned}
$$
Since $\|\mathbf{v}(t)\|$ is uniformly bounded, we have $\|\mathbf{v}(t)\|^{2(1 - \nu)} \leqslant C \|\mathbf{v}(t)\|$. By the Poincaré's inequality, we obtain
$$
\left( \mathcal{E}(t) - E_1(\theta_\infty) \right)^{1 - \nu} \lesssim \|\nabla \mathbf{v}(t)\| + \left\| \Delta \theta(t) + \frac{\mathrm{H}^2}{2} \sin 2 \theta(t) \right\|.
$$
From the basic energy equality, for all $t \in [t_n, \bar{t}_n]$, we have
$$
\begin{aligned}
\frac{\mathrm{d}}{\mathrm{d} t} \left( \mathcal{E}(t) - E_1(\theta_\infty) \right)^\nu &= \nu \left( \mathcal{E}(t) - E_1(\theta_\infty) \right)^{-(1 - \nu)} \frac{\mathrm{d}}{\mathrm{d} t} \mathcal{E}(t) \\
&\leqslant -\nu \left( \|\nabla \mathbf{v}(t)\| + \left\| \Delta \theta(t) + \frac{\mathrm{H}^2}{2} \sin 2 \theta(t) \right\| \right)^{-1} C_1 E_H(t) \\
&\leqslant -C_2 \left( \|\nabla \mathbf{v}(t)\| + \left\| \Delta \theta(t) + \frac{\mathrm{H}^2}{2} \sin 2 \theta(t) \right\| \right),
\end{aligned}
$$
where $C_1$ and $C_2$ are positive constants. Therefore, for all $t \in [t_n, \bar{t}_n]$, we obtain
$$
\begin{aligned}
C_2 \int_{t_n}^t \|\nabla \mathbf{v}(s)\| + \left\| \Delta \theta(s) + \frac{\mathrm{H}^2}{2} \sin 2 \theta(s) \right\| \, \mathrm{d}s 
\leqslant \left( \mathcal{E}(t_n) - E_1(\theta_\infty) \right)^\nu - \left( \mathcal{E}(t) - E_1(\theta_\infty) \right)^\nu.
\end{aligned}
$$
Thus,
\begin{equation}\label{iso_ine_thm_second}
\int_{t_n}^t \|\nabla \mathbf{v}(s)\| + \left\| \Delta \theta(s) + \frac{\mathrm{H}^2}{2} \sin 2 \theta(s) \right\| \, \mathrm{d}s \leqslant \frac{1}{C_2} \left( \mathcal{E}(t_n) - E_1(\theta_\infty) \right)^\nu.
\end{equation}
From (\ref{iso_ine_thm_second}) and the \L ojasiewicz--Simon inequality, for all $n \geqslant N$, we derive
$$
\begin{aligned}
\left\| \theta(\bar{t}_n) - \theta_{\infty}(n) \right\| &\leqslant \left\| \theta(t_n) - \theta_{\infty}(n) \right\| + \int_{t_n}^{\bar{t}_n} \left\| \partial_\tau \theta \right\| \, \mathrm{d}\tau \\
&\leqslant \left\| \theta(t_n) - \theta_{\infty}(n) \right\| + C \left( \mathcal{E}(t_n) - E_1(\theta_\infty) \right)^\nu.
\end{aligned}
$$
As $n \rightarrow \infty$, we have
$$
\left( \left\| \theta(t_n) - \theta_{\infty}(n) \right\| + C \left( \mathcal{E}(t_n) - E_1(\theta_\infty) \right)^\nu \right) \rightarrow 0.
$$
Therefore,
$$
\lim_{n \rightarrow \infty} \left\| \theta(\bar{t}_n) - \theta_{\infty}(n) \right\| = 0.
$$

Since $\|\nabla \theta(t)\|_{H^1_p} = \|\nabla \mathbf{d}(t)\|_{H^1_p} \leqslant C$, the sequence $\nabla \theta(t)$ is relatively compact in $L^2_p$. Thus, there exists a subsequence $\{\nabla \theta(\bar{t}_n)\}$, still denoted as $\{\nabla \theta(\bar{t}_n)\}$, which converges in $L^2_p$ to $\nabla \theta_{\infty}(n) = \nabla \theta_{\infty}$. Hence, we obtain
\begin{equation}\label{iso_thm3.10_third}
\lim_{n \rightarrow \infty} \left\| \theta(\bar{t}_n) - \theta_{\infty}(n) \right\|_{H^1} = 0.
\end{equation}

Since $\bar{t}_n$ satisfies
$$
\left\| \theta(t_n) - \theta_{\infty}(n) \right\|_{H^1_p} = \frac{\rho}{2},
$$
by (\ref{iso_thm3.10_third}), there exists a sufficiently large positive integer $n$ such that
$$
\left\| \theta(\bar{t}_n) - \theta_{\infty}(n) \right\|_{H^1} < \frac{\rho}{4}.
$$
This contradicts the definition of $\bar{t}_n$, hence there exists an integer constant $N_1$ such that $\bar{t}_{N_1} = \infty$.
Taking $n=N_1$, $\theta_\infty = \theta_{\infty}(N_1)$
 and $T_1 = t_{N_1}$, the inequality (\ref{iso_inequality_boundness_first_1_theta}) is proved.

Additionally, $|\theta(x, t)|$ is bounded in $Q \times (0, \infty)$.

\paragraph{}
(iii)
Directly compute \( \left|\mathbf{d} - \mathbf{d}_\infty\right|^2 \),  
$$
\left|\mathbf{d} - \mathbf{d}_\infty\right|^2 = \left|\sin\theta - \sin\theta_\infty\right|^2 + \left|\cos\theta - \cos\theta_\infty\right|^2 \leqslant 2 \left|\theta - \theta_\infty \right|^2,
$$
from which we can derive  
\begin{equation}\label{iso_the3.10_ine_forth}
\left\|\mathbf{d}(t) - \mathbf{d}_\infty\right\| \leqslant C \left\|\theta(t) - \theta_\infty\right\|.
\end{equation}

Since the function \( \theta_\infty \) is smooth, we have \( \left|\nabla\theta_\infty\right| \leqslant C \). For \( k = 1, 2 \), we can obtain  
$$
\begin{aligned}
\left|\partial_k\left( \sin\theta - \sin\theta_\infty \right)\right| &\leqslant \left|\cos\theta \left( \partial_k \theta - \partial_k \theta_\infty \right)\right| + \left|\left(\cos\theta - \cos\theta_\infty\right) \partial_k \theta_\infty \right| \\
&\leqslant C \left|\theta - \theta_\infty\right| + C \left|\partial_k \theta - \partial_k \theta_\infty \right|.
\end{aligned}
$$  
Similarly, we can derive  
$$
\left|\partial_k\left( \cos\theta - \cos\theta_\infty \right)\right| \leqslant C \left|\theta - \theta_\infty - 2 l_n \pi\right| + C \left|\partial_k \theta - \partial_k \theta_\infty \right|.
$$  
Thus,  
$$
\left|\nabla(\mathbf{d}(t) - \mathbf{d}_\infty)\right| \leqslant C \left|\nabla(\theta(t) - \theta_\infty)\right|.
$$  
Combining this with (\ref{iso_the3.10_ine_forth}), we obtain  
$$
\left\|\mathbf{d}(t) - \mathbf{d}_\infty\right\|_{H^1} \leqslant C \left\|\theta(t) - \theta_\infty \right\|_{H^1}.
$$

Therefore,
there exists $T_2$ such that for $t > T_2$,
$$
\left\| \mathbf{d} - \mathbf{d}_{\infty} \right\|_{H^1} \leqslant C \left\| \theta - \theta_{\infty} \right\|_{H^1} < \rho.
$$
Thus, we can obtain (\ref{iso_inequality_boundness_first_2_d}). $\Box$

\begin{remark}
(i)
Directly compute \( \nabla^2 (\sin\theta - \sin\theta_{\infty}) \),  
$$
\begin{aligned}
&\partial_i\partial_k (\sin\theta - \sin\theta_{\infty}) \\
=& -\sin\theta (\partial_i\theta \partial_k\theta) + \sin\theta_\infty (\partial_i\theta_\infty \partial_k\theta_\infty) + \cos\theta (\partial_i\partial_k\theta) - \cos\theta_\infty (\partial_i\partial_k\theta_\infty) \\
=& -\sin\theta \left[ (\partial_i\theta \partial_k\theta) - (\partial_i\theta_\infty \partial_k\theta_\infty) \right] - (\sin\theta - \sin\theta_\infty)(\partial_i\theta_\infty \partial_k\theta_\infty) \\
&+ \cos\theta (\partial_i\partial_k\theta - \partial_i\partial_k\theta_\infty) + \partial_i\partial_k\theta_\infty (\cos\theta - \cos\theta_\infty).
\end{aligned}
$$
From this, we can derive  
$$
\left\| \nabla^2 \left( \mathbf{d}(t) - \mathbf{d}_{\infty} \right) \right\| \leqslant C \left\| \theta(t) - \theta_{\infty}  \right\|_{H^2}.
$$

(ii)
By directly computing \( \left| \nabla (\mathbf{d} - \mathbf{d}_{\infty}) \right|^2 \) and \( \left| \nabla (\theta - \theta_{\infty}) \right|^2 \), we obtain  
$$
\begin{aligned}
\left| \nabla (\mathbf{d} - \mathbf{d}_{\infty}) \right|^2 =& \left( \partial_{1}\theta \right)^2 + \left( \partial_{1}\theta_\infty \right)^2 + \left( \partial_{2}\theta \right)^2 + \left( \partial_{2}\theta_\infty \right)^2 \\
&- 2 \cos(\theta - \theta_{\infty}) \left( \partial_{1}\theta \partial_{1}\theta_\infty + \partial_{2}\theta \partial_{2}\theta_\infty \right), \\
\left| \nabla (\theta - \theta_{\infty}) \right|^2 =& \left( \partial_{1}\theta \right)^2 + \left( \partial_{1}\theta_\infty \right)^2 + \left( \partial_{2}\theta \right)^2 + \left( \partial_{2}\theta_\infty \right)^2 
- 2 \left( \partial_{1}\theta \partial_{1}\theta_\infty + \partial_{2}\theta \partial_{2}\theta_\infty \right).
\end{aligned}
$$
Thus,  
$$
\begin{aligned}
&2 \left| \nabla (\mathbf{d} - \mathbf{d}_{\infty}) \right|^2 - \left| \nabla (\theta - \theta_{\infty}) \right|^2 \\
=& \left( \partial_{1}\theta \right)^2 + \left( \partial_{1}\theta_\infty \right)^2 + \left( \partial_{2}\theta \right)^2 + \left( \partial_{2}\theta_\infty \right)^2 
- \left[ 4 \cos(\theta - \theta_{\infty}) - 2 \right] \left( \partial_{1}\theta \partial_{1}\theta_\infty + \partial_{2}\theta \partial_{2}\theta_\infty \right) \\
\geqslant& \left( 1 - \left| 2 \cos(\theta - \theta_{\infty}) - 1 \right| \right) \left[ \left( \partial_{1}\theta \right)^2 + \left( \partial_{1}\theta_\infty \right)^2 + \left( \partial_{2}\theta \right)^2 + \left( \partial_{2}\theta_\infty \right)^2 \right].
\end{aligned}
$$
where 
$$
1 - \left| 2 \cos(\theta - \theta_{\infty}) - 1 \right| = 2 - 2 \cos(\theta - \theta_{\infty}) \geqslant 0.
$$
Therefore,
\begin{equation}\label{important_1}
2 \left| \nabla (\mathbf{d} - \mathbf{d}_{\infty}) \right|^2 
\geqslant
\left| \nabla (\theta - \theta_{\infty}) \right|^2.    
\end{equation}
\end{remark}

\subsection{Convergence to an equilibrium}
If there exists a constant \( T \geqslant 0 \) such that for all \( t \geqslant T \),  we have
\(\mathcal{E}(t) = E_1(\theta_\infty)\),
then the solution \( (\mathbf{v}, \mathbf{d}) \) of system (\ref{iso_Ericksen_Leslie_new_formulation}) reaches the steady state \( (0, \mathbf{d}_\infty) \) in finite time.  
This case is trivial,
and we just need to discuss the non-trivial case  for all \( t \geqslant 0 \),  we have
\(\mathcal{E}(t) > E_1(\theta_\infty)\).

\begin{lemma}\label{iso_lem3.11}
Suppose that \( A(t) \) and \( B(t) \) are continuously non-negative differentiable functions. 
Let \( a \) and \( b \) be positive constants with \( a < b \). The function \( A(t) \) satisfies
$$
\frac{\mathrm{d}}{\mathrm{d}t} A(t) + a A(t) \leqslant b A(t).
$$

(i) If \( \nu \in \left(0, \frac{1}{2}\right) \), and for \( t > 0 \),
$$
\frac{\mathrm{d}}{\mathrm{d}t} B(t) + a B(t) + a A(t) \lesssim (1 + t)^{-\frac{2\nu}{1 - 2\nu}},
$$
then
\begin{equation}\label{iso_lem3.7_inequality_1}
A(t) + B(t) \lesssim (1 + t)^{-\frac{2\nu}{1 - 2\nu}}.
\end{equation}

(ii) If \( \nu = \frac{1}{2} \), and there exists a constant \( c > 0 \) such that for all \( t > 0 \),
$$
\frac{\mathrm{d}}{\mathrm{d}t} B(t) + a B(t) + a A(t) \lesssim e^{-c t},
$$
then there exists a constant \( b_1 > 0 \) such that
\begin{equation}\label{iso_lem3.7_inequality_2}
A(t) + B(t) \lesssim e^{-b_1 t}.
\end{equation}
\end{lemma}

\paragraph{Proof.}
(i) 
For \( \nu \in \left(0, \frac{1}{2}\right) \), 
we have
$$
\frac{\mathrm{d}}{\mathrm{d}t} B(t) + a B(t) \lesssim (1 + t)^{-\frac{2\nu}{1 - 2\nu}}.
$$
Thus,
$$
\frac{\mathrm{d}}{\mathrm{d} t}\left[e^{a t} B(t)\right] \leqslant C e^{a t}(1 + t)^{-\frac{2 \nu}{1 - 2 \nu}}.
$$
Integrating both sides, we get
$$
\begin{aligned}
B(t) & \leqslant e^{-a t} B(0) + C \int_0^t e^{-a(t - \tau)}(1 + \tau)^{-\frac{2 \nu}{1 - 2 \nu}} \mathrm{d} \tau \\
& \lesssim e^{-a t} + (1 + t)^{-\frac{2 \nu}{1 - 2 \nu}} \int_0^t e^{-a(t - \tau)}\left(\frac{1 + t}{1 + \tau}\right)^{\frac{2 \nu}{1 - 2 \nu}} \mathrm{d} \tau \\
& \lesssim e^{-a t} + (1 + t)^{-\frac{2 \nu}{1 - 2 \nu}} \int_0^t e^{-a(t - \tau)}(t - \tau + 1)^{\frac{2 \nu}{1 - 2 \nu}} \mathrm{d} \tau.
\end{aligned}
$$
For \( \nu \in [0, t] \), we have
$$(1 + \tau)(t - \tau + 1) \geqslant 1 + t.$$ 
Since \( a > 0 \), for all \( t > 0 \), we obtain
$$
B(t) \lesssim (1 + t)^{-\frac{2 \nu}{1 - 2 \nu}}.
$$

Next, let \( a_0 = \frac{a}{b} \), 
and we derive
$$
\frac{\mathrm{d}}{\mathrm{d}t} \left(B(t) + a_0 A(t)\right) + a \left(B(t) + a_0 A(t)\right) \lesssim (1 + t)^{-\frac{2\nu}{1 - 2 \nu}}.
$$
Similarly,
$$
B(t) + a_0 A(t) \lesssim (1 + t)^{-\frac{2 \nu}{1 - 2 \nu}}.
$$
Thus, inequality (\ref{iso_lem3.7_inequality_1}) is established.

\paragraph{}
(ii) Let \( b_1 = \frac{1}{2} \min(a, c) \). For \( \nu = \frac{1}{2} \), we have
$$
\frac{\mathrm{d}}{\mathrm{d}t} B(t) + b_1 B(t) \lesssim e^{-c t}.
$$
Thus,
$$
\frac{\mathrm{d}}{\mathrm{d} t}\left[e^{b_1 t} B(t)\right] \leqslant C e^{-(c - b_1) t}.
$$
Integrating both sides, we get
$$
\begin{aligned}
B(t) & \leqslant e^{-b_1 t} B(0) + C e^{-b_1 t} \int_0^t e^{-(c - b_1) \tau} \mathrm{d} \tau \\
& \lesssim e^{-b_1 t}.
\end{aligned}
$$

Next, since \( a_0 = \frac{a}{b} \), 
 we derive
$$
\frac{\mathrm{d}}{\mathrm{d}t} \left(B(t) + a_0 A(t)\right) + a \left(B(t) + a_0 A(t)\right) \lesssim e^{-c t}.
$$
Similarly,
$$
B(t) + a_0 A(t) \lesssim e^{-b_1 t}.
$$
Thus, inequality (\ref{iso_lem3.7_inequality_2}) is established.

\begin{theorem}\label{iso_final_theorem}
Suppose that \((\mathbf{v}_0, \mathbf{d}_0) \in V \times H^2_{sp}\),
then the following inequalities hold.

(i) If \(\nu \in \left(0, \frac{1}{2}\right)\), then for all \(t > 0\),
\begin{equation}\label{iso_final_theorem_inequality_1}
\|\mathbf{v}(t)\|_{H^1} + \left\|\theta(t) - \theta_{\infty}\right\|_{H^2} \lesssim (1 + t)^{-\frac{\nu}{1 - 2\nu}}.
\end{equation}

(ii) If \(\nu = \frac{1}{2}\), then for all \(t > 0\), there exists a constant \(c > 0\) such that
\begin{equation}\label{iso_final_theorem_inequality_2}
\|\mathbf{v}(t)\|_{H^1} + \left\|\theta(t) - \theta_{\infty}\right\|_{H^2} \lesssim e^{-ct}.
\end{equation}
The constant \(c\) depends on $Q$,  $ |\mathbf{H}|$, $\left\|\mathbf{d}_0\right\|_{H^2}$,  
$\left\|\mathbf{v}_0\right\|_{H^1}$
and
coefficients of the system (\ref{iso_Ericksen_Leslie_new_formulation}).
\end{theorem}

\paragraph{Proof.}
(i) First, we prove the convergence to an equilibrium regarding $\theta$ in $L^2$ norm.

By the basic energy equality, for any \(t > 0\), we have
\begin{equation}\label{iso_thm_final_proof_inequ_1}
\begin{aligned}
0 \geqslant \frac{\mathrm{d}}{\mathrm{d}t}\left(\mathcal{E}(t) - E_1(\theta_\infty)\right) + C E_H(t) 
\geqslant \frac{\mathrm{d}}{\mathrm{d}t}\left(\mathcal{E}(t) - E_1(\theta_\infty)\right) + C\left(\mathcal{E}(t) - E_1(\theta_\infty)\right)^{2(1 - \nu)}.
\end{aligned}
\end{equation}

Let \(T_1\) be the constant from Theorem \ref{iso_thm_boundness_theta}. 
If \(\nu \in \left(0, \frac{1}{2}\right)\), then for all \(t > T_1\),
$$
\begin{aligned}
\left\|\theta(t) - \theta_{\infty}\right\| &\leqslant \int_t^{\infty} \left\|\partial_{\tau} \theta\right\| \, \mathrm{d}\tau \\
&\lesssim \int_t^{\infty} \left(\|\nabla \mathbf{v}(s)\| + \left\|\Delta \theta(s) + \frac{\mathrm{H}^2}{2} \sin 2\theta(s)\right\|\right) \, \mathrm{d}s \\
&\lesssim \left(\mathcal{E}(t) - E_1(\theta_\infty)\right)^\nu \\
&\lesssim (1 + t)^{-\frac{\nu}{1 - 2\nu}}.
\end{aligned}
$$
Therefore, we can obtain
\begin{equation}\label{iso_Thm3201}
\left\|\theta(t) - \theta_{\infty}\right\| \lesssim (1 + t)^{-\frac{\nu}{1 - 2\nu}}.
\end{equation}

If \(\nu = \frac{1}{2}\), by (\ref{iso_thm_final_proof_inequ_1}),
we derive for all \(t > 0\),
$$
\frac{\mathrm{d}}{\mathrm{d}t}\left(\mathcal{E}(t) - E_1(\theta_\infty)\right) + C\left(\mathcal{E}(t) - E_1(\theta_\infty)\right) \leqslant 0.
$$
Thus, 
$$
\mathcal{E}(t) - E_1(\theta_\infty) \lesssim e^{-Ct}.
$$
Therefore, 
for all \(t > T_1\),
\begin{equation}\label{iso_thm_final_proof_inequ_2}
\left\|\theta(t) - \theta_{\infty}\right\| \lesssim \left(\mathcal{E}(t) - E_1(\theta_\infty)\right)^\nu \lesssim e^{-\nu Ct}.
\end{equation}

One the other hand,
since $\theta(x,t)$ is bounded in $Q\times \left(0,T_1\right]$ and $\theta_\infty$ is  bounded in $Q$,
then for $t\in \left(0,T_1\right]$,
$$
\left\|\theta(t) - \theta_{\infty}\right\| 
\leqslant
C
\lesssim 
\min\{(1 + T_1)^{-\frac{\nu}{1 - 2\nu}},
e^{-\nu CT_1}\}
\lesssim 
\min\{(1 + t)^{-\frac{\nu}{1 - 2\nu}},
e^{-\nu Ct}\}.
$$
Therefore, inequalities (\ref{iso_Thm3201}) and (\ref{iso_thm_final_proof_inequ_2}) hold for $t>0$.

\paragraph{}
(ii)
From  (\ref{iso_Ericksen_Leslie_new_formulation})  and (\ref{elliptic_first}), we derive
\begin{equation}\label{iso_Ericksen_Leslie_new_formulation_elliptic_first_1}
\left\{
\begin{aligned}
\partial_t \mathbf{v} =& -\mathbf{v} \cdot \nabla \mathbf{v} - \nabla \widetilde{p} + \frac{\gamma}{\operatorname{Re}} \Delta \mathbf{v} + \frac{1 - \gamma}{\operatorname{Re}} \nabla \cdot \sigma_1  - \frac{1 - \gamma}{\operatorname{Re}} \left( \Delta \theta \cdot \nabla \theta \right) \\
&+ \frac{1 - \gamma}{2 \operatorname{Re}} \nabla \cdot \left[ \left( (1 - \mu_2) \mathbf{d} \otimes \mathbf{d}^{\perp} - (1 + \mu_2) \mathbf{d}^{\perp} \otimes \mathbf{d} \right) (\Delta \theta + \frac{\mathrm{H}^2}{2} \sin 2\theta) \right], \\
\nabla \cdot \mathbf{v} =&\; 0, \\
\partial_t \theta =& -\mathbf{v} \cdot \nabla \theta + \mu_1 \left( \Delta (\theta - \theta_{\infty}) + \frac{\mathrm{H}^2}{2} \sin 2\theta - \frac{\mathrm{H}^2}{2} \sin 2\theta_{\infty} \right) \\
&\quad + \left( \boldsymbol{\Omega} : \mathbf{d} \otimes \mathbf{d}^{\perp} + \mu_2 \mathbf{D} : \mathbf{d} \otimes \mathbf{d}^{\perp} \right),
\end{aligned}
\right.
\end{equation}
where
$$
\widetilde{p} = p + \frac{1 - \gamma}{2 \operatorname{Re}} |\nabla \theta|^2,
$$
and \(\theta - \theta_{\infty}\) is periodic in \(x\). 

Taking the inner product of the first  equation in (\ref{iso_Ericksen_Leslie_new_formulation_elliptic_first_1}) with \(\mathbf{v}\) and 
taking the inner product of the first  equation in (\ref{iso_Ericksen_Leslie_new_formulation_elliptic_first_1}) with\(\frac{1 - \gamma}{\operatorname{Re}} \left( \Delta (\theta - \theta_{\infty}) + \frac{\mathrm{H}^2}{2} \sin 2\theta - \frac{\mathrm{H}^2}{2} \sin 2\theta_{\infty} \right)\).
We can obtain
\begin{equation}\label{iso_Ericksen_Leslie_new_formulation_elliptic_first_2}
\begin{aligned}
& \frac{\mathrm{d}}{\mathrm{d} t} \int_Q |\mathbf{v}|^2 + \frac{1 - \gamma}{\operatorname{Re}} \left| \nabla (\theta - \theta_{\infty}) \right|^2 \, \mathrm{d}x \\
&+ \frac{\mathrm{d}}{\mathrm{d} t} \int_Q \frac{(1 - \gamma) \mathrm{H}^2}{2 \operatorname{Re}} \left[ \cos 2\theta - \cos 2\theta_{\infty} + 2 \sin 2\theta_{\infty} (\theta - \theta_{\infty}) \right] \, \mathrm{d}x \\
=& -2 \int_Q \frac{\gamma}{\operatorname{Re}} |\nabla \mathbf{v}|^2 + \frac{(1 - \gamma)\mu_1}{\operatorname{Re}} \left[ |\mathbf{h}|^2 - \mu_1 (\mathbf{h} \cdot \mathbf{d})^2 \right] \, \mathrm{d}x \\
&- 2 \int_Q \frac{1 - \gamma}{\operatorname{Re}} \left[ \beta_1 (\mathbf{d} \otimes \mathbf{d} : \mathbf{D})^2 + \beta_2 \mathbf{D} : \mathbf{D} + \beta_3 |\mathbf{D} \cdot \mathbf{d}|^2 \right] \, \mathrm{d}x.
\end{aligned}
\end{equation}

From the third equation in (\ref{iso_Ericksen_Leslie_new_formulation_elliptic_first_1}), we have
\begin{equation}\label{iso_Ericksen_Leslie_new_formulation_elliptic_first_3}
\begin{aligned}
& \frac{1}{2} \frac{\mathrm{d}}{\mathrm{d} t} \int_Q \left| \theta - \theta_{\infty} \right|^2 \, \mathrm{d}x + \mu_1 \int_Q \left| \nabla (\theta - \theta_{\infty}) \right|^2 \, \mathrm{d}x \\
=& - \int_Q \mathbf{v} \cdot \nabla \theta (\theta - \theta_{\infty}) \, \mathrm{d}x + \frac{\mathrm{H}^2 \mu_1}{2} \int_Q (\sin 2\theta - \sin 2\theta_{\infty}) (\theta - \theta_{\infty}) \, \mathrm{d}x \\
&+ \int_Q \left( \boldsymbol{\Omega} : \mathbf{d} \otimes \mathbf{d}^{\perp} + \mu_2 \mathbf{D} : \mathbf{d} \otimes \mathbf{d}^{\perp} \right) (\theta - \theta_{\infty}) \, \mathrm{d}x.
\end{aligned}
\end{equation}
Thus, for a sufficiently small constant \(\epsilon\),
$$
\frac{1}{2} \frac{\mathrm{d}}{\mathrm{d} t} \int_Q \left| \theta - \theta_{\infty} \right|^2 \, \mathrm{d}x + \int_Q \left| \nabla (\theta - \theta_{\infty}) \right|^2 \, \mathrm{d}x \leqslant \epsilon \|\nabla \mathbf{v}\|^2 + C \left\| \theta - \theta_{\infty} \right\|^2.
$$

Define
$$
\begin{aligned}
G(t) :=& \frac{1}{2} \int_Q |\mathbf{v}|^2 + \left| \theta - \theta_{\infty} \right|^2 + \frac{1 - \gamma}{\operatorname{Re}} \left| \nabla (\theta - \theta_{\infty}) \right|^2 \, \mathrm{d}x \\
&+ \frac{1}{2} \int_Q \frac{1 - \gamma}{\operatorname{Re}} \left( \frac{\mathrm{H}^2}{2} \left[ \cos 2\theta - \cos 2\theta_{\infty} + 2 \sin 2\theta_{\infty} (\theta - \theta_{\infty}) \right] \right) \, \mathrm{d}x.
\end{aligned}
$$
It is straightforward to show that
$$
\left| \int_{\Omega} \left[ \cos 2\theta - \cos 2\theta_{\infty} + 2 \sin 2\theta_{\infty} (\theta - \theta_{\infty}) \right] \right| \leqslant 2 \left\| \theta - \theta_{\infty} \right\|.
$$

\paragraph{}
(iii) If \(\nu \in \left(0, \frac{1}{2}\right)\), from (\ref{iso_Ericksen_Leslie_new_formulation_elliptic_first_2}) and (\ref{iso_Ericksen_Leslie_new_formulation_elliptic_first_3}), there exists a constant \(a > 0\) such that for all \(t > 0\),
$$
\frac{\mathrm{d}}{\mathrm{d} t} G(t) + a G(t) + a E_H(t) \leqslant C \left\| \theta - \theta_{\infty} \right\|^2 \leqslant C (1 + t)^{-\frac{2\nu}{1 - 2\nu}}.
$$
By the higher-order estimate (\ref{iso_higher-order_energy_estimate}), there exists a constant \(b > 0\) such that for all \(t > 0\),
$$
\frac{\mathrm{d}}{\mathrm{d} t} E_H(t) + a E_H(t) \leqslant b E_H(t).
$$
Applying Lemma \ref{iso_lem3.11}, we obtain
$$
G(t) + E_H(t) \lesssim (1 + t)^{-\frac{2\nu}{1 - 2\nu}}.
$$
Thus, for all \(t > 0\),
$$
\begin{aligned}
&\|\mathbf{v}(t)\|^2 + \left\| \theta(t) - \theta_{\infty} \right\|_{H^1}^2 \leqslant 2 G(t) + C \left\| \theta - \theta_{\infty} \right\|^2 \lesssim (1 + t)^{-\frac{2\nu}{1 - 2\nu}}, \\
&\|\nabla \mathbf{v}(t)\|^2 + \left\| \Delta \theta - \Delta \theta_{\infty} \right\|^2 \lesssim E_H(t) + \left\| \theta(t) - \theta_{\infty} \right\|^2 \lesssim (1 + t)^{-\frac{2\nu}{1 - 2\nu}}.
\end{aligned}
$$
Inequality (\ref{iso_final_theorem_inequality_1}) is established.

If \(\nu = \frac{1}{2}\), from (\ref{iso_Ericksen_Leslie_new_formulation_elliptic_first_2}) and (\ref{iso_Ericksen_Leslie_new_formulation_elliptic_first_3}), there exist constants \(a, C > 0\) such that for all \(t > 0\),
$$
\frac{\mathrm{d}}{\mathrm{d} t} G(t) + a G(t) + a E_H(t) \leqslant C \left\| \theta - \theta_{\infty} \right\|^2 \leqslant C e^{-ct}.
$$
By the higher-order estimate (\ref{iso_higher-order_energy_estimate}), there exists a constant \(b > 0\) such that for all \(t > 0\),
$$
\frac{\mathrm{d}}{\mathrm{d} t} E_H(t) + a E_H(t) \leqslant b E_H(t).
$$
Applying Lemma \ref{iso_lem3.11}, there exists a constant \(c_1 > 0\) such that
$$
G(t) + E_H(t) \lesssim e^{-c_1 t}.
$$
Thus, for all \(t > 0\),
$$
\begin{aligned}
&\|\mathbf{v}(t)\|^2 + \left\| \theta(t) - \theta_{\infty} \right\|_{H^1}^2 \leqslant 2 G(t) + C \left\| \theta - \theta_{\infty} \right\|^2 \lesssim e^{-c_1 t}, \\
&\|\nabla \mathbf{v}(t)\|^2 + \left\| \Delta \theta - \Delta \theta_{\infty} \right\|^2 \lesssim E_H(t) + \left\| \theta(t) - \theta_{\infty} \right\|^2 \lesssim e^{-c_1 t}.
\end{aligned}
$$
Inequality (\ref{iso_final_theorem_inequality_2}) is established.
$\Box$

\paragraph{}
Define
\begin{equation}\label{L1}
L_{\theta_{\infty}} := -\Delta - \mathrm{H}^2 \cos \left(2 \theta_{\infty}\right).
\end{equation}
From Theorem 2.1 in \cite{Haraux2003}, if the kernel of \( L_{\theta_{\infty}} \) is trivial, meaning the smallest eigenvalue of \( L_{\theta_{\infty}} \) is greater than 0, then \( \nu = \frac{1}{2} \). This leads to the following proposition, inspired by Corollary 2.5 in \cite{CKY2018SIAM}.

\begin{proposition}
If \( \mathrm{H}^2 < \lambda_2 \), then \( \nu = \frac{1}{2} \).    
\end{proposition}

\paragraph{Proof.}
Let \( \theta_{\infty} \) be a solution to (\ref{elliptic_first}). The operator \( L_{\theta_{\infty}} = -\Delta - \mathrm{H}^2 \cos \left(2 \theta_{\infty}\right) \), defined by (\ref{L1}), is the linearized operator around \( \theta_{\infty} \).

If \( \mathrm{H}^2 < \lambda_2 \), then the smallest eigenvalue \( b \) of \( L_{\theta_{\infty}} \) is \( \lambda_2 - \mathrm{H}^2 > 0 \). Hence, the kernel of \( L_{\theta_{\infty}} \) is trivial, implying \( \nu = \frac{1}{2} \).
$\Box$

\paragraph{}
In summary, Theorem \ref{Theorem1.3} is proved.

\section{Appendix}\label{Appendix}

The parametric equations of the torus $\mathbb{T}^2$ embedded in \( \mathbb{R}^3 \) are:
\[
\left\{
\begin{aligned}
x =&\; (R + r \cos v) \cos u, \\
y =&\;  (R + r \cos v) \sin u, \\
z =&\;  r \sin v,
\end{aligned}
\right.
\]
where constant \( R \) is the major radius (distance from the center of the torus to the center of the tube), and
constant \( r \) is the minor radius (radius of the tube), with \( R >> r>0 \).
Angle \( u \in [0, 2\pi) \) is  around the major circle, and 
angle \( v \in [0, 2\pi) \) is  around the minor circle.

Consider a circular loop of radius \( R \) (located in the \( xy \)-plane with its center at the origin), carrying a constant current \( I \). According to {Ampère's Law} and the {Biot--Savart Law}, the loop induces a constant magnetic field on the torus \( \mathbb{T}^2 \). 
The parametric equation of the loop is
\[
\mathbf{r}(\theta) = (R \cos \theta, R \sin \theta, 0), \quad \theta \in [0, 2\pi).
\]
The current density can be expressed using the Dirac delta function
\[
\mathbf{J}(\mathbf{r}') = I \int_0^{2\pi} \delta(\mathbf{r}' - \mathbf{r}(\theta)) \, \mathbf{t}(\theta) \, \mathrm{d}\theta,
\]
where \( \mathbf{t}(\theta) = (-\sin \theta, \cos \theta, 0) \) is the tangential unit vector.

By Biot--Savart Law,
the magnetic field \( \mathbf{H} \) at any point \( \mathbf{r} \) is given by
\[
\mathbf{H}(\mathbf{r}) = \frac{\mu_0}{4\pi} \int \frac{\mathbf{J}(\mathbf{r}') \times (\mathbf{r} - \mathbf{r}')}{|\mathbf{r} - \mathbf{r}'|^3} \, \mathrm{d}^3\mathbf{r}'.
\]
Substituting the current density, we obtain:
\[
\mathbf{H}(\mathbf{r}) = \frac{\mu_0 I}{4\pi} \int_0^{2\pi} \frac{\mathbf{t}(\theta) \times (\mathbf{r} - \mathbf{r}(\theta))}{|\mathbf{r} - \mathbf{r}(\theta)|^3} \, \mathrm{d}\theta.
\]

If the loop is placed along the primary circle direction, which means the \( u \)-direction of the torus \( \mathbb{T}^2 \), the magnetic field forms closed field lines along the poloidal direction (the \( v \)-direction). 
Within the torus (the tubular region), the magnetic field is approximately uniform, near the loop
  \[
  \mathrm{H} \approx \frac{\mu_0 I}{2\pi r} \quad,
  \]
which is along the \( v \)-direction (the tangential direction of the secondary circle).

Select a closed poloidal path \( C \) on the torus (fixed \( u \), \( v \in [0, 2\pi) \)), and apply Ampère's Law:
\[
\oint_C \mathbf{H} \cdot \mathrm{d}\mathbf{l} = \mu_0 I.
\]
Since \( C \) encloses the current \( I \), we  can obtain
\[
\mathbf{H} \cdot 2\pi r = \mu_0 I, \quad   \mathbf{H}= \frac{\mu_0 I}{2\pi r},
\]
which is consistent with the Biot--Savart result.

 A constant current \( I \) in the circular loop generates a  constant magnetic field \( \mathbf{H} \), whose direction is determined by the right-hand rule.
The geometry of the torus \( \mathbb{T}^2 \) ensures that the magnetic field lines are closed, forming a poloidal field.
Consider the liquid crystal flow on \( \mathbb{T}^2 \) with \( \mathbf{d}\in T(\mathbb{T}^2) \), which is the model studied in this paper.

\paragraph{Acknowledge}
The author thanks Professor Hao Wu 
for helpful suggestions.

\bibliography{ref}	
\bibliographystyle{plain}

\end{document}